\documentclass[12pt]{article} %
\usepackage{fleqn,amssymb,amsmath,amscd,epsfig,version,theorem}
\newcommand{\Section}[1]{\section{#1}\setcounter{equation}{0}}
\newtheorem{theorem}{Theorem} [section]
\newtheorem {lemma}[theorem]{Lemma}
\newtheorem {proposition}[theorem]{Proposition}
\newtheorem {corollary}[theorem]{Corollary}
\theorembodyfont{\normalfont}
\newtheorem {definition}[theorem]{Definition}

\newtheorem {example}[theorem]{Example}

\newtheorem {remark}[theorem]{Remark}
\newtheorem {remarks}[theorem]{Remarks}
\newcommand{\beq}{\begin{equation}}
\newcommand{\eeq}{\end{equation}}
\newcommand{\Leq}[1]{\label{#1}\end{equation}}
\newcommand{\beqn}{\begin{eqnarray}}
\newcommand{\eeqn}{\end{eqnarray}}
\newcommand{\beqno}{\begin{eqnarray*}}
\newcommand{\eeqno}{\end{eqnarray*}}
\newcommand{\bem}{\l(\! \begin{array}}
\newcommand{\eem}{\end{array}\!\ri)}
\newcommand{\bsm}{\left(\begin{smallmatrix}} 
\newcommand{\esm}{\end{smallmatrix}\right)}  
\renewcommand {\l}{\left}
\newcommand {\ri}{\right}
\newcommand {\q}{{\vec{q}}}
\newcommand {\p}{{\vec{p}}}
\newcommand {\Q}{{\vec{Q}}}
\renewcommand {\P}{{\vec{P}}}
\newcommand {\pq}{{(\p,\q)}}
\newcommand {\pqs}{(\p_{0},\q_{0})}
\newcommand {\s}{{\vec{s}}}
\newcommand {\w}{{\vec{w}\,}}
\newcommand {\y}{{\vec{y}\,}}
\newcommand {\z}{{\vec{z}\,}}

\newcommand {\qs}{\q-\s_{l}}

\newcommand {\Rmin}{R_{\rm min}}
\newcommand {\Rvir}{R_{\rm vir}}       
\newcommand {\amin}{\alpha_{\rm min}}  
\newcommand {\rmin}{r_{\rm min}}       
\newcommand {\dmin}{d_{\rm min}}       
\newcommand {\dmax}{d_{\rm max}}

\newcommand {\Vmax}{V_{\rm max}}       
\newcommand {\Zmax}{Z_{\rm max}}       
\newcommand {\thmin}{\Theta_{\rm min}} 
\newcommand {\htop}{{h_{\rm top}}}
\newcommand {\htopi}{{h_{\rm top}^\infty}}
\newcommand {\diam}{{\rm diam}}
\newcommand {\dist}{{\rm dist}}
\newcommand {\dimH}{{\rm dim}_{\cH}}
\newcommand {\dimB}{{\rm dim}_{\cB}}
\newcommand {\Mu}{{\bR^3_\q}}
\newcommand {\Muh}{\hat{M}}
\newcommand {\Mo}{{\bf M}}

\newcommand {\Pin}{P_{\infty}}
\newcommand {\Pinp}{P_{\infty,+}}
\newcommand {\Zi}{Z_{\infty}}
\newcommand {\Vsr}{V_{\rm sr}}
\newcommand {\Hh}{\hat{H}}
\newcommand {\Hi}{H_{\infty}}
\newcommand {\Hhi}{\hat{H}_{\infty}}
\newcommand {\Pt}{\Phi^{t}}
\newcommand {\PutE}{\Phi^{t}_{E}} 
\newcommand {\Pit}{\Phi_{\infty}^{t}}

\newcommand {\pptE}[1]{\mbox{\boldmath $\Phi$}^{#1}_{E}}
\newcommand {\Opm}{\Omega^{\pm}}
\newcommand {\Op}{\Omega^{+}}
\newcommand {\Om}{\Omega^{-}}
\newcommand {\B}[1]{\l|#1\ri|}
\newcommand {\vep}{\varepsilon}
\newcommand {\eps}{\epsilon}   
\newcommand {\vv}{\varphi}     

\newcommand {\LA}{\left\langle}
\newcommand {\RA}{\right\rangle}
\newcommand {\pa}{\partial}
\newcommand {\pih}{\hat{\pi}}
\newcommand {\eh}{{\textstyle \frac{1}{2}}}

\newcommand {\ev}{{\textstyle \frac{1}{4}}}
\newcommand {\ea}{{\textstyle \frac{1}{8}}}
\newcommand {\inInd}{ \in \{1,\ldots,n\} }
\newcommand {\ar}{\rightarrow}
\newcommand {\sign}{{\rm sign}}
\newcommand {\Id}{{\rm Id}}
\newcommand {\tr}{{\rm tr}}

\newcommand {\Ule}{U_{l}^{\vep}}
\newcommand {\Uhle}{{\hat{U}_{l}^{\vep}}}
\newcommand {\Sle}{S_{l}^{\vep}}
\newcommand {\Shle}{{\hat{S}_{l}^{\vep}}}
\newcommand {\SuE}{\Sigma_{E}}

\newcommand {\SoE}{\mbox{\boldmath $\Sigma$}_{E}}

\newcommand {\TuE}{T_{E}}

\newcommand {\buE}{b_{E}}
\newcommand {\Adm}{{\bf X}}
\newcommand {\Ho}{\mbox{\boldmath ${\cF}$} } 

\newcommand {\Eth}{{E_{\rm th}}}

\newcommand {\htheta}{{\hat{\theta}}}
\newcommand {\DCS}{\frac{d\sigma}{d\htheta^{+}}}
\newcommand {\DCSP}{\DCS(E,\htheta^{-},\htheta^{+})}
\newcommand {\vet}{\hat{\vv}_{E,\htheta^{-}}}

\newcommand {\set}{\sigma(E,\htheta^{-})}

\newcommand {\Liou}{\lambda}
\newcommand {\NuE}{N_{E}}
\newcommand {\bR}{{\mathbb R}}
\newcommand {\bN}{{\mathbb N}}
\newcommand {\bZ}{{\mathbb Z}}
\newcommand {\bC}{{\mathbb C}}

\newcommand {\bH}{{\mathbb H}}

\newcommand {\bQ}{{\mathbb Q}}
\newcommand {\ImH}{{{\rm Im}\bH}}
\newcommand {\hH}{{\hat{\bH}}}
\newcommand {\hpi}{{\hat{\pi}}}
\newcommand {\btheta}{\mbox{\boldmath $\theta$} }
\newcommand{\rstr}{{\upharpoonright}}
\newcommand{\idty}{{\rm 1\mskip-4mu l}} 
\newcommand{\vidty}{{\vec{\idty}}}      
\newcommand{\cB}{{\cal B}} 
\newcommand{\cC}{{\cal C}} 
\newcommand{\cD}{{\cal D}} 
\newcommand{\cE}{{\cal E}} %
\newcommand{\cF}{{\cal F}} 
\newcommand{\cH}{{\cal H}}
\newcommand{\cI}{{\cal I}} 
\newcommand{\cL}{{\cal L}} 
\newcommand{\cM}{{\cal M}}
\newcommand{\cO}{{\cal O}} 
\newcommand{\cP}{{\cal P}} 
\newcommand{\cQ}{{\cal Q}} 
\newcommand{\cS}{{\cal S}} 
\newcommand{\cW}{{\cal W}} 
\newcommand {\cDo}{\mbox{\boldmath $\cD$}}
\newcommand{\Mid}{\vec{m}}
\newcommand{\Ax}{A}
\newcommand{\Poi}{{\cH}} 
\newcommand {\Po}{{\cal P}_E} 
\newcommand{\Cyl}{{\cal Z}}
\newcommand{\ov}{\overline}
\newcommand{\NN}{\nonumber}
\newcommand {\Pol}{\mbox{\boldmath $\Phi$}_{l}}    
\newcommand {\PoLl}{\mbox{\boldmath $\Phi$}_{L,l}} 
\newcommand {\Psol}{{\mbox{\boldmath $\Psi$}_{l}}}    
\newcommand {\PsolN}{{\mbox{\boldmath $\Psi$}_{l}^0}}  
\newcommand {\PsoLl}{\mbox{\boldmath $\Psi$}_{L,l}} 
\newcommand {\Psul}{\Psi_{l}}                      
\newcommand {\PsuLl}{\Psi_{L,l}}                   
\newcommand {\Tol}{{\bf T}_l}      
\newcommand {\ToLl}{{\bf T}_{L,l}} 
\newcommand {\TuR}{T_E}  
\newcommand {\Dol}{{{\cDo}_{l}}}    
\newcommand {\DoLl}{{{\cDo}_{L,l}}} 
\newcommand {\DoLlb}{{\breve{{\cDo}}_{L,l}}} %
\newcommand {\Hol}{{{\mbox{\boldmath $\cH$}_{l}}}}    
\newcommand {\HoLl}{{{\mbox{\boldmath $\cH$}_{L,l}}}}  
\newcommand {\Hul}{{\cH}_{l}}        
\newcommand {\qmbox}[1]{\quad\mbox{#1}\quad}
\newcommand {\AD}{{\cal A\!D}} 
\newcommand {\CT}{{\cal C\!T}} 
\newcommand {\NCT}{{\cal N\!C\!T}} 
\newcommand {\IZ}{{\cal I\!Z}} 
\newcommand {\NC}{{\cal N\!C}} 
\newcommand {\TR}{{\cal T\!R}} 
\newcommand {\Xio}{\mbox{\boldmath $\Xi$}}
\newcommand {\Ind}{{\rm Index}} 
\newcommand {\uk}{{\underline k}} 
\newcommand {\ul}{{\underline l}} 
\newcommand {\um}{{\underline m}} 
\newcommand {\CC}{{\rm const.\,}} %
\begin{document}
\title {The $n$-Centre Problem of Celestial Mechanics for Large Energies}
\author{Andreas Knauf\thanks{Mathematisches Institut der
Universit\"at Erlangen-N\"urnberg. Bismarckstr. 1 1/2, 91054 Erlangen. 
e-mail: knauf@mi.uni-erlangen.de
}}
\date{August 2001}
\maketitle
\begin{abstract}
We consider the classical three-dimensional motion in a potential which is 
the sum of $n$ attracting or repelling Coulombic potentials.
Assuming a non-collinear configuration of the $n$ centres,
we find a universal behaviour for all energies $E$ above a positive threshold.

Whereas for $n=1$ there are no bounded orbits, and for $n=2$ there is
just one closed orbit,
for $n\geq 3$ the bounded orbits form a Cantor set.
We analyze the symbolic dynamics and estimate Hausdorff dimension and
topological entropy of this hyperbolic set.

Then we set up scattering theory, including symbolic dynamics of the scattering
orbits
and differential cross section estimates.

The theory includes the $n$--centre problem of celestial mechanics,
and prepares for a geometric understanding of a class of restricted 
$n$-{\em body} problems.

To allow for applications in semiclassical molecular scattering,
we include an additional smooth (electronic) potential 
which is arbitrary except its Coulombic decay at infinity.
Up to a (optimal) relative error of order $1/E$, all estimates are 
independent of that
potential but only depend on the relative positions
and strengths of the centres.

Finally we show that different, non-universal, phenomena occur for collinear
configurations. 
\end{abstract}
\newpage
\tableofcontents
%
\Section{Introduction} %
%
The {\em $n$-body problem} of celestial mechanics is the problem of solving
the Hamilton equation for the Hamiltonian function 
\[H:T^*M\ar\bR\qmbox{,}
H(\p_1,\ldots,\p_n,\q_1,\ldots,\q_n):=
\sum_{i=1}^n \frac{\p_i^{\,2}}{2m_i} +v\sum_{i\neq k} \frac{ Z_i Z_k}{|\q_i-\q_k|}\]
on the phase space $T^*M$ over the configuration space 
\[M:= \l\{(\q_1,\ldots,\q_n)\in(\bR^3)^n\mid 
\q_i\neq\q_k\mbox{ for } i\neq k\ri\}.\] 
In celestial mechanics $v=-1$ and
the $Z_i$ coincide with the positive masses $m_i>0$
(in units where the gravitational constant equals one).
However in an electrostatical context $v=+1$, 
the $Z_i$ are interpreted as charges and may be positive or negative.

Whereas the one-body problem corresponds to free motion and the two-body
problem was solved by Newton, the $n\geq3$-body problem is 
analytically non-integrable.

If one of $n+1$ bodies is much faster than the others then one may
approximate its motion by considering the {\em $n$--centre problem.} There 
the Hamiltonian function $\Hh:T^*\Muh\ar\bR$ on the phase space $T^*\Muh$ 
over the configuration space 
\[\Muh:=\bR^3\setminus \{\s_{1},\ldots,\s_{n}\}\] 
of that body is given by
\beq 
\Hh(\p,\q):=\eh\p^{\,2} +V(\q)\qmbox{with}
V(\q):=-\sum_{l=1}^n \frac{Z_l}{|\q-\s_l|}.
\Leq{ncp}
The 2--body problem thus reduces after separation of the centre of mass motion
to the 1--centre problem, that is the Kepler problem.

Moreover, the 2--centre problem is integrable and has been solved by Jacobi
(see Appendix B). This solution is particularly relevant 
for calculating the motion of artificial satellites, since the gravitational
field of the earth can be approximated by the one of two centres 
(analytically continued, since the earth is oblate and not prolate).
This has been used by Vinti in \cite{Vi}.  See
\cite{GKM} for an application. 

However, like the $n$--body problem, 
the $n$--centre problem for $n\geq3$ is analytically 
non-integrable, see Bolotin \cite{Bol}.

The analogy between the Coulombic and the gravitational interaction is
not perfect, since in {\em electrostatics} 
one has repulsive as well as attractive 
forces. This allows us to approximate molecules by static configurations of
nuclei surrounded by an electronic cloud of opposite charge.
Thus the electronic potential $V$ of the molecule has the form
\beq
V(\q)=  \sum_{l=1}^{n}\frac{-Z_{l}}{\B{\q-\s_{l}}}+ W(\q), 
\Leq{form:V}
$\s_l\in\bR^3$ being the position, and $-Z_l\neq0$ the charge of the $l$th
nucleus, multiplied by the test charge. 
The smooth electronic potential $W:\bR^3\ar\bR$ may partially shield
the nucleonic charge.
We model this by assuming the existence of
a net asymptotic charge $\Zi\in\bR$ with 
\beq
V(\q) = \frac{-\Zi}{\B{\q}} +\cO(|\q|^{-1-\eps})\qquad (\q\ar\infty).
\Leq{def:Zi}
Thus we consider the generalized $n$--centre problem (\ref{ncp}) with these
Coulombic potentials $V$.

The understanding of the motion in an $n$--centre potential, $n\geq3$,  
is very limited if one considers {\em negative} energies $E$. 
There one expects a complicated mixture of ergodic components
and motion on KAM tori.

However, by combining known techniques of celestial mechanics, 
we show in this article that the {\em high energy
motion} allows for a more or less complete qualitative 
and even quantitative description.

Up to an error of relative order $1/E$, 
the quantitative aspects treated here do
not depend on the precise form of 
the smooth potential $W$ but only on the
charges $Z_l$ and positions $\s_l$ of the nuclei. 

The qualitative structures do not even depend on these data but
only on the {\em number} $n$ of nuclei.
Thus the case $n=1$ of an atom resembles the Rutherford case, 
whereas the $n=2$--atomic molecule is similar to the two--centre problem
solved by Jacobi.
Here we are mainly interested in the case $n\geq 3$ where the dynamics is no
longer analytically integrable.

One motivation for this work is to establish the basis
for a geometric understanding of certain restricted $(n+1)$--body problems
of celestial mechanics, where one is interested in the motion of a fast
test body in the force field of $n$ bodies, whose motion is assumed to be known.

Such an analysis should be based on perturbation theory for the 
$n$--centre problem. 
In a joint work \cite{DK} with T.\ Dierkes, we intend to 
show that indeed the ($2$--dim.) $n$--centre problem is structurally stable in 
the following sense. For a sufficiently small local perturbation there exists
a homeomorphism conjugating the two phase portraits and leaving the
asymptotic initial and final directions of the scattering orbits invariant.
This homeomorphism is unique up to changes in the flow direction.

\paragraph{Survey of Results}

We now describe the techniques and results of this paper.
\medskip

The class of potentials $V$ under consideration
is introduced, together with other basic definitions,
in {\bf Sect.\ \ref{sect:defi}}. 
An important consequence of the fall-off of $V$ is
the existence of an {\em interaction zone} $\IZ$, 
a ball in configuration space which contains the points $\s_l$ 
and which, once left, 
cannot be reentered by an orbit, see (\ref{def:IZ}).
In Def.\ \ref{defi:noncol} we formulate the 
standing assumption that the centres are {\em non-collinear}, {\em i.e.}
no three centres are on a line.
\medskip

If $V$ contains a negative singularity ($Z_l>0$ for some $l$),
then the Hamiltonian flow generated by (\ref{ncp}) is incomplete.
In the Kepler problem it is well-known that the only sensible way to 
continue a {\em collision orbit} is to reflect it at the 
singularity in configuration space, for this 
is the limit behaviour for the Keplerian conic sections in the limit
of vanishing angular momentum.

But for our purposes we need to control smoothness and energy dependence
of the resulting flow, for potentials of the form (\ref{form:V}).
Due to the singularity at $\s_l$ the usual comparison techniques
for o.d.e.\ are not effective for proving such a result directly.

Instead, near $\s_l$ we apply in {\bf Sect.\ \ref{sect:KSdef} } 
the so-called Kustaanheimo-Stiefel transformation.
The {\em Hopf map} 
\[\bC^2\ar\bR^3,\qquad z\mapsto \bsm 
\LA z,\sigma_1 z\RA\\
\LA z,\sigma_2 z\RA\\
\LA z,\sigma_3 z\RA\esm\] 
with the Pauli matrices
$\sigma_l$, extends the Hopf fibration $S^3\ar S^2$. 
A cognate map of the phase spaces (cotangent bundles) is known 
\cite{KS,StS} 
as the {\em Kustaanheimo-Stiefel (KS) transformation}. 
It relates the positive energy 
Keplerian dynamics with the one of a particle in an inverted
harmonic potential.

The KS transformation was used in celestial mechanics (see, e.g.
Aarseth and Mikkola \cite{Aa} and articles by Aarseth and by Heggie
cited therein), applied
to spectral problems of Schr\"{o}dinger operators (see, e.g. \cite{HKSW}),
and shown in \cite{GK} to regularize the semiclassical dynamics. 

\medskip
The perturbation estimates of  {\bf Sect.\ \ref{sect:KS:Appl} } 
work in the covering phase space of the KS transformation.

Although the leading contribution to $V$ near $\s_l$ is the Kepler potential
\beq
-\frac{Z_l}{|\q-\s_l|},
\Leq{lc}
the following argument indicates that we cannot simply
approximate a collision orbit of (\ref{ncp}) by the Kepler hyperbola for
the potential (\ref{lc}) which has the same initial conditions $\pqs$.

The effective scattering region of (\ref{lc}) is a ball of radius $\cO(1/E)$
around $\s_l$, see Lemma \ref{lem:defl}.
On the other hand, the time the particle needs to reach the singularity
at $\s_l$ from its initial position $\q_0$ is of order $1/\sqrt{E}$.

Due to the (bounded) difference 
\[\nabla\l[V(\q)-\frac{-Z_l}{|\q-\s_l|}\ri]\]
of forces exerted on the two particles,
within this time the distance between the true orbit and the Kepler orbit
can grow to something of the order $1/E$. So by the above the Kepler orbit
starting at $\pqs$ may miss the effective scattering region.

But this means that shortly after collision the distance between 
the two trajectories does not necessarily decrease as $E\nearrow\infty$.

So instead we effectively
show that for every orbit with initial condition $\pqs$
there is a Kepler orbit, whose initial condition differs only by $\cO(1/E)$
(in a natural metric) from $\pqs$, 
and which remains in a $\cO(1/E)$--neighbourhood through
the whole collision process.
In Prop.\ \ref{propo:key}
the corresponding statements about perturbations
of the inverted harmonic oscillator are formulated.
More precisely, it is shown that these perturbed solutions
are $C^1$-near to the ones of the harmonic oscillator.

Moreover the true scattering process is approximated by Kepler scattering
in the $C^1$ sense, so that we may use the Kepler solutions if we want to
linearize the true flow (see Prop.\ \ref{prop:near:u}) . 
\medskip

Whereas the KS transformation is particularly suited  to understand the
topological and geometrical structures, its
disadvantages consist in
its local nature, the introduction of a
higher dimensional phase space and the reparameterization of time.
The first problem, the locality near one singularity, should not be
insurmountable. In fact Helffer and Siedentop found in 
\cite{HS} a generalization of the KS transformation to two centres.

Time reparameterization, however, is unwanted since we are
not only interested in the orbits but also in quantities like the time
delay of scattering orbits.
Thus we describe in {\bf Sect.\ \ref{sect:RL} } 
another regularization method.
In Theorem \ref{P:omega} we show that the incomplete Hamiltonian system 
\[(T^{*}\Muh,\omega_0,\Hh)\] 
may be uniquely extended to a {\em smooth} complete Hamiltonian system
\[(P,\omega,H)\]
whose phase space 
$P\supset T^{*}\Muh$ is a smooth six-dimensional manifold,
$\omega$ is a smooth symplectic two-form on $P$ restricting to the canonical 
symplectic two-form 
$\omega_0=\sum_{i=1}^3 dq_{i}\wedge dp_{i}$ on $T^{*}\Muh$,
and $H:P\ar \bR$ is a smooth Hamiltonian function with
$H\rstr_{T^{*}\Muh} = \Hh$.

In fact we linearize the flow near each 
negative singularity and then add a copy 
of $\bR\times S^2$ to $T^{*}\Muh$.
Here $\bR$ parametrizes the energy and $S^2$ the incoming (and outgoing)
direction of the collision orbit.

Thus we may henceforth work with the complete smooth flow
\[ \Phi:\bR\times P\ar P \qquad (t\in \bR)\]
generated by $H$. 

By Theorem \ref{P:omega} we may think of the Coulomb singularity
as an artefact of the use of the inappropriate phase space coordinates
$\pq$ which, however, leads to a non-trivial topology of $P$.

\medskip
In {\bf Sect.\ \ref{sect:moeller} } we introduce the M\o ller transformation 
by comparing the flow $\Phi^t$ on $P$ with the flow $\Pit:\Pin\ar \Pin$
generated by the Kepler Hamiltonian function 
\[\Hhi\pq := \eh \p^{\,2} - \frac{\Zi}{\B{\q}}\]
on its phase space $\Pin$, with $\Zi\in\bR$ defined by (\ref{def:Zi}).

The {\em M\o ller transformations} are then given by
\[\Opm := \lim_{t\ar\pm\infty} \Phi^{-t}\circ\Id\circ\Pit,\]
where $\Id$ canonically identifies the two phase spaces outside 
a region near the singularities.

The $\Opm$ are continuous and, under mild conditions on $V$ near infinity,
smooth canonical transformations, see Theorem \ref{thm:both:moeller}
and \ref{thm:smooth:moeller}. So the flow is {\em asymptotically complete},
see Corr.\ \ref{coro:complete}.

However, we are not primarily interested in these typical properties,
but in the specific traits of the multi-scattering dynamics.

It is a general fact that to a large extend the unbounded motion, i.e.,
{\em scattering} is determined by the {\em bounded} orbits of positive energy.
We denote by $\buE$ the set of such bound states of energy~$E$.

To control these orbits, we combine the above perturbative results for
the single scattering process with an analysis of the motion 
inside the interaction zone $\IZ$,  but away from the singularities
({\bf Sect.\ \ref{sect:between}}).
This turns out to be a 
$C^1$-perturbation of relative order $\cO(1/E)$ of the free motion. 
Using this result and the one of Sect.\ \ref{sect:si:sc},
we may approximate the true dynamics by 
a combination of free and of Keplerian motion.
\medskip

In {\bf Sect.\ \ref{sect:si:sc} } the perturbative results of 
Sect.\ \ref{sect:KS:Appl}
for the flow in KS space are used to obtain Prop.\ \ref{prop:near:u} 
which says that the single scattering process is $C^1$-near to pure
Kepler scattering, up to a relative error $\cO(1/E)$.
\medskip

Now in {\bf Sect.\ \ref{sect:long:paths} }
it is shown that
if a trajectory is not strongly reflected in uniformly bounded 
time by the singularities, it soon leaves the interaction zone 
(Prop.\ \ref{propo:drei:R}). 
\medskip

The bound states
are then analyzed in {\bf Sect.\ \ref{sect:P:map} }
by a {\em Poincar\'{e} section} technique.
We erect Poincar\'{e} surfaces which, in their projections to configuration
space, sit between pairs of singularities and then define in
(\ref{Poincare:map}) a Poincar\'{e} map $\Po$.
\medskip

In Proposition \ref{propo:cf} of
{\bf Sect.\ \ref{sect:ICF} } we estimate the linearization of $\Po$,
which, up to a relative error $\cO(1/E)$, only depends on the scattering angle
at $\s_l$ and the charge~$Z_l$.
\medskip

This allows us to establish in {\bf Sect.\ \ref{sect:symbol}}, 
Thm.\ \ref{thm:homeo} a symbolic dynamics for $\buE$,  $E>\Eth$.
The bounded orbits are described by their sequence of near-collisions,
(which is well defined, since we assumed the configuration of centres to
be non-collinear).
 
Thus for $n=1$ there are no bounded orbits, 
for $n=2$ the set $\buE$ consists of one closed orbit, which is closed.
For the case $n\geq3$ of primary interest, 
$\buE$ is locally homeomorphic to the product of a Cantor set and
an interval.

All bounded orbits are hyperbolic. $\buE$ has measure zero w.r.t.\ Liouville
measure on the energy shell.
The Morse index of a periodic trajectory equals its number of near-collisions.

\medskip
In {\bf Sect.\ \ref{sect:fractal} } the fractal dimension of this set 
$\buE$ of bounded orbits is estimated. More precisely, we consider its 
{\em Hausdorff dimension} 
$\dimH$ and its upper {\em box counting dimension} $\dimB$, 
since it is known that most other 
definitions lead to numerical values between these two.

Theorem \ref{thm:fractal} says that for energies $E$ above 
the threshold $\Eth$ and $n\geq 3$ centres
\beqno
\lefteqn{\hspace{-8mm}1+2d(E)\cdot\l(1-\cO((E\ln E)^{-1})\ri)
\leq\dimH(\buE)\leq}\NN\\
&\leq& \dimB(\buE)
\leq 1+2d(E)\cdot\l(1+\cO((E\ln E)^{-1})\ri),
\eeqno
with the solution $d(E)$ of a finite matrix eigenvalue problem (\ref{ev:eq}).
In particular they meet the rough estimate
\[\dimH(\buE)=1+\frac{2\ln(n-1)}{\ln(E)}+\cO\l((\ln E)^{-2}\ri)=
\dimB(\buE)\]
\medskip
In {\bf Sect.\ \ref{sect:top:ent} }
the topological entropy of the flow $\Phi_E^t$ of energy $E$
is estimated.

Topological entropy is a quantity which, roughly speaking, measures
the information loss per time unit about the state of the system.
Here for positive energy $E$  the energy shell $\SuE$ is non-compact,
so that we use Bowen's definition of entropy.
Prop.\ \ref{prop:restriction} states that 
\[\htop(\Phi_E^{1}) = \htop(\Phi_E^{1}\rstr_{\buE}),\]
{\em i.e.} that the source of information loss is the intricate structure
of the set $\buE$ of bounded orbits, whereas the scattering orbits
only have a sub-exponential dependence on initial conditions.

So the estimates in the proof of Thm.\ \ref{thm:top:ent}
can be based on symbolic dynamics. It states that for $E$ large 
$\htop(\Phi_E^{1})=0$ for $n=1$ or $2$ centres, whereas for 
$n\geq3$
\[\htop(\Phi_E^{1})=
\htopi\cdot\sqrt{2E}\cdot
\l(1+\frac{\ln(E)}{E}C_\htop + \cO(1/E)\ri).\]
Here the constants $\htopi$ and $C_\htop$
are determined by solving a finite matrix eigenvalue problem.

Whereas the factor $\sqrt{2E}$ is of kinematical nature, 
in the simplest case of an equilateral triangle ($n=3$) resp.\ tetrahedron 
($n=4$) of side length $d$  
the constant $\htopi$ equals $\ln(n-1)/d$.

This divergence of topological entropy is not in contradiction with 
integrability (in the sense of independent $C^\infty$ integrals of motion) 
of the hamiltonian dynamics above the energy threshold (compare with
Bolsinov and Taimanov \cite{BT}).

\medskip

{\bf Sect.\ \ref{sect:sc} } 
is devoted to the classification of the scattering orbits.
After excluding asymptotic directions in cones of aperture
$\cO(1/\sqrt{E})$ around the axes through
two nuclei, and near the forward direction, 
one obtains such a universal classification.
Thm.\ \ref{thm:sc} states that the orbits of given energy 
and asymptotic directions are enumerated by the succession of the 
nuclei they visit.
In particular they form a Cantor set if $n\geq3$.

\medskip
The differential cross section $\DCSP$ of the scattering process
is analyzed in {\bf Sect.\ \ref{sect:diff:cross}}. Roughly speaking,
this experimentally accessible quantity is the 
`probability' that a particle of energy $E$ and initial
direction $\htheta^{-}$ has final direction $\htheta^{+}$.
For general potential scattering in $\bR^d$ 
one cannot even expect the cross section measure on 
$S^{d-1}\setminus\htheta^{-}$
(defined in Def.\ \ref{defi:diff:cross}) to be absolutely continuous  
w.r.t.\ Haar measure, see \cite{Kn3}. Here, however, this is the case 
for energies $E>\Eth$, and $\DCSP$
is defined as the Radon-Nikodym derivative.

Moreover (after excluding cones of $E$--independent aperture)
by Thm.\ \ref{theorem:classify:scattering} 
it differs from the well-known {\em Rutherford cross section} in $\bR^d$
\[\l(\DCSP\ri)_{\rm Ru} =\l(\frac{Z}{4E\sin^2(\eh\Delta\theta)}\ri)^{d-1}\]
of the single Coulomb potential
with charge $Z =\sqrt{{\textstyle\sum_{l=1}^n Z_l^2}}$ only by 
\[\DCSP =\l(\DCSP\ri)_{\rm Ru} \cdot(1+\cO(1/E)).\]
So for these asymptotic data the intricate structure of the scattering orbits 
is not showing up in the cross section. In fact, 
Thm.\ \ref{theorem:classify:scattering} 
also states that the differential cross section is approximated by
the Rutherford cross section outside the (much smaller) 
cones of aperture $\cO(1/\sqrt{E})$, though with less accuracy.

\medskip
In the final {\bf Sect.\ \ref{sect:collinear}} we show by counterexamples 
that many results of the paper do not generalize if we drop the
assumption of non-collinearity.

\medskip
The first {\bf Appendix} is devoted to a comparison between the two-dim.\
case treated in \cite{KK} and the three-dim.\ situation of the present
paper. Basically, the analysis of \cite{KK} is the one of geodesic motion
on smooth many-handled surfaces of negative curvature, whereas 
here we apply perturbation techniques around infinite energy.\medskip\\ 
In the second Appendix we describe the (known) solution of
the purely Coulombic two-centre problem, and its bifurcation
diagramme.\medskip\\[5mm]
{\bf Notation.} 
Transposed vectors are used according to typographical, not
mathematical needs.   \\[5mm]
{\bf Acknowledgement.}
I thank the Max Planck Institute for Mathematics in the Sciences
(Leipzig), where a part of this paper was written, 
for hospitality and support.
%
\Section{Basic Definitions} \label{sect:defi}
%
We consider the time evolution generated by a Hamiltonian function
\beq
\Hh\pq := \eh \p^{\, 2}+ V(\q)
\Leq{def:hamiltonian}
with $n$ Coulombic singularities of the potential $V$
situated at the points 
\[\label{eq:sing:l}
\s_{1},\ldots,\s_{n}\in \Mu \qquad
(\s_{i}\neq \s_{k}\qmbox{for} i\neq k).\] 
To control the asymptotic behaviour,
we assume that $V$ decomposes into the sum of a purely Coulombic potential
and a short range potential.
By this we mean the following:
\begin{definition} \label{defi:coulombic} 
A smooth, real-valued function $V$ on the {\em configuration space}
\beq
\Muh := \Mu\setminus \{\s_{1},\ldots,\s_{n}\} 
\Leq{M1}
is called {\em Coulombic} if
\begin{enumerate}
\item 
There exist $Z_{l}\neq 0$, $l\inInd$, such that
\beq
V(\q) = \sum_{l=1}^{n}\frac{-Z_{l}}{\B{\q-\s_{l}}}+ W(\q) \qquad (\q\in\Muh)
\Leq{def1:V}
with $W:\Mu\ar \bR$ smooth.
$Z_{l}$ is called the {\em charge} of the $l$th {\em nucleus}, and we set
\[\Zmax:=\max\{|Z_1|,\ldots,|Z_n|\}.\]
\label{def:Zl}
\item The potential vanishes at infinity, {\em i.e.}\
\beq
\lim_{\B{\q}\ar\infty}V(\q)=0,
\Leq{V:ar:0}
and there exist $\Zi\in\bR$, called the {\em asymptotic charge},
$\eps\in(0,1]$ and 
\[\Rmin>2\max(|\s_1|,\ldots,|\s_n|)\]
such that for some $C_1>0$ 
\beq
\B{\nabla V(\q) - \Zi\frac{\q}{\B{\q}^{3}}} < 
\frac{C_1\cdot\Rmin}{\B{\q}^{2+\eps}}\qquad(\B{\q}\geq \Rmin)
\Leq{def1:u1}
and
\beq
\B{\nabla V(\q_{1}) - \nabla V(\q_{2})} < 
C_1\frac{\B{\q_{1} - \q_{2}}}{\min(\B{\q_{1}},\B{\q_{2}})^{2+\eps}}
\qquad(\B{\q_1},\B{\q_2}\geq \Rmin).
\Leq{def1:u2}
\end{enumerate}
\end{definition}
\begin{example} 
For the class of {\em purely Coulombic} potentials
\beq
V(\q) := \sum_{l=1}^{n} \frac{-Z_{l}}{\B{\q-\s_{l}}},
\Leq{simple}
the asymptotic charge $\Zi=\sum_{l=1}^{n} Z_{l}$.
Here $V$ meets (\ref{def1:u1}) and (\ref{def1:u2})
with $\eps=1$ (and $C_1=\frac{17}{2}n\Zmax$).
\end{example}
The question of dynamics in a Coulombic potential is called the 
{\em $n$--centre problem} of classical mechanics.

\begin{remarks} 
{\bf 1)}
For $Z_l>0$ and $W=0$
this corresponds to the idealization of a celestial body
in the force field of $n$ other bodies of masses $Z_l$ with  
fixed positions $\s_l$.

But our definition also covers the physical situation
of (classical) scattering by the potential of partially ionized, 
neutral ($\Zi=0$) or negatively charged quantum molecules. There
$\Zi$ does not coincide with $\sum_{l=1}^{n} Z_{l}$.
Again, one expects potentials with $\eps=1$, 
due to exponential decay of the quantum mechanical
eigenfunctions of the bound electrons, which in turn leads to an
electronic charge distribution which decays exponentially
(see Agmon \cite{Ag}).\\
{\bf 2)}\label{rem:symm}
Note that $V(\q)$ and the asymptotic potential $-\Zi/\B{\q}$
appear symmetrically in (\ref{def1:u1}), and that (\ref{V:ar:0}) 
and (\ref{def1:u2}) are met by the asymptotic potential 
(the last one even with $\eps=1$). 
We choose $C_1$ large enough so that with the given constant
$\eps$ (\ref{def1:u1}) is met by both potentials. 
This will be used in Sect.\ \ref{sect:moeller} for simplifying
the existence proof of the M\o ller transformations.\\
{\bf 3)}
By a suitable translation of the origin, one could for $n\geq2$ assume that 
$\Rmin$ is, say, twice the maximal distance between 
the centres (but this assumption will not be used).
\end{remarks}
For $n\geq 2$ we set
\beq
\dmin := \min_{k\neq l}d^{k,l}\qmbox{and}
\dmax := \max_{k\neq l}d^{k,l}\qmbox{with} d^{k,l}:=\B{\s_{k}-\s_{l}}.
\Leq{d:min}
In the case of a single atom ($n=1$) we fix $\dmin$ by setting
$\dmin := 2\Rmin$, say.

A threshold of the energy is 
\beq
\Vmax :=  
\sup \l( \l\{V(\q)\mid \q\in\Muh,\nabla V(\q)=\vec{0}\ri\}\cup\{0\}\ri).
\Leq{def:Vmax}
One has $0\leq \Vmax < \infty$, and in many cases $\Vmax = 0$.

We will generally assume that there are no more than two
nuclei on one line: 
\begin{definition} \label{defi:noncol} %
The configurations of singularities in
$$\NC := \l\{(\s_1,\ldots,\s_n)\in (\Mu)^n\mid 
\forall i\neq j\neq k\neq i: 
(\s_i-\s_j)\times(\s_j-\s_k)\neq \vec{0} \ri\}$$
are called {\em non--collinear} (NC) configurations.
\end{definition}
In three dimensions this is a weaker assumption 
than that the nuclei are in general position,
since the latter also means that there are no more than three 
nuclei in one plane. 

The space $\NC$ of configurations
is a smooth connected manifold.

Most statements in this article 
are shown for values of $H$ above some
threshold energy $\Eth$. For given $W$ (in particular for the 
purely Coulombic case  
$W=0$) and charges $Z_l$ this is a function $\Eth:\NC\ar\bR$
of the positions $(\s_1,\ldots,\s_n)$ of singularities.
$\Eth$ diverges near the set $(\Mu)^n\setminus\NC$ of partly collinear
configurations. 
This is not an artefact of our methods of proof.
Instead, we will see in Sect.\ \ref{sect:collinear} that several 
assertions proved in previous sections become wrong if one drops the NC condition.
 
One parameter measuring the degree of non--collinearity of a configuration
of $n\geq 3$ nuclei is the {\em minimal angle} 
\beq
\amin:\NC\ar (0,\pi/3]\qquad \amin := \min_{i\neq j\neq k\neq
i}\alpha(i,j,k),
\Leq{def:amin}
where for $i\neq j\neq k$
\beq
\alpha(i,j,k)\in [0,\pi),\qquad
\cos(\alpha(i,j,k)):=\LA\hat{s}^{j,i},\hat{s}^{j,k} \RA,
\Leq{a:ijk}
with the directions 
\[\hat{s}^{i,j}:= \frac{\s_j-\s_i}{|\s_j-\s_i|},\]
is the angle between $\s_i$ and $\s_k$, seen from $\s_j$.
For $n=2$ we set $\amin:= \pi/3$. \\[5mm]
The Hamilton equations
\[\dot{\p}=-\nabla V(\q)\qmbox{,}\dot{\q}=\p\]
lead to solutions $t\mapsto(\p(t),\q(t))\equiv\Pt(x_0)$ of the initial
value problem  with initial values $x_0=\pqs$ which exist uniquely 
up to eventual collisions with the $\s_l$. In Sect.\ \ref{sect:RL}
we extend $\Pt$ uniquely to a smooth complete flow.

The {\em virial identity}
\beq
\frac{d}{dt} \LA\q(t),\p(t) \RA = 2(E-V(\q(t)))- \LA\q(t),\nabla V(\q(t))\RA 
\Leq{virial}
holds true for any trajectory $t\mapsto (\p(t),\q(t))$ with 
energy $E := H(x_0)$ (whenever $\q(t)\neq \s_{l}$).
Let us choose a function
$\Rvir:(0,\infty)\ar\bR$ of the energy, called 
the {\em virial radius}, with
\beq
\max\l(|V(\q)|,\frac{|\Zi|}{|\q|}\ri) < \frac{E}{2} \qmbox{and} |\LA\q, \nabla V(\q)\RA| < E/2
\qquad(\B{\q}\geq \Rvir(E)). 
\Leq{V:small}
As a consequence of part 2 of Def.\ \ref{defi:coulombic} of 
Coulombic potentials, 
such a function exists. Without loss of generality (w.l.o.g.) we assume
$\Rvir$ to be continuous, nonincreasing and constant for energies $E>\Eth$
above some threshold. Property (\ref{V:small}) already implies
that $\Rvir(E)>2\Zi/E$. Technically we choose $\Rvir$ with
\beq
\Rvir(E)\geq\max(2\Rmin, C_2/E)\qmbox{with}
C_2:= 31(1+1/\eps)\Rmin^{1-\eps} C_1
\Leq{R:vir}

Then by (\ref{virial}) and (\ref{V:small})
\beq
\frac{d}{dt} \LA\q(t),\p(t) \RA > \frac{E}{2} > 0 \qquad\mbox{if }
 \B{\q(t)}\geq \Rvir(E).
\Leq{qp}
We conclude that a configuration space 
trajectory $t\mapsto\q(t)$ of energy $E$ 
leaving the ball of radius $\Rvir(E)$
cannot reenter this ball in the future but must go to spatial infinity:

Namely assume w.l.o.g.\ that $\LA\q(0),\p(0)\RA\geq 0$. By (\ref{qp})
\[\frac{d^2}{dt^2}\q^{\,2}(t) = 
2\frac{d}{dt}\LA\q(t),\p(t)\RA > E\qquad (t\geq 0)\]
so that 
\beq
\q^{\,2}(t) \geq \q^{\,2}_0 + \eh Et^2 \qquad (t\geq 0).
\Leq{q:grows}
We shall mainly deal with energies $E>\Eth$ so that we may consider the 
$E$-independent {\em interaction zone}
\beq
\IZ := \{\q\in \Mu\mid \B{\q} \leq \Rvir(\Eth)\}.
\Leq{def:IZ}
Sometimes we will use outside $\IZ$ the shorthand
\beq
\Vsr:=V-V_\infty \qmbox{with} V_\infty(\q):=-\frac{\Zi}{|\q|}
\Leq{V:sr}
In addition to the condition (\ref{def1:u1}) on $\nabla\Vsr$
this potential is of short range in the sense $\Vsr(\q)=\cO(|\q|^{-1-\eps})$,
too, as follows from (\ref{V:ar:0}) and  (\ref{def1:u1}).

Due to collisions with the nuclei situated at $\s_{l}$, the flow on the 
phase space $T^{*}\Muh$ of the particle is incomplete. 
There are several ways to regularize the collision orbits which are all 
essentially equivalent. In Section \ref{sect:KSdef} we now
introduce the Kustaanheimo-Stiefel (KS)
regularization method, and use it in Sect.\ \ref{sect:KS:Appl} for 
comparison estimates which control the deviation from the Keplerian motion 
near a singularity.
Later, in Sect.\ \ref{sect:RL}, the flow is then regularized without
a time change.
%
\Section{The Kustaanheimo--Stiefel Transformation} \label{sect:KSdef}
%
In \cite{KS} Kustaanheimo and Stiefel related the Kepler motion
in three spatial dimensions to the motion of a resonant harmonic oscillator
in four dimensions, thus linearizing the dynamics, see the book \cite{StS}
by Stiefel and Scheifele. 

Although the authors had used spinor theory in order to derive their results,
the emphasis of their article was more on the application to perturbation
theory then on the geometry of the problem.
In \cite{Ku} Kummer presented this aspect of the KS-transformation and related 
it to the approach \cite{Mo} by Moser. 

Our presentation will be based on the quaternion algebra over $\bR$
\[\bH := \l\{\l.\bem{cc} z_1 & -z_2\\ \bar{z}_2 &
\bar{z}_1\eem 
\ri| z_1,z_2\in \bC\ri\} \cong\bR^4\]  
with matrix multiplication (see, e.g., \cite{KR}). We use the basis 
\beqno
(\idty,I_1,I_2,I_3)&:=&\l(\bsm 1& 0\\0& 1\esm,\bsm 0& i\\i& 0\esm,
\bsm 0& 1\\-1& 0\esm, \bsm i& 0\\0& -i\esm\ri)\\
&=& (\idty,i\sigma_1,i\sigma_2,i\sigma_3)
\eeqno 
of $\bH$, the $\sigma_l$ being the Pauli matrices.
The direct sum decomposition 
\[\bH= \bR\cdot \idty \oplus {\rm Im} \bH\] 
with 
\beqno
\ImH &:=& \{Z\in\bH\mid Z^2 = \lambda \cdot \idty\mbox{ with }
\lambda \leq 0 \}\\
&=& {\rm Span}_\bR (I_1,I_2,I_3)
\eeqno
into real and imaginary space
is orthogonal w.r.t.\ the inner product
\[\bH\times\bH \ar\bR \quad,\qquad \LA X,Y\RA := \eh \tr(XY^*),\]
$X\mapsto X^*:=\bar{X}^t$ being the conjugation.
The norm $|X|:=\LA X,X\RA^{\eh}$ is multiplicative:
\[|XY| = |X|\,|Y|\qquad (X,Y\in\bH).\]
The vector product $\times:\ImH\times\ImH\ar \ImH$ is given by
\[X\times Y := \eh(XY-YX),\]
and we have 
\[XY= - \LA X,Y \RA \idty+ X\times Y \qquad(X,Y\in \ImH).\]
We consider the {\em Hopf map} 
\beq
\pi_0:\bH\ar\ImH,\quad \pi_0(Z):=Z^*I_3 Z=
i\bem{cc} 
z_1\bar{z}_1-z_2\bar{z}_2 & -2\bar{z}_1 z_2 \\ 
-2z_1 \bar{z}_2      & z_2\bar{z}_2-z_1\bar{z}_1\eem 
\Leq{Hopf:map}
which is a surjection $\bR^4\ar\bR^3$ whose preimages are the orbits of the
isometric group action
\[\alpha_0:S^1\ar{\rm Aut}(\bH),\qquad\alpha_0(\vv)(Z):= \exp(I_3\vv)Z.  \]
This action is free on $\hH:= \bH\setminus \{0\}$.
We have the canonical symplectic 
one-forms 
\[\btheta := \eh\tr(P^*dQ)=\Re(\bar{p}_1dq_1+\bar{p}_2dq_2)
\qquad ((P,Q)\in T^*\bH)\]
and 
\[\theta := \eh\tr(P^*dQ)\qquad ((P,Q)\in T^*\ImH)\]
on the cotangent bundles, and denote by 
$\hat{\btheta}:= \btheta\rstr_{T^*\hH}$, resp.\  
$\hat{\theta}:= \theta\rstr_{T^*{\rm Im}\hH}$ their restrictions.

The restricted Hopf map $\pih_0:=\pi_0\rstr_{\hH}$ is then related to
the {\em KS--trans\-for\-ma\-tion} 
\beqn
\hpi &:&T^*\hH\ar T^*{\rm Im}\hH \NN\\ 
\hpi(P,Q) &=& 
\l({\textstyle \frac{-1}{4|Q|^2}}\l(Q^{*}I_3 P+P^*I_3Q\ri),Q^*I_3 Q\ri)
\label{def:KS}
\eeqn
of the cotangent bundles.

We consider the quadric surface 
\[S:= \cI^{-1}(0)\subset T^*\bH\]
for the bilinear form
\[\cI:T^*\bH\ar \bR,\qquad \cI(P,Q):= \eh\tr(Q^*I_3 P) = 
\Im(q_1\bar{p}_1+q_2\bar{p}_2) \]
and its restriction $\hat{S} := S\cap T^*\hH$.
Then 
\beq
\hat{\btheta}\rstr_{\hat{S}} = \hpi^*\hat{\theta}\rstr_{\hat{S}},
\Leq{lift:one:form}
since on $\hat{S}$ 
\beq
Q^*I_3 P=P^*I_3 Q\in \ImH.
\Leq{eq:QIP}

$\cI$, $\btheta$ and $\hpi$ are all invariant w.r.t.\ the group action
\beq
\alpha:S^1\ar{\rm Aut}(T^*\bH),\qquad\alpha(\vv)(P,Q):=
(\exp(I_3 \vv)P,\exp(I_3 \vv)Q),
\Leq{circle:action}
and $\cI$ is a Hamiltonian function generating that time--$\vv$--flow.

\begin{lemma} \label{lem:lifted:fcts}
Restricted to the quadric surface $\hat{S}\subset T^*\bH$, the lift of a 
Hamiltonian function $\Hh\pq=\eh |\p|^2-Z/|\q|+W(\q)$ equals
\beq
\Hh\circ\hpi(P,Q) = 
|Q|^{-2}\l( \ea|P|^2-Z\ri)+W(Q^*I_3 Q)\qquad ((P,Q)\in\hat{S}).
\Leq{lift:H}
For the angular momentum $\vec{L}\pq = \q\times \p$
\beq
\vec{L}\circ \hpi(P,Q) = \ev(Q^*P-P^*Q) \qquad ((P,Q)\in\hat{S}),
\Leq{lift:L}
whereas the Runge-Lenz vector $\vec{F}\pq=\p\times\vec{L}-Z\frac{\q}{|\q|}$
transforms into
\beq
\vec{F}\circ \hpi(P,Q) = 
-\ea P^*I_3 P + (\Hh\circ\hpi(P,Q)-W(Q^*I_3Q))\cdot Q^*I_3 Q  
\quad ((P,Q)\in\hat{S}).
\Leq{lift:F}
\end{lemma}
{\bf Proof.} (see also \cite{Ku}).\\
By (\ref{eq:QIP}) the restricted KS--transformation 
$\pi\rstr_{\hat{S}}$ maps $|Q|^2=|Q^*I_3 Q|$ onto $|q|$
and $|2Q|^{-2} \cdot|P|^2=|(2|Q|^2)^{-1}Q^*I_3 P|^2$ onto $|p|^2$, implying 
(\ref{lift:H}).

Since $q,p\in \ImH$, 
\[\vec{L} = q\times p= \eh(q p - pq).\]
The KS-transformation $\hpi$ gives, using (\ref{eq:QIP})
\beqno
\vec{L}\circ \hpi(P,Q) &=& \frac{(Q^*I_3P+P^*I_3Q) Q^*I_3 Q-Q^*I_3 Q(Q^*I_3P+P^*I_3Q) }{8|Q|^2}\\
&=&  \frac{P^*I_3QQ^*I_3 Q-Q^*I_3 QQ^*I_3P }{4|Q|^2} =\ev(Q^*P-P^*Q),
\eeqno 
This shows (\ref{lift:L}).
To prove (\ref{lift:F}), we use the identity
\[p\times(q\times p) = \eh(pqp-qpp)=\eh(pqp+q|p|^2)\qquad (p,q\in\ImH)\]
which implies that 
\[\vec{F}\pq = \eh pqp + q\cdot(H(p,q)-W(q)).\]
The term $pqp$ lifts to
\[ (2|Q|^2)^{-2} (P^*I_3 Q)(Q^*I_3 Q)(Q^*I_3 P) = -\ev P^*I_3 P, \]
finishing the proof.
\hfill $\Box$
%
\Section{Application of the KS--Transformation} \label{sect:KS:Appl}
%
In a configuration space ball 
\[B_l(r) :=\{\q\in \Mu\mid |\q-\s_l| \leq r \}\]
of radius $r:=c_q<\eh\dmin$ around $\q=\s_l$ we want to regularize the motion 
generated by the restricted Hamiltonian function 
\[\Hh(\p,\q) = \eh|\p|^2-\frac{Z_l}{|\q-\s_l|}+ W_l(\q),\qquad 
\l((\p,\q)\in T^*(B_l(c_q)\setminus\{\vec{0}\})\ri),\]
\beq
W_{l}(\q) := \sum_{i\neq l} \frac{-Z_{i}}{\B{\q-\s_{i}}} + W(\q)
\Leq{Wl}
being a {\em smooth} function on $B_l(c_q)\subset\bR^3\cong \ImH$.
The radius $c_q\in(0,\eh\dmin)$ is chosen so that 
\beq
|W_l(\q)|\leq\eh \frac{|Z_l|}{|\q-\s_l|}\qquad (l=1,\ldots,n\,,\,|\q-\s_l|\leq c_q).
\Leq{cq:small}
For simplicity of notation we assume $\s_l=\vec{0}$. 
For regular values $E$ of the energy $\Hh$ the orbits on $\Hh^{-1}(E)$ 
coincide with the ones of the
zero surface of the Hamiltonian function 
\beq
|q|\cdot(\Hh(p,q)-E),
\Leq{to:lift}
since they are already determined by the form of the submanifold $\Hh^{-1}(E)$.

By (\ref{lift:H}) the lift of (\ref{to:lift}) with the KS-transformation equals
\beq
{\bf H}_E(P,Q):= \ea\LA P,P\RA +\LA Q,Q\RA (-E+W_l(Q^*I_3 Q)) - Z_l.
\Leq{def:HE}
The orbits of its Hamiltonian vector field on 
$\hat{S}\cap({\bf H}_E)^{-1}(0)$ project to the ones of $\Hh$ on $\Hh^{-1}(E)$,
since the symplectic one-form $\hat{\theta}$ on the phase space region 
$T^*(B_l(c_q)\setminus\{\vec{0}\})$ lifts according to (\ref{lift:one:form}).

For $W_l\equiv 0$ the Hamiltonian function ${\bf H}_E$ 
is the one of a four-dimensional harmonic oscillator (with negative potential
for $E>0$). In the general case the additional potential 
$Q\mapsto W_l(Q^*I_3 Q)$ will be a small perturbation if $E$ is large, 
at least for $|Q|\leq \sqrt{c_q}$.

From now on we assume initial conditions to lie in $S\subset T^*\bH$. 

Instead of considering the Hamilton equations
\beqno
\l\{\begin{array}{rcl}
\frac{d}{ds} P &=& 2EQ +2R(Q)\\
\frac{d}{ds} Q &=& \ev P 
\end{array} \ri.
\eeqno
of (\ref{def:HE}) with perturbation
\beq
R(Q):= -Q W_l(Q^*I_3 Q) - \eh|Q|^2 \nabla W_l(Q^*I_3 Q),
\Leq{def:R}
we want the parameter $E$ to appear only in the perturbative term. Thus we set 
\[X := {\tilde{P} \choose Q} \qquad{\rm with }\quad \tilde{P}:=P/\sqrt{8E}\]
and use the time variable $\tau:= \sqrt{E/2}\cdot s$. Then
\beq
\frac{d}{d\tau} X = \bem{cc}0& \idty\\ \idty& 0\eem X
+ \tilde{R}(X)/E\quad,\qquad X(0)=X_0  \equiv {\tilde{P}_0\choose Q_0}  
\Leq{init:val}
with $\tilde{R}(X) := {R(Q)\choose 0}$.

The initial conditions $X_0 := {P_0/\sqrt{8E}\choose Q_0}$
meet the energy constraint
\beq
|\tilde{P}_0|^2-|Q_0|^2(1-W_l(Q_0^*I_3 Q_0)/E)= \frac{Z_l}{E}.
\Leq{en:con}
So with 
\[c_Q:=\sqrt{c_q}\] 
and 
\beq
\cE(\tilde{P},Q) :=\frac{Z_l-|Q|^2 W_l(Q^*I_3 Q)}{|\tilde{P}|^2-|Q|^2}
\Leq{En:o} 
the region 
\[\Dol := \l\{ \l. (\tilde{P},Q)\in\bR^4\times\bR^4\ri| |Q|\leq c_Q, 
\cE(\tilde{P},Q) > \Eth \ri\} \]
of phase space points meeting (\ref{en:con}) for some $E> \Eth$ is bounded.

It is natural to compare the solution $X(\tau)$ of (\ref{init:val})
with the solution 
\beq
Y(\tau) =  
\bem{cc}\cosh(\tau)\idty & \sinh(\tau)\idty\\ 
\sinh(\tau)\idty& \cosh(\tau)\idty\eem  X_0
\Leq{l:init:val}
of the linear initial value problem 
\beq
\frac{d}{d\tau} Y = \bem{cc}0& \idty\\ \idty& 0\eem Y
\qmbox{,} Y(0)=X_0   .
\Leq{init:val2}
We denote this linear flow 
(corresponding to $W_l\equiv 0$), restricted to the invariant phase space domain
\[\DoLlb := \l\{ \l. (\tilde{P},Q)\in\bR^4\times\bR^4\ri| \ |\tilde{P}|\neq
|Q|\,,\,
\cE_L(\tilde{P},Q) >  \eh\Eth \ri\} \]
with 
\beq
\cE_L(\tilde{P},Q):= Z_l/(|\tilde{P}|^2-|Q|^2)
\Leq{En:Lo} by  
\beq
\PoLl:\bR\times\DoLlb \ar\DoLlb,
\Leq{linear:flow} 
and set
\[\DoLl := \{ (\tilde{P},Q)\in\DoLlb\mid |Q|\leq c_Q\}.\]
Our condition (\ref{cq:small}) on the radius $c_q$ implies that
\[\DoLl\supset \Dol.\]
Let 
\[\Tol^\pm: \Dol\ar\bR\cup\{\pm\infty\},\quad 
\Tol^\pm(X_0) := \pm\inf \{t > 0\mid |Q(\pm t,X_0)|\geq c_Q\}.\]
be the exit times from $\Dol$,
\[\Pol:{\bf U}_l \ar\Dol, \quad 
{\bf U}_l := \l\{ (t,X_0) \in \bR \times \Dol \mid t\in [\Tol^-(X_0),\Tol^+(X_0)]\ri\}\]
the maximally extended KS flow on $\Dol$, and
\[\Psol^\pm:\Dol\ar\pa\Dol \qmbox{,} X\mapsto\Pol(\Tol^\pm(X),X)\]
the map to the exit points.

Basically we are interested in the Poincar\'{e} map
\beq
\Psol:\pa\Dol\ar\pa\Dol ,\quad X\equiv (\tilde{P},Q)\mapsto 
\l\{\begin{array}{ll}\Psol^+(X) & ,\LA \tilde{P},Q\RA \leq 0\\
\Psol^-(X) & ,\LA \tilde{P},Q\RA > 0\end{array} \ri.
\Leq{exit}
that permutes incoming and outgoing data, but up to now we do not even know 
whether this is defined everywhere.
 
Therefore we compare with the linear flow (\ref{linear:flow}) 
and thus introduce in analogy
its exit times $\ToLl^\pm: \DoLl\ar\bR$. They are finite and smooth, 
and we set $\PsoLl^\pm(X):=\PoLl(\ToLl^\pm(X),X)$. 

The analog
\beq
\PsoLl:\pa\DoLl\ar\pa\DoLl ,\quad X\equiv (\tilde{P},Q)\mapsto 
\l\{\begin{array}{ll}\PsoLl^+(X) & ,\LA \tilde{P},Q\RA \leq 0\\
\PsoLl^-(X) & ,\LA \tilde{P},Q\RA > 0\end{array} \ri.
\Leq{exit:L}
of (\ref{exit}) for the linear flow is a involutive diffeomorphism:
\beq
\PsoLl\bem{c}\tilde{P}\\ Q\eem = \frac{1}{\sqrt{1-u^2}}
\bem{c} \tilde{P}-u Q\\ Q - u\tilde{P}\eem\qmbox{with} 
u:=\frac{2\LA\tilde{P},Q\RA}{|\tilde{P}|^2+|Q|^2},
\Leq{PsoLl}
since $|u|<1$ on $\DoLl$.

Ideally, one would like to prove that (\ref{exit}) is approximated by
the map (\ref{exit:L}) so that the $C^1$-norm of
$\Psol\circ(\PsoLl)^{-1}-\Id$ is of order $\cO(1/\cE)$.  
This would be true if the trajectories of $\Pol$, resp.\ $\PoLl$ 
would spend a uniformly bounded time
inside $\Dol$, resp.\ $\DoLl$. But this is not the case, and 
in order to keep the error terms small 
we will compare the two flows for initial conditions on the pericentric
hypersurface
$$\Hol := \l\{ (\tilde{P},Q)\in\Dol\mid \langle \tilde{P},Q\rangle=0\ri\}
\subset \HoLl := 
\l\{ (\tilde{P},Q)\in\DoLlb\mid \langle \tilde{P},Q\rangle=0\ri\}\!.$$
Note that by transversality of the linear flow to that hypersurface 
the pericentric time
\beq
\ToLl^0:\DoLlb \ar\bR \qmbox{with} \PsoLl^0(X):=\PoLl(\ToLl^0(X),X)\in\HoLl
\Leq{ToLl:0} 
is uniquely defined and smooth.
\begin{proposition} \label{propo:key}
For $\Eth$ large there is a unique pericentric 
time $\Tol^0:\Dol\ar\bR$ with 
\[\PsolN(X):=\Pol(\Tol^0(X),X)\in\Hol,\qquad (X\in \Dol).\]
The functions $\Tol^-\leq\Tol^0\leq \Tol^+:\Dol\ar\bR$ are smooth, and
\beq
|\Tol^\pm(X_0)- \ToLl^\pm(X_0)|= \cO(1/\cE(X_0))\qquad (X_0\in\Hol).
\Leq{time:int}
The exit times $\ToLl^\pm$ of the linear flow are estimated by
\beq
\exp\l(\pm 2\ToLl^\pm(X_0)\ri) = 
\frac{4c_q}{|Z_l| e(X_0)} \cE(X_0) +\cO(\cE^0(X_0))\qquad
(X_0\in\Hol),
\Leq{L:time:o}
$e(X_0):= \sqrt{1+\frac{2E|\vec{L}(X_0)|^2}{Z_l^2}}$ 
being the eccentricity of the corresponding Kepler hyperbola.
The diffeomorphism 
\[\Xio_l: \pa \Dol\ar\pa\DoLl,\quad X\equiv (\tilde{P},Q)\mapsto 
\l\{ \begin{array}{ll}\PsoLl^-\circ\PsolN(X) & ,\LA \tilde{P},Q\RA \leq 0\\
\PsoLl^+\circ\PsolN(X) & ,\LA \tilde{P},Q\RA > 0\end{array}  \ri.\] 
onto its image which conjugates the maps
($\Psol= (\Xio_l)^{-1}\circ\PsoLl\circ\Xio_l$)
is $C^0$-near to the identity in the sense that 
\beq
|\Xio_l(X)-X| = \cO(1/\cE(X)),\qquad (X\in\pa\Dol),
\Leq{est:xio}
and the solution $\PsoLl$ of the linear problem is $C^1$-near to 
$\Psol$ in the sense
\beqn
|\Psol(X)-\PsoLl\circ\Xio_l(X)| &=& \cO(1/\cE(X))
\label{Po:one}\\
 \l\|D\Psol(X) - D\PsoLl\circ\Xio_l(X) \ri\| &=& \cO(1/e(X)).
 \label{Po:two}
\eeqn
\end{proposition}
\begin{remarks}
{\bf 1)} The energy-independent estimate (\ref{Po:two}) may seem to be poor
but is in fact optimal in its energy dependence since, relative
to the optimal estimates
\beq
\|D\Psol(X)\|=\cO(\cE(X)/e(X))=\|D\PsoLl(X)\|,
\Leq{o:E} 
it is of order $\cO(1/E)$
(The r.h.s.\ of (\ref{o:E}) is obtained by inserting the time bound
(\ref{L:time:o}) into the linearization
of (\ref{l:init:val})).\\[2mm]
{\bf 2)} 
For pericentric initial data $X_0\in\Hol$ the total time spent inside the ball
equals $\ToLl^+(X_0)-\ToLl^-(X_0)\equiv\pm 2\ToLl^\pm(X_0)$.
Estimate (\ref{L:time:o}) for that time 
is presented in a form needed to evaluate the
term $D\PsoLl\circ\Xio_l(X)$ in (\ref{Po:two}).

In polar coordinates $(r,\vv)$ the Kepler hyperbola 
has the parametric form (see e.g.\ \cite{Th}, Chapter 4.2)
\[r(\vv)=\frac{\vec{L}^2}{|Z|e\cos(\vv-\vv_0)+Z}.\]
The denominator has the zeros $\vv^\pm$, and $\Delta\vv:=\vv^+-\vv^-$
is the angle under which the hyperbola is seen from the origin.
Thus
\[\cos(\eh\Delta\vv) = -\frac{\sign(Z)}{e},\]
so that  $\Delta\vv\in (\pi,2\pi]$ for $Z>0$ and 
$\Delta\vv\in [0,\pi)$ for $Z<0$.

On the other hand, the total change in direction $\Delta\psi$
of the velocity vector equals
\[\Delta\psi=\sign(Z)\cdot(\Delta\vv-\pi)= 2\arcsin(1/e(X_0))\in (0,\pi].\]
Thus
(\ref{L:time:o}) can be rewritten as
\beq
\exp\l(\pm 2\ToLl^\pm(X_0)\ri) = 
\frac{4c_q\sin(\eh\Delta\psi)}{|Z_l|}\cE(X_0) +\cO(\cE^0(X_0))\qquad
(X_0\in\Hol).
\Leq{L:time:o2}
This equation will be useful for the study of orbit instability, since
$\exp\l(t\ri)$ equals the expansion of the unstable manifold of 
(\ref{l:init:val}) after time $t$.
\end{remarks}
{\bf Proof.}
We set $(\tilde{P}(t),Q(t)):=\Pol(t,X_0)$ for $X_0=(\tilde{P}_0,Q_0)\in \Dol$.
For $\Eth$ large 
the squared distance $t\mapsto |Q(t)|^2$ is a strictly convex function
of time, since 
\beqn
\eh\frac{d^2}{dt^2}|Q(t)|^2 &=& 
\frac{d}{dt}\LA\tilde{P}(t),Q(t)\RA \label{second:der:o} \\
&=& |\tilde{P}(t)|^2+|Q(t)|^2 - \LA R(Q(t)),Q(t)\RA/\cE(X_0)\NN \\
&\geq& 
|\tilde{P}(t)|^2+|Q(t)|^2 \cdot(1-\cL_1 c_q/\cE(X_0))\NN \\
&\geq&\eh\l( |\tilde{P}(t)|^2+|Q(t)|^2\ri)>0\NN
\eeqn
with Lipschitz constant
\beq
\cL_1:=\sup \l\{\l. \frac{|R(Q_1)-R(Q_2)|}{|Q_1-Q_2|}\ri| 
|Q_1|,|Q_2|\leq c_Q, Q_1\neq Q_2\ri\},
\Leq{Lip:const}
(one notes from inspection of (\ref{def:R}) that $R(0)=0$).

We can bound (\ref{second:der:o}) more precisely from below 
by using the inequality
\beq
|X|^2\equiv |\tilde{P}|^2+|Q|^2\geq 
\eh\l|  |\tilde{P}|^2-|Q|^2(1-W_l(Q^* I_3 Q)/\cE(X)) \ri| =\frac{|Z_l|}{2\cE(X)} 
\Leq{X:larger}
which follows from (\ref{en:con}) and 
is valid for $\Eth>\max_{\q\in B_l(c_q)} |W_l(\q)|$. 
 
Thus $\Tol^\pm$ and $\Tol^0$ are uniquely defined finite functions.
By transversality of $\Hol$ w.r.t.\ the flow the pericentric
time $\Tol^0$ is smooth. The hypersurface $\pa \Dol$ is transversal to the
flow, too, except at $\pa \Dol\cap \Hol$. Thus it is only there that 
we have to control smoothness of the exit times $\Tol^\pm$.
 
The maps
\[\Dol\ar\bR\times\Hol,\qquad X\mapsto (\Tol^0(X),\PsolN(X))\]
and
\[\bR\times\Hol\ar \bR^4\times\bR^4,\qquad (t,X)\mapsto \PoLl(-t,X)\]
are diffeomorphisms onto their images, since $\DoLl$ does not contain 
the (single) stationary point $0$ of the linear flow $\PoLl$. 
Hence the composition 
\beq
\rho_l:\Dol\ar\DoLlb,\quad \rho_l(X) := \PoLl(-\Tol^0(X),\PsolN(X))
\Leq{def:rho:l}
of these diffeos is a diffeomorphism onto its image. 

Thus in order to compare the flows $\Pol$ and $\PoLl$,
it suffices to compare the trajectories 
$$X(t):=(\tilde{P}(t),Q(t)):=\Pol(t,X_0)\mbox{ and }
Y(t):=(\tilde{P}_L(t),Q_L(t)):=\PoLl(t,X_0)$$
for {\em pericentric} initial conditions 
$X_0=(\tilde{P}_0,Q_0)\in \Hol$.

We partition $\Hol$ into the regions 
$\Hol^< := \{ (\tilde{P},Q)\in\Hol\mid|Q|\leq\eh c_Q\}$
and 
\[\Hol^> := \{ (\tilde{P},Q)\in\Hol\mid\eh c_Q\leq |Q|\leq c_Q\}.\]
The flow through $\Hol^<$ is uniformly transversal to 
$\pa\DoLl$, whereas the exit times are uniformly bounded on $\Hol^>$.
\\[2mm]
{\bf 1)} 
On $\Hol^>$ the second equation in (\ref{second:der:o}) together
with (\ref{en:con}) yields for 
\beq
g_1(t):=|Q(t)|^2\qquad\qquad \eh\frac{d^2}{dt^2}g_1(t) - 2g_1(t) = 
h_1(t)/\cE(X_0)
\Leq{sinh2}
with 
\[h_1(t):=Z_l- \LA R(Q(t)),Q(t)\RA)-|Q(t)|^2W_l(Q^*(t) I_3 Q(t)) 
\ \ (|t|\leq\Tol^+(X_0)),\]
$g_1(0)=|Q_0|^2$ and $g_1'(0)=0$
so that 
\[g_1(t)= |Q_0|^2+\int_0^t\int_0^s (4g_1(u)+2h_1(u)/\cE(X_0))du\,ds.\]
As $C_3:= \sup_{|Q|\leq c_Q}| Z_l- \LA R(Q),Q\RA)-|Q|^2W_l(Q^* I_3 Q)|<\infty$,
for $\Eth$ large and times $t$ between $0$ and $\Tol^\pm(X_0)$ the integrand
is bounded below by $3|Q_0|^2$ and above by $5c_Q^2$ so that
\[\frac{3}{2}|Q_0|^2 t^2\leq g_1(t)\leq \frac{5}{2}c_Q^2 t^2. \]
As $g_1(\Tol^\pm(X_0))=c_Q^2$,
\beq
\eh\sqrt{1-|Q_0|^2/c_Q^2}\leq\pm\Tol^\pm(X_0)\leq 2\sqrt{1-|Q_0|^2/c_Q^2}
\qquad \l(X_0\equiv (\tilde{P}_0,Q_0)\in \Hol^> \ri),
\Leq{est:M} 
and the same estimate holds for $\ToLl^\pm$.
Est.\ (\ref{est:M}) implies in particular the uniform bound
$|\Tol^\pm|,|\ToLl^\pm|\leq 2$ on $\Hol^>$.

The difference between these times is much smaller:
\beq
|\Tol^\pm(X_0)-\ToLl^\pm(X_0)| = \cO\l(\sqrt{1-|Q_0|^2/c_Q^2}/\cE(X_0)\ri)
\qquad (X_0\in\Hol^>).
\Leq{T:diff:g}
Namely setting
$g_2(t):=|Q(t)|^2-|Q_L(t)|^2$, so that $g_2(0)=g_2'(0)=0$ and, by (\ref{sinh2}), 
\[\eh\frac{d^2}{dt^2}g_2(t)- 2g_2(t) =\frac{h_2(t)}{\cE(X_0)}\qmbox{with}
h_2(t):=h_1(t)-Z_l,\]
we get
\[g_2(t)=\int_0^t\int_0^s (4g_2(u)+2h_2(u)/\cE(X_0))du\,ds\]
or
\[|g_2(t)|\leq \frac{C_4}{\cE(X_0)}\sinh^2(t)\]
for $C_4:=C_3+\Zmax$.

As $c_Q^2=|Q_L(\ToLl^\pm(X_0))|^2=\sinh^2(\ToLl^\pm(X_0))|Q_0|^2$
and
\[g_2(\Tol^\pm(X_0))=c_Q^2-|Q_L(\Tol^\pm(X_0))|^2=
c_Q^2-\sinh^2(\Tol^\pm(X_0))|Q_0|^2,\]
\beqno
&&\sinh(\Tol^\pm(X_0))\cdot \sqrt{1-\frac{C_4}{|Q_0|^2\cE(X_0)}}\leq
\sinh(\ToLl^\pm(X_0))\\
&& \hspace*{5cm}\leq
\sinh(\Tol^\pm(X_0))\cdot \sqrt{1+\frac{C_4}{|Q_0|^2\cE(X_0)}}.
\eeqno
In view of (\ref{est:M}) this gives (\ref{T:diff:g}). 

In turn (\ref{T:diff:g}) implies
\beqn
\lefteqn{
|\PoLl(\ToLl^\pm(X_0),X_0)-\PoLl(\Tol^\pm(X_0),X_0)|\ +}\NN\\
&&\hspace*{3cm}\l\|D\PoLl(\ToLl^\pm(X_0),X_0)-D\PoLl(\Tol^\pm(X_0),X_0) \ri\|
\NN\\
&=&
\cO\l(\sqrt{1-|Q_0|^2/c_Q^2}/\cE(X_0)\ri).\label{part1}
\eeqn

A comparison between the initial value problems (\ref{init:val}) and
(\ref{init:val2}) using a Gronwall estimate
on the uniformly bounded time interval $|t|\leq 2$ yields
\beqn
|\Pol^t(X_0)-\PoLl^t(X_0)|+ \l\| D\Pol^t(X_0)-D\PoLl^t(X_0)\ri\|
&\leq& \frac{C|t|}{\cE(X_0)}\label{part2}\\
& &\hspace{-4cm}(X_0\in \Hol^>,t\in[\Tol^-(X_0),\Tol^+(X_0)]).\NN
\eeqn
Setting $t:=\Tol^\pm(X_0)$ in (\ref{part2}), the triangle inequality and
(\ref{part1}) leads to
\beqno
\lefteqn{|\Psol^\pm(X_0)-\PsoLl^\pm(X_0)|+\l\|
D\Psol^\pm(X_0)-D\PsoLl^\pm(X_0)\ri\|
\equiv}& &\\
& &|\Pol(\Tol^\pm(X_0),X_0)-\PoLl(\ToLl^\pm(X_0),X_0)|\ +\\
& &\hspace{-1cm}\l\|D\Pol(\Tol^\pm(X_0),X_0)-D\PoLl(\ToLl^\pm(X_0),X_0) \ri\|=
\cO\l(\sqrt{1-|Q_0|^2/c_Q^2}/E\ri).
\eeqno
This does not only prove the estimate (\ref{est:xio}) on
$\Psol^+(\Hol^>)\cup \Psol^-(\Hol^>)\subset \pa\Dol$, but also shows 
that $\Xio$ is continuously differentiable at the submanifold
\[\Hol\cap\pa\Dol = \Psol^+(\Hol^>)\cap \Psol^-(\Hol^>)\]
of phase space points where the flows $\Pol$ and $\PoLl$ are tangential to
$\pa\Dol$ (for $X$ in this set $\Xio(X)=X$ and $D\Xio(X)= \idty$).
\\[2mm]
{\bf 2)}
$\bullet$ On $\Hol^<$ we begin with a rough estimate. For time 
\beq
\tau\equiv\tau(X_0):= \ln\l(3c_Q/|X_0|\ri)
\Leq{def:tau}
\beqn
\hspace{-8mm} |Q_L(\pm\tau)| &=& |\cosh(\tau)Q_0+\sinh(\tau)\tilde{P}_0|
=\sqrt{\sinh^2(\tau)|\tilde{P}_0|^2 + \cosh^2(\tau)|Q_0|^2}\NN\\
&\geq& \sqrt{\eh |X_0|^2\cdot(\cosh(2\tau)-1)}
\geq \sqrt{\eh |X_0|^2\cdot(\eh\exp(2\tau)-1)}\NN\\
&=&\sqrt{{\textstyle \frac{9}{4}}c_q-\eh|X_0|^2}\geq \frac{\sqrt{7}}{2}c_Q
\label{Q:tau}
\eeqn
since by (\ref{en:con})   
\[|X_0|^2 = 2|Q_0|^2+ (|\tilde{P}_0|^2 - |Q_0|^2) \leq\eh c_q
+\frac{Z_l-|Q_0|^2 W_l(Q_0^*I_3 Q_0)) }{E}\leq c_q.\]
Therefore 
any trajectory of the {\em linear} flow $\PoLl$ 
with these initial conditions leaves the region $\Dol$ before time $\tau$
(and enters it after time $-\tau$).\\[2mm] 
$\bullet$ We now do perturbation theory around these linear solutions, setting
\[Z(t) := X(t)-Y(t) .\]
$Z(t)$ meets the integral equation 
\[Z(t) = \int_0^t \l(
  \bem{cc}0 & \idty\\ \idty& 0\eem Z(s)
  + E^{-1}\tilde{R}(X(s))\ri) ds\]
with $E:=\cE(X_0)$.
Gronwall's inequality says that
\beq
f(t)\leq A\exp\l(\int_0^t g(s)ds\ri) 
\Leq{gron}
if  $f(t)\leq A+\int_0^t g(s)f(s)ds $. 
Applied to $f(t):= |Z(t)|$, we get for $0\leq t\leq \Tol^+(X_0)$
\beqn
f(t)&\leq& \int_0^t \l( f(s) + 
E^{-1}|\tilde{R}(X(s))-\tilde{R}(Y(s))+\tilde{R}(Y(s))|\ri) ds
\label{A:est}\\
&\leq& E^{-1}\int_0^t |\tilde{R}(Y(s))|ds +
\int_0^t \l(1+\frac{\cL_1}{E}\ri) f(s) ds \NN\\
&\leq& E^{-1}\cL_1\int_0^t |Q_L(s)|ds +
\int_0^t \l(1+\frac{\cL_1}{E}\ri) f(s) ds \NN
\eeqn
with the Lipschitz constant $\cL_1$ from (\ref{Lip:const}). But
\beqno
\int_0^t |Q_L(s)|ds 
&=& \int_0^t |\sinh(s)\tilde{P}_0+\cosh(s)Q_0|ds\NN\\
&\leq& |X_0|\int_0^t e^s ds \leq |X_0|e^t.\NN
\eeqno
Thus the constants in (\ref{gron}) can be chosen as 
\[A\equiv A(X_0,t):=\frac{\cL_1|X_0|e^{|t|}}{\cE(X_0)}\qmbox{,}
g\equiv g(X_0):= 1+\cL_1/\cE(X_0),\]
and we obtain for $|t|\leq \tau(X_0)$ with $X(t)\in \Dol$
\beqn
|Z(t)| &\equiv f(t)& \leq \cL_1|X_0|e^{|t|} \frac{\exp((1+\cL_1/E)\cdot\tau(X_0))}{E}\NN\\
&=& |X_0|e^{|t|}\cdot \cO(E^{-\eh}) = 3c_Q\cdot \cO(E^{-\eh})
\label{ineq:f}
\eeqn
for $\Eth$ large since by (\ref{def:tau}) and (\ref{X:larger}) 
\[\tau(X_0) = \ln\l(3c_Q/|X_0|\ri)\leq
\ln\l(3\sqrt{2}c_Q\sqrt{E/|Z_l|}\ri) = \ln( c E^\eh).\]
In particular we conclude from (\ref{Q:tau}) that
\[0\leq \pm \Tol^\pm(X_0)\leq \tau(X_0),\]
since otherwise by (\ref{ineq:f})
\[ |Q(\pm\tau)| \geq |Q_L(\pm\tau)| - |Z(\pm\tau)| \geq 
\frac{\sqrt{7}}{2}c_Q-3c_Q\cdot \cO \l(E^{-\eh}\ri) > c_Q.\]

We use (\ref{ineq:f}) as an input for a refined estimate
which will imply $|Z(t)|=\cO(1/E)$. To that end we note that by (\ref{ineq:f})
\[ |Q(t)| \leq |Q_L(t)| + |Z(t)| \leq 
e^{|t|}|X_0|+|Z(t)|\leq 2e^{|t|}|X_0|\quad (|t|\leq \tau(X_0)).\]
We write
\beq
Z(t) = E^{-1}
\int_0^t \exp\l(\bsm0& \idty\\ \idty& 0\esm(t-s)\ri)\tilde{R}(X(s))\, ds.
\Leq{Z:inte}
W.l.o.g.\ we may assume that $W_l(\s_l)=0$, since otherwise we may shift
the energy $E$ by that constant, producing an error term of relative order 
$\cO(1/E)$. Then instead of the Lipschitz estimate (\ref{Lip:const}) 
for $R$ we use
\[R(Q)\leq \cL_2|Q|^2 \qquad (|Q|\leq c_Q).\]
Inserting this into (\ref{Z:inte})
we get for $0\leq t\leq \tau(X_0)$, using (\ref{def:tau}), 
\[|Z(t)| \leq \frac{4\cL_2}{E}\int_0^t \exp(t-s)e^{2s}|X_0|^2\, ds\leq 
\frac{12 \cL_2 c_Q}{E}(e^t-1)|X_0|\leq 
\frac{24 \cL_2}{E}|Q_L(t)|, \]
since 
\[|Q_L(t)|=\sqrt{\sinh^2(t)|\tilde{P}_0|^2+\cosh^2(t)|Q_0|^2}\geq
\sinh(t)|X_0|.\] 
A similar estimate holds for $0\geq t\geq -\tau(X_0)$.

But $|Q(t)|\geq|Q_L(t)|-|Z(t)|$, so that 
\[|Z(t)|\leq \l( \frac{E}{24\cL_2}-1 \ri)^{-1}|Q(t)|,\]
showing that the diffeomorphism $\rho_l$ of (\ref{def:rho:l})
onto its image is $C^0$--near to the identity in the sense 
\beq
\frac{|\rho_l(X)-X|}{|Q|}=\cO(1/\cE(X))\qquad (X=(\tilde{P},Q)\in\Dol).
\Leq{C0:near}
If we assume $X\in\pa\Dol$ so that $|Q|=c_Q$, then (\ref{C0:near})
shows that 
\beq
|Q_L|-c_Q=\cO(1/\cE(X))
\qmbox{for} (\tilde{P}_L,Q_L):=\rho_l(X).
\Leq{near:pa}
 Since $X=\Psol^\pm(X_0)$ with $X_0\in \Hol^<$,
and $\pa\DoLl$ is uniformly transversal to the flow $\PoLl$ through $\Hol^<$    
we obtain from (\ref{near:pa})
\beq
|\Tol^\pm(X_0)-\ToLl^\pm(X_0)| = \cO(1/\cE(X_0))\qquad (X_0\in\Hol^<).
\Leq{T:diff:s}
In turn, this and (\ref{C0:near}) imply (\ref{est:xio}).\\[2mm]
$\bullet$ 
We also get (\ref{time:int}) from the time estimates (\ref{T:diff:g}) and 
(\ref{T:diff:s}) in the two regions.\\[2mm]
$\bullet$ 
Estimate (\ref{L:time:o}) for the exit time of the linear flow is derived as follows.
By (\ref{En:Lo}) the `linear' energy parameter
equals 
\[\cE_L(X_0)=\frac{Z_l}{|\tilde{P}_0|^2-|Q_0|^2}.\]
The (lifted) angular momentum (see (\ref{lift:L}) 
of the pericentric initial data $X_0=(\tilde{P}_0,Q_0)$ equals
\[|\vec{L}(X_0)| = \ev|Q_0^*P_0-P_0^*Q_0| =\eh |{\rm Im}(Q_0^*P_0)|=
\eh |P_0|\,|Q_0| = \sqrt{2\cE_L(X_0)} |\tilde{P}_0|\,|Q_0|,\]
since $\LA P_0, Q_0\RA=0$, see (\ref{lift:L}).
Thus
\[\frac{|Z_l|\cdot e(X_0)}{\cE_L(X_0)}=
\sqrt{(Z_l/\cE_L(X_0))^2+4|\tilde{P}_0|^2|Q_0|^2}
=|\tilde{P}_0|^2+|Q_0|^2.\]
On the other hand the exit times $\ToLl^+(X_0)$ and $\ToLl^-(X_0)=-\ToLl^+(X_0)$ 
are implicitly given by the equation
\[\l|Q_L(\ToLl^+(X_0),X_0)\ri|=c_Q\]
with $Q_L(t,X_0)=\sinh(t)\tilde{P}_0+\cosh(t)Q_0$, whence
\beqno
\exp(2\ToLl^+(X_0)) &=& \frac{4c^2_Q}{|\tilde{P}_0|^2+|Q_0|^2} +\cO(\cE_L^0(X_0))\\
&=& \frac{4c_q \cE_L(X_0)}{|Z_l|\, e(X_0)}+\cO(\cE_L^0(X_0)).
\eeqno
This proves (\ref{L:time:o}).
\\[2mm]
$\bullet$ The $C^0$-estimate (\ref{Po:one}) follows immediately from 
(\ref{est:xio}).\\[2mm]
$\bullet$ In order to obtain the $C^1$-estimate (\ref{Po:two}), one considers
\[DZ(t):=D\Pol^t(X)-D\PoLl^t(\Xio(X))\qquad(X_0\in\Psol^\pm(\Hol^<)).\]
$DZ(t)$ solves the integral equation
\[DZ(t) = \int_0^t \l(
  \bem{cc}0 & \idty\\ \idty& 0\eem DZ(s)
  + E^{-1}D\tilde{R}(X(s)) DX(s)\ri) ds.\]
W.l.o.g.\ we may again assume that $W_l(\s_l)=0$, since otherwise we may shift
the energy $E$ by that constant, producing an error term of relative order 
$\cO(1/E)$. Then by inspection of (\ref{def:R}) one notes that
the matrix $D\tilde{R}((\tilde{P},Q))=0$ at $Q=0$. We may thus estimate
\beq
\|D\tilde{R}((\tilde{P},Q))\|\leq \cL_3|Q| \qquad ((\tilde{P},Q)\in \Dol),
\Leq{Lip}
where $\cL_3>0$ is the Lipschitz constant of $DR$ for $|Q|\leq c_Q$.

From (\ref{Lip}) we obtain (\ref{Po:two}) by a Gronwall estimate similar to 
(\ref{ineq:f}), which we apply to
\[DZ(t) = \int_0^t\l[ \l(
  \bsm0 & \idty\\ \idty& 0\esm + E^{-1}D\tilde{R}(X(s))\ri)DZ(s)
  + E^{-1}D\tilde{R}(X(s)) DY(s)\ri] ds,\]
with $t\leq \Tol^+(X_0)$,
$f(t):=\|DZ(t) \|$, 
\[g(s):=1+\cL_3c_Q/E\geq\|\bsm0 & \idty\\ \idty& 0\esm + E^{-1}D\tilde{R}(X(s))\|\]
and 
\[A:=\frac{4c_q\cL_3}{\sqrt{|Z_l|}}\sqrt{E/e}\geq 
E^{-1}\l\|\int_0^{\Tol^+(X_0)} D\tilde{R}(X(s)) DY(s) ds \ri\|.\]
To obtain the last estimate, the time bounds
(\ref{time:int}) and (\ref{L:time:o}) are inserted. Thus for $t\leq \Tol^+(X_0)$
\[f(t)\leq A\exp\l(\int_0^t g(s)ds\ri)=
\cO\l(E^{-\eh}e(X_0)^{-\eh}\ri)\cdot\cO\l(E^{\eh}e(X_0)^{-\eh}\ri)
.\hspace*{5mm} \Box\]
%
\Section{Regularization by Phase Space Extension} \label{sect:RL}
%
In Sect.\ \ref{sect:KS:Appl} the motion in configuration space 
$\bR^3_\q$ near a singularity at $\s_l$
was regularized using the KS transform. This will enable us in Sect.\
\ref{sect:si:sc} to compare that motion with the motion in the 
Kepler potential $\q\mapsto -Z_l/|\q-\s_l|$.

However, we would like to apply these local estimates
to a complete Hamiltonian flow on a phase space which arises by
a completion of $T^*\Muh$.
Therefore in the case of attracting singularities ($Z_l>0$) 
we now employ a different regularization.

To preserve continuity of the motion with respect to the initial conditions,
a particle colliding with a nucleus at $\s_{l}\in\Mu$
must be reflected backwards. Then we parametrize the state of the 
colliding particle by its 
energy and by its incoming (or outgoing) direction. 
That is, we complete phase space by adjoining manifolds
$\bR\times S^{2}$, one for each attracting singularity.

Note that for all energies $E$ 
the energy surfaces $\Hh^{-1}(E)$ could be completed topologically
using only the Kustaanheimo-Stiefel construction (by taking the quotient of  
the quadric surface $\cI^{-1}(0)$ w.r.t.\ the circle action $\alpha$
defined in (\ref{circle:action})). This, however would 
lead to a time change which is unwanted here.
\begin{theorem} \label{P:omega} %
There exists a unique smooth extension $(P,\omega,H)$
of the 
Ha\-miltonian system $(T^{*}\Muh,\hat{\omega},\Hh)$,
where the {\bf phase space} \label{phase:sp}
$P$ is a smooth six-dimensi\-onal manifold with
\[P := T^{*}\Muh\cup
\bigcup_{\stackrel{1\leq l\leq n}{Z_l>0}} \l(\bR \times S^{2}\ri)\]
as a set, 
$\omega$ is a smooth symplectic two-form on $P$ with
\[\omega\rstr_{T^{*}\Muh} = \hat{\omega}:= \sum_{i=1}^3 dq_{i}\wedge dp_{i},\]
and $H:P\ar \bR$ is a smooth Hamiltonian function with
$H\rstr_{T^{*}\Muh} = \Hh$.

The smooth Hamiltonian flow 
\beq
\Phi:\bR\times P\ar P 
\Leq{eq:Pt}
generated by $H$ is complete (and we often write $\Phi^t(x)$ instead of
$\Phi(t,x)$).

For all energies $E > \Vmax$ (defined in (\ref{def:Vmax})) the 
{\em energy shell} 
\beq
\SuE := \{ x\in P\mid H(x) = E \}
\Leq{def:SuE}
is a smooth, five-dimensional manifold. 
\end{theorem}
{\bf Proof.} 
It is clear that the particle cannot escape to 
spatial infinity in finite time.
If $Z_l<0$ there is no need here to regularize, 
since then by (\ref{cq:small})
the minimal distance is bounded by
$|\q-\s_l|\geq|Z_l|/(2E)$ if $E>0$ resp.\ $|\q-\s_l|\geq c_q$
if $E\leq0$.

For the remaining case of an attracting singularity ($Z_l>0$) 
we linearize the motion near
collision by using as phase space coordinates the 
angular momentum components, the direction of the
Runge-Lenz vector $\vec{F}_{l}$, 
energy and the time passed since the pericentre of the
orbit. The first five of these six functions are
constant on the Kepler orbit. Then we add the collision manifold 
of phase space points with
time and angular momentum both equal to zero. This manifold is then
parametrized  by energy and by $\vec{F}_{l}/|\vec{F}_{l}|$ 
and is thus diffeomorphic to $\bR \times S^{2}$.

In \cite{KK} the case of $d=2$ dimensions is treated,
the calculations being more detailed than here.
Comparing with the Delaunay coordinates
(see, e.g., \cite{AM}, Chapter 9.3), instead of the semi-major 
axis we use energy, and instead of the mean anomaly (which is only defined
for negative energies) we use time.

If the potential $V$ is not centrally symmetric around $\s_l$, then 
at collision some of the
former constants of motion cease to be smooth functions of time.
To remedy this is we redefine them by using their value at the pericentre
of the orbit. \\
{\bf 1)}
More specifically, we introduce adapted coordinates to regularize the flow
in the phase space neighbourhood $\Uhle$, $0<\vep\leq c_q$, 
of the $l$th nucleus, with
\beq
\Uhle := \l\{\pq\in T^{*}\Muh \l| \B{\q-\s_{l}}<\vep,\,\,
|\p|^{2}>\frac{3}{2}\frac{Z_{l}}{\B{\q-\s_{l}}} \ri.\ri\}.
\Leq{U:l}
On $\Uhle$, the Hamiltonian function has the form
$\Hh\pq = \Hh_{l}\pq + W_{l}(\q)$ with
\beq
\Hh_{l}\pq := \eh \p^{\,2} - \frac{Z_{l}}{\B{\q-\s_{l}}}
\Leq{pure:kepler}
and the {\em smooth} additional potential 
$W_{l}$ on $B_l(c_q)$,
\[W_{l}(\q) = \sum_{i\neq l} \frac{-Z_{i}}{\B{\q-\s_{i}}} + W(\q).\]
One basic estimate on $\Uhle$, valid for $\vep$ small, is
\beqn
\frac{d}{dt}((\q-\s_l)\cdot\p)
&\geq&
\eh\frac{Z_l}{|\q-\s_l|}-(\q-\s_l)\cdot\nabla W_l(\q)>0.
\label{elapse}
\eeqn
Every collision orbit with $\s_l$ enters $\Uhle$, as
\[|\p|^2-\frac{3}{2}\frac{Z_l}{\B{\q-\s_l}}=
\eh\frac{Z_l}{\B{\q-\s_l}}+2(E-W_l(\q))\ar\infty\]
as $\q$ approaches $\s_l$.
\\[2mm]
{\bf 2)} 
We first treat the Keplerian case $W_l\equiv 0$.
The angular momentum $\hat{L}_{l}:\Uhle\ar \bR^3$ relative to the position 
of the $l$th nucleus equals
\beq
\hat{L}_{l}\pq := (\qs)\times\p . 
\Leq{def:hat:Ll}

Let $\hat{T}_{l}:\Uhle\ar\bR$ be the time elapsed since the closest 
encounter of the Kepler
solution with the nucleus. By (\ref{elapse}) 
there is only one such {\em pericentre} of the
orbit, with distance $\rmin$.
$\hat{T}_{l}$ is given by
\beq
\hat{T}_{l}\pq := \int^{\B{\qs}}_{\rmin\pq}
\frac{r\,dr}{\sqrt{2r^{2}\Hh_{l}\pq + 2Z_{l}r - \hat{L}_{l}^{2}\pq}}\cdot
\sign((\q-\s_{l})\cdot\p)
\Leq{peric:time}
with
\beq
\rmin\pq := 
\l\{ \begin{array}{cc}
     \frac{-Z_{l}+\sqrt{Z_{l}^{2}+2\Hh_{l}\pq\hat{L}_{l}^{2}\pq}}{2\Hh_{l}\pq}\!\! 
         &, \Hh_{l} \neq 0 \\
     \hat{L}_{l}^{2}\pq/2Z_{l} \!\! & , \Hh_{l} = 0
     \end{array} \ri. .
\Leq{rmin}
$\hat{T}_{l}$ is a smooth function, which can be seen by 
explicit evaluation of 
the integrals:
\begin{eqnarray}
\lefteqn{\int \frac {r} {\sqrt{2r^{2}E + 2Z r - L^{2}}} dr = }
\label{eq:explicit:int}\\ 
& &\hspace{-14mm}
\frac{r}{\sqrt{2E}} \sqrt{1 +\frac{Z}{rE} - \frac{L^{2}}{2r^{2}E}} -
              \frac{Z}{(2E)^{3/2}}
              \ln \l( Er + \eh Z + \sqrt{ E(r^{2}E + Zr - \eh L^{2}) }\ri)
\nonumber
\end{eqnarray}
for $E>0$ and $Z>0$ (see Thirring \cite{Th}, for more information).

The Runge-Lenz vector $\vec{F}_{l}: \Uhle\ar \bR^3$
relative to the $l$th centre is given by
\beq
\vec{F}_{l}\pq :=  \p\times\hat{L}_{l}\pq - Z_{l}\frac{\qs}{\B{\qs}}.
\Leq{def:vv:h:l}
On its domain $\Uhle$ of definition $\vec{F}_{l}$ is non-zero:
$ |\vec{F}_{l}|^{2} = 
2|\hat{L}_{l}|^{2}\Hh_{l} + Z_{l}^{2}>Z_{l}^{2}/4$.
Thus we may define the {\em pericentral direction} 
$\hat{F}_l:\Uhle\ar S^2$ by $\hat{F}_l:=\vec{F}_{l}/|\vec{F}_{l}|$.
The angular momentum vector $\hat{L}_l$ is perpendicular to that direction:
$\hat{L}_l\cdot \hat{F}_l = 0$.
The map 
\[\hat{{\cal Y}}:\Uhle\ar T^*(\bR\times S^2)\setminus \bar{0},\quad
\pq\mapsto (\hat{T}_{l},\hat{L}_l; \Hh_{l},\hat{F}_l)\]
is a diffeomorphism onto its image, ($\bar{0}$ denoting the zero section
of the cotangent bundle $T^*(\bR\times S^2)$).\\[2mm]
{\bf 3)} 
The Poisson brackets between the above variables are given by
\[\{\hat{T}_{l},\Hh_{l}\}=1,\mbox{ , }\{\hat{L}_{l},\Hh_{l}\}=\vec{0}
\mbox{ , }\{\hat{F}_{l},\Hh_{l}\}=\vec{0}
\mbox{ , }\{\hat{F}_{l},\hat{T}_{l}\}=\vec{0}
\mbox{ and } \{\hat{L}_{l},\hat{T}_{l}\}=\vec{0}
,\] 
so that in particular
the angular momentum and the pericentral direction are constants of the
Kepler flow. 
The components of angular momentum and the asymptotic direction
have the Poisson brackets
$$ \{(\hat{F}_{l})_i,(\hat{F}_{l})_j\}=0
,\quad 
\{(\hat{L}_{l})_i,(\hat{L}_{l})_j\}=\vep_{ijk}(\hat{L}_{l})_k
,\quad 
\{(\hat{L}_{l})_i,(\hat{F}_{l})_j\}=\vep_{ijk}(\hat{F}_{l})_k,$$
using the Poisson brackets 
$\{(\vec{F}_{l})_i,(\vec{F}_{l})_j\}=-2\Hh_{l}\vep_{ijk}(\vec{F}_{l})_k$.\\[2mm]
{\bf 4)} 
Because of the above Poisson brackets with $\Hh_{l}$, 
by introducing the above coordinates, we obtain a
chart in $\Uhle$ which explicitly linearizes the (incomplete) Kepler flow.

The motion is then regularized in the following way. One defines a
completion of $\Uhle$ by setting $\Ule := \Uhle\cup \l(\bR \times
S^{2}\ri)$ as a set, and one introduces a topology on $\Ule$ by
extending the map
$\hat{{\cal Y}}=(\hat{T}_{l},\hat{L}_{l};\Hh_{l},\hat{F}_{l})$ to  
\[ {\cal Y} := (T_{l},{\bf L}_{l};H_{l},{\bf F}_{l}):
\Ule \ar T^*(\bR\times S^{2})\]
by mapping $(h,{\bf f}) \in \bR\times S^{2}$
onto the point $(0,0;h,{\bf f})$ of the zero section. 
By that procedure we obtain the
topological manifold $P$ and, by taking limits, we extend the
Hamiltonian $\Hh$ to a continuous function $H:P\ar \bR$. 
The topology of $P$ is thus determined by
the purely Coulombic local Hamiltonians $\Hh_{l}$, and, by taking
limits for the collision orbits, we are able to extend the flow 
generated by $\Hh$ to a
complete continuous flow $\Pt$ on $P$. 
 
Moreover, the calculation of all the Poisson brackets shows that we
may continuously extend the symplectic form 
$\l(\sum_{i=1}^3 dq_{i}\wedge dp_{i}\ri)\rstr_\Uhle$ to $\Ule$ and obtain 
a nondegenerate two-form, which is smooth in the new coordinates.\\[2mm]
{\bf 5)} 
To generalize the construction to the case of the flow generated by
$H$ which is of the local form $H_{l}+W_{l}$, 
we define similar canonical coordinates in $\Ule$
which linearize the $H$-flow $\Pt$.

By (\ref{elapse}) for $\vep>0$ small enough, 
the orbits are
transversal to the hypersurface $\Shle$.
We extend that hypersurface to the topological submanifold
\beq
\Sle := \Shle \cup \l(\bR\times S^{2}\ri) \subset \Ule.
\Leq{def:Sle}

If $W_{l}\equiv 0$, we define the differential structure near the $l$th
collision manifold $\bR\times S^2$ by pull-back with ${\cal Y}$.
In general we define a map
\[\tilde{{\cal Y}} := (\tilde{T}_{l},\tilde{L}_{l};\tilde{H}_{l},\tilde{F}_{l}):
\Ule\ar T^*(\bR\times S^{2})\]
by letting $\tilde{H}_{l} := H\rstr_{\Ule}$ be the energy and 
\beq
\tilde{T}_{l}(\Pt(x)) := t\qquad (x\in \Sle)
\Leq{as:time:passes:by}
the time passed since the passage of the pericentre. 
Note that (\ref{as:time:passes:by}) defines $\tilde{T}_{l}$ everywhere
on $\Ule$, since every orbit in $\Ule$
passes $\Sle$ exactly once.
\[\tilde{L}_{l}(\Pt(x)) := {\bf L}_{l}(x)\qmbox{and}
\tilde{F}_{l}(\Pt(x)) := {\bf F}_{l}(x)\qquad
(x\in\Sle)\]
are then the angular momentum and the asymptotic direction
at the pericentre $x$ of the orbit.

Clearly, {\em by fiat}, $\tilde{H}_{l},\tilde{L}_{l}$ and $\tilde{F}_{l}$ are
constant on one orbit $\Pt\pq$, whereas
$\tilde{T}_{l}(\Pt\pq)$ is an affine function of time $t$.
Therefore we have linearized the full motion. We must show that the
functions 
$(\tilde{T}_{l},\tilde{L}_{l};\tilde{H}_{l},\tilde{F}_{l})$,
restricted to $\Uhle$, are indeed
smooth coordinates. 
But this follows from the smoothness of the lifted functions in
Lemma \ref{lem:lifted:fcts}, smoothness of time change, and smoothness of the
KS flow (\ref{init:val}) (which of course exists for all energies $E\in\bR$).

Since $\tilde{{\cal Y}}\equiv 
(\tilde{T}_{l},\tilde{L}_{l};\tilde{H}_{l},\tilde{F}_{l})$
defines a homeomorphism of $\Ule$ onto its image, we use
it to define a differential structure on the whole of
$\Ule$ and thus on $P$.

Smoothness of the energy shells $\SuE$ for $E > \Vmax$ 
follows by noticing that these $E$ are regular values of $H:P\ar
\bR$. \hfill $\Box$
%
%
\Section{M\o ller Transformations} \label{sect:moeller}
%
Next we define the M\o ller and scattering transformations
which compare the asymptotics of the motion with a `free motion'. 
We base ourselves on the articles  \cite{Hu} of Hunziker and \cite{Sim} 
of Simon.
The recent monograph \cite{DG} by Derezi\'{n}ski and G\'{e}rard 
treats these questions in the context of classical and quantum mechanical
$n$-body scattering.

Due to the long--range character of the Coulomb interaction we cannot
in general use the flow generated by the 
Hamiltonian function $\p^{\,2}/2$ as `free motion'.
Instead, we compare with the Kepler motion generated by 
\beq
\Hhi\pq := \eh \p^{\,2} + V_\infty(\q)\qmbox{with}
V_\infty(\q)= - \frac{\Zi}{\B{\q}}.
\Leq{free:hamiltonian}
Thus we consider the smooth complete flow 
\beq
\Pit:\Pin\ar \Pin
\Leq{eq:Pit}
generated by (\ref{free:hamiltonian}).
\begin{itemize}
\item
If $\Zi=0$, $\Pit$ is the free flow on the phase space 
$\Pin:=T^{*}\bR^{3}$.
\item
If $\Zi>0$, we regularize
$T^{*}(\bR^{3}\setminus\{0\})$ in the way described in Thm.\ \ref{P:omega} to obtain $\Pin$. 
\item
If $\Zi<0$, then $\Pit$ is already
complete on $\Pin := T^{*}(\bR^{3}\setminus\{0\})$, since particles
of finite energy cannot meet the origin at $\q=\vec{0}$.
\end{itemize}
Thus we are to compare motions $\Pt$ and $\Pit$
on the {\em different} phase spaces $P$ and $\Pin$. We cannot just
identify $P$ with $\Pin$ by neglecting the measure zero sets projecting
to the singularities, since later on we will be interested in 
certain sets of measure zero like the bound states 
(moreover, these bound states will turn
out to be crucial in our analysis of scattering, too). 

We overcome the above difficulty by observing that it suffices to
identify $\Pin$ with $P$ in a neighbourhood of spatial infinity.
More precisely, let 
\beq
\Pinp := \{ x\in \Pin\mid \Hi(x)>0\}
\Leq{def:Pinp}
be the set of phase space points with positive `free' energy.
Then the orbit $\Pit(x)$ starting at $x\in \Pinp$ goes to
spatial infinity for large positive and negative times. The ball 
$\{ \q\in\Mu \mid \B{\q}\leq \Rmin\}$ contains all singularities
of $V$ and the singularity at the origin of the Kepler Hamiltonian
$\Hi$. Therefore, we can canonically identify points 
$\pq\in \Pin$ with points $\pq\in P$ if $\B{\q}>\Rmin$, and we
denote this identification by $\Id$. Thus the 
{\em M\o ller transformations}
\beq
\Opm := \lim_{t\ar\pm\infty} \Phi^{-t}\circ\Id\circ\Pit
\Leq{Moeller}
are formally maps $\Opm :\Pinp\ar P$. They exist as
pointwise limits, see Thm.~\ref{thm:both:moeller} below.

First some standard definitions (see \cite{Hu}):
\begin{definition} \label{defi:bound:scattering} %
{\rm
\begin{eqnarray*}
b^{\pm} &:=& \{x\in P\mid \q\,(\pm\bR^+,x)\mbox{ is bounded }\}
\hspace{5mm}\qmbox{,}b^\pm_E:=b^\pm\cap \SuE\\
b &:=& b^{+}\cap b^{-} \quad\qmbox{(the {\em bound
states})}\hspace{7mm}\qmbox{,}
\buE:=b\cap \SuE\\
s^{\pm} &:=& \{ x\in P \mid x\not\in b^{\pm}\mbox{ and } H(x)>0\}
\hspace{8.5mm}\qmbox{,} s^\pm_E:=s^\pm\cap \SuE \\
s &:=& s^{+}\cap s^{-} \quad\qmbox{(the {\em scattering states})}\qmbox{,}
s_E:=s\cap \SuE.
\end{eqnarray*}  }
\end{definition}
We shall show that $s^{\pm} = \Opm(\Pinp)$ so that the 
term `scattering states' is really justified. 
\begin{remarks}
{\bf 1)}
By continuity of $\Pt$, $b^\pm$ can be represented as the union
$b^\pm=\cup_{k=1}^\infty b^{\pm,k}$ of compact sets $b^{\pm,k}$. Hence 
$b^{\pm}$ and $s^{\pm}$ are measurable w.r.t.\ Liouville measure
\[\Liou := \frac{1}{3!} \omega\wedge\omega\wedge\omega\]
on the symplectic manifold $(P,\omega)$.\\
{\bf 2)}
The sets $b^\pm$ and $s^{\pm}$ are also $\Pt$-invariant.

For all $E>0$ the sets $\buE^\pm$ are closed, and $\buE$ is compact, 
being a subset of the compact region in $\SuE$ projecting
to the ball $\{\q\in\bR^3\mid |\q|\leq \Rvir(E)\}$ in configuration space.
\end{remarks}

\begin{theorem} \label{thm:both:moeller} %
The limits 
\[\Opm = \lim_{t\ar\pm\infty} \Phi^{-t}\circ\Id\circ\Pit\]
exist pointwisely on $\Pinp\subset \Pin$ and thus define the
M\o ller transformations $\Opm:\Pinp\ar s^{\pm}$. These
are measure-preserving homeomorphisms and intertwine $\Pt$ and
$\Pit$:
\[\Opm\circ\Pit=\Pt\circ\Opm.\]
The asymptotic limits 
$\p^{\pm}:s^{\pm}\ar\bR^{3}$ and $\vec{L}^{\pm}:s^{\pm}\ar\bR^3$ of 
the momentum and the angular momentum $\vec{L}\pq := \q\times \p$
\beq
\p^{\pm}(x_0) := \lim_{t\ar\pm\infty}\p\circ\Phi^t(x_0) \qmbox{and}  
\vec{L}^{\pm}(x_0) := \lim_{t\ar\pm\infty}\vec{L}\circ\Phi^t(x_0)
\Leq{asympto}
are continuous functions. If $\pqs\equiv x_0$ with 
$q_0:=|\q_0|>\Rvir(E)$ and $\pm\LA\q_0,\p_0\RA \geq 0$,
then
\beq
\p^{\pm}(x_0)=\p_0+\cO\l(1/(q_0\sqrt{E})\ri)\qmbox{,}
\vec{L}^{\pm}(x_0)=\vec{L}(x_0)+\cO\l(1/(q_0^\eps\sqrt{E})\ri),
\Leq{O:pL}
and for $(\P_0,\Q_0):=\Opm\pqs$ and $E>\Eth$
\beq
\P_0=\p_0+\cO\l(1/(q_0^{1+\eps}\sqrt{E})\ri)\qmbox{,}
\Q_0=\q_0+\cO\l(1/(q_0^\eps E)\ri).
\Leq{Moe:E}
If $\eps=1$ in Def.\ \ref{defi:coulombic} of Coulombic potentials, 
then for all energies $E>0$
\[\B{\P_0-\p_0}\leq
\frac{\sqrt{E}\Rmin}{q_0}\frac{C_2}{Eq_0}
\qmbox{,}\B{\Q_0-\q_0}\leq\Rmin\frac{C_2}{Eq_0}\]
(remark that by Def.\ (\ref{R:vir}) of $C_2$ one always has $\frac{C_2}{Eq_0}\leq1$).
\end{theorem}
{\bf Proof.}
First we show (\ref{asympto}). Since the system is reversible, we consider only
the case $t\ar +\infty$.
For initial conditions $x_0\in s_E^{+}$ there exists a time $t_{0}$ with
\[\B{\q(t_{0})}\geq \Rvir(E) \qmbox{and} \LA\q(t_{0}),\p(t_{0})\RA \geq 0\]
for $(\p(t),\q(t)):=\Phi^t(x_0)$, 
since otherwise $\B{\q(t)}$ would be uniformly bounded as $t\ar\infty$.

W.l.o.g.\ we assume $t_0=0$. Then by (\ref{q:grows})
\[\q^{\,2}(t) \geq \q^{\,2}_0 + \eh Et^2 \qquad (t\geq 0).\]
Thus by (\ref{def1:u1})
\[\B{\frac{d}{dt} \p(t)} = \B{\nabla V(\q(t))} \leq
\frac{|\Zi|+C_1\Rmin^{1-\eps}}{q_0^2+\eh E t^2}\qquad (t\geq 0), \]
the limit $\p^+(x_0)= \lim_{t\ar\infty}\p(t)$ exists, and
for this choice of $x_0$
\[\p^+(x_0)-\p_0=\cO(1/(q_0\sqrt{E})).\]
Being a locally uniform limit of the continuous functions
$x_0\mapsto \p(t,x_0)$, $\p^+$ is continuous.
Next we show that the asymptotic limit $\vec{L}^{+}(x_0)$ of the 
angular momentum exists. For $t\geq 0$ we can estimate
\beqno
\B{\frac{d}{dt}\vec{L}(\p(t),\q(t))} &=& \B{\nabla V(\q(t))\times \q(t)}\\ 
   &=&    \B{ \l(\nabla V(\q(t))-\Zi\frac{\q(t)} {\B{\q(t)}^{3}} \ri)
   \times \q(t) } \\
   &\leq& C_1\Rmin \B{\q(t)}^{-1-\eps} \leq 
          C_1 \Rmin\l(q_0^2+\eh Et^2\ri)^{-\eh(1+\eps)}
\eeqno
using (\ref{def1:u1}) and (\ref{q:grows}), 
which shows the existence of the limit $\vec{L}^{+}(x_0)$, and (\ref{O:pL}).
Continuity of $\vec{L}^{+}$ follows as above.

We now seek a Kepler hyperbola $(\P(t),\Q(t)):=\Pit(X_0)$ 
which is positive asymptotic to $(\p(t),\q(t))$
and write 
\[\vec{r}(t):=\q(t)-\Q(t).\]
Then $\vec{r}$ is a solution of the differential equation
\beq
\ddot{\vec{r}}(t)=
\Zi\frac{\q(t)-\vec{r}(t)}{\B{\q(t)-\vec{r}(t)}^3}-\nabla V(\q(t))
\qmbox{with}\lim_{t\ar\infty} \vec{r}(t)=\vec{0}.
\Leq{AWP}
Setting for $E>0$
\beq 
\cC_E := \l\{ \vec{r}\in C\l([0,\infty),\bR^{3}\ri) \l| \,
\|\vec{r}\|:=\sup_{t\geq0} \B{\vec{r}(t)} < \min(\Rmin,C_2/E) \ri.\ri\}
\Leq{C:E}			
(with $C_2= 31(1+1/\eps)\Rmin^{1-\eps} C_1$ from (\ref{R:vir})), 
by (\ref{q:grows}) and (\ref{R:vir})
\beq
({\cal F}\vec{r})(t) := \int_{t}^{\infty} ds \int_{s}^{\infty} d\tau
\l(\Zi\frac{\q(\tau)-\vec{r}(\tau)}{\B{\q(\tau)-\vec{r}(\tau)}^{3}}-
\nabla V(\q(\tau)) \ri)
\Leq{integral:eq}
is well-defined for $\vec{r}\in \cC_E$, noting that by (\ref{R:vir}) 
\[\B{\q(\tau)-\vec{r}(\tau)}\geq\Rmin.\]

We estimate $|({\cal F}\vec{r})(t)|$ as follows. 
The integrand of (\ref{integral:eq}) is bounded by
\beqno
\lefteqn{\l|\Zi\frac{\q(\tau)-\vec{r}(\tau)}{\B{\q(\tau)-\vec{r}(\tau)}^{3}}-
\nabla V(\q(\tau)) \ri|}\\
&\leq&
\l|\nabla V(\q(\tau)) - \Zi\frac{\q(\tau)}{\B{\q(\tau)}^{3}}
\ri|+
\l|\Zi\frac{\q(\tau)-\vec{r}(\tau)}{\B{\q(\tau)-\vec{r}(\tau)}^{3}}
-\Zi\frac{\q(\tau)}{\B{\q(\tau)}^{3}}\ri|\\
&\leq& \frac{C_1\Rmin}{\B{\q(\tau)}^{2+\eps}}+
\frac{C_1\Rmin}{\min(\B{\q(\tau)},\B{\q(\tau)-\vec{r}(\tau)})^{2+\eps}}
\leq \frac{9C_1\Rmin}{\B{\q(\tau)}^{2+\eps}}
\eeqno
by the decay assumptions (\ref{def1:u1}) and (\ref{def1:u2})
(which are valid for the asymptotic potential, too, see Remark
\ref{rem:symm}.2),  and 
Def.\ (\ref{C:E}). For the last inequality we used 
$\B{\vec{r}(\tau)}\leq \eh\B{\q(\tau)}$, following from (\ref{R:vir}).

Estimating $\B{\q(\tau)}$ with (\ref{q:grows}) gives
\beqn
\lefteqn{\l|\int_{s}^{\infty} d\tau
\l(\Zi\frac{\q(\tau)-\vec{r}(\tau)}{\B{\q(\tau)-\vec{r}(\tau)}^{3}}-
\nabla V(\q(\tau)) \ri)\ri|}\label{first:int}\\
&\leq &
9C_1\Rmin\int_{s}^{\infty} \l(\q_0^{\,2}+
\eh E\tau^2\ri)^{-\eh(2+\eps)}d\tau \NN\\
&\leq& \frac{9(2+\sqrt{2})C_1\Rmin}{\sqrt{E}} \max\l(\sqrt{E}s,q_0\ri)^{-1-\eps}\NN.
\eeqn
Hence a second integration yields 
\beq
\|{\cal F}\vec{r}\|\leq 9(2+\sqrt{2})C_1\Rmin\frac{1+1/\eps}{E q_0^\eps}
\leq \frac{C_2\Rmin^\eps}{E q_0^\eps},
\Leq{F:r}
and thus for $\eps=1$ or $E>C_2/\Rmin$ or $\B{\q_0}\geq\Rvir(E)$ large 
${\cal F}$ maps $\cC_E$ into itself.
Similarly, for $\vec{r}_i\in\cC_E$
\beqn
\lefteqn{\|{\cal F}\vec{r}_{1}-{\cal F}\vec{r}_{2}\|}\NN\\
&\leq & C_1
\int_{0}^{\infty} ds \int_{s}^{\infty} d\tau
\frac{\B{\vec{r}_{1}(\tau) - \vec{r}_{2}(\tau)}}
{\min(\B{\q(\tau)-\vec{r}_{1}(\tau)},
      \B{\q(\tau)-\vec{r}_{2}(\tau)})^{2+\eps}}\NN\\
&\leq & 8C_1  \|\vec{r}_{1}-\vec{r}_{2}\| 
\int_{0}^{\infty} ds \int_{s}^{\infty} d\tau \B{\q(\tau)}^{-2-\eps}\NN\\   
&\leq & 8(2+\sqrt{2})C_1 \frac{1+1/\eps}{E q_0^\eps}\|\vec{r}_{1}-\vec{r}_{2}\|
\leq\frac{8}{9}\frac{C_2\Rmin^{\eps-1}}{E q_0^\eps}
\|\vec{r}_{1}-\vec{r}_{2}\|. 
\label{contraction}
\eeqn
By (\ref{R:vir}) 
for $\eps=1$ or $E>C_2/\Rmin$ or $\B{\q_0}\geq\Rvir(E)$ large 
the r.h.s.\ of (\ref{contraction}) is bounded above by
$\frac{8}{9}\|\vec{r}_{1}-\vec{r}_{2}\|$,
so that the map ${\cal F}:\cC_E\ar\cC_E$
is a contraction, and thus has a unique fixed point 
$\vec{r}\in\cC_E$.  
The first estimate in (\ref{Moe:E}) follows from (\ref{first:int})
and the second from (\ref{F:r}).

After having shown unique existence of $(\Opm)^{-1}:s^{\pm}\ar \Pinp$,
we show unique existence of the M\o ller transforms by the same method
of integral equations, interchanging the roles of the two potentials.
This is indeed possible with the same constants
(see Remark \ref{rem:symm}.2).

The proof of the remaining statements is the same as in \cite{Sim}.
\hfill $\Box$\\[2mm]
By energy conservation and (\ref{V:ar:0}) the modulus 
asymptotic momenta $\p^{\pm}(x)$ equals $\sqrt{2H(x)}$, 
which is non-zero by Def.\ \ref{defi:bound:scattering} of $s^{\pm}$.
So the {\em asymptotic directions} 
\beq
\hat{p}^\pm: s^\pm\ar S^2\quad ,\quad 
\hat{p}^\pm(x) := \frac{\p^{\pm}(x)}{\sqrt{2H(x)}}
\Leq{as:dir}
are well-defined.

\begin{corollary} \label{coro:complete} 
\begin{enumerate}
\item 
$s^{\pm} = \Opm(\Pinp)$, that is, every positive energy orbit
which is unbounded in positive (negative) time is positively
(negatively) asymptotic to a Kepler hyperbola. 
\item The motion $\Pt$ generated by the Hamiltonian function
$H:P\ar \bR$ is {\bf asymptotically complete}, \label{asymptotically:complete}
that is, up to a subset
of Liouville measure zero the phase space consists of bound
states and scattering states: 
\[\Liou(P\setminus(b\cup s)) = 0.\]
Similarly for the Liouville measure $\lambda_E$ on the energy
shell $\SuE$
\[\Liou_E(\SuE\setminus(\buE\cup s_E)) = 0\qquad (E>0).\]
\item The {\bf scattering transformation}
\beq
S := \Op_{*}\circ \Om: D\ar \Pin
\Leq{eq:S}
with domain $D:=\Om_{*}(s)\subset \Pinp$ and range 
$\Op_{*}(s)$ is continuous and $\Pit$-invariant, {\em i.e.}\
\[S\circ \Pit=\Pit\circ S.\]
\end{enumerate} 
\end{corollary}
{\bf Proof.} \begin{enumerate}
\item 
The equality of $s^{\pm}$ with $\Opm(\Pinp)$
is a consequence of Theorem \ref{thm:both:moeller}, since $\Opm_{*} =
(\Opm)^{-1}$ is defined on $s^{\pm}$.
\item 
This is a consequence the fact that
$\Pt(b^{+,k})\subset b^{+,k}$ for $t\geq 0$ but 
$\Liou(\Pt(b^{+,k})) = \Liou(b^{+,k})$, since $\Pt$ is
canonical (see \cite{Hu}). 
\item 
Follows from Theorem \ref{thm:both:moeller}.   \hfill $\Box$
\end{enumerate}
\begin{theorem}\label{thm:smooth:moeller} %
Let $V$ be a Coulombic potential whose partial derivatives decay at
infinity according to 
\beq
\pa^{\beta}_q  \l( V(\q) + \frac{\Zi}{\B{\q}} \ri)
\stackrel{\q\ar\infty}{=} \cO\l(\B{\q}^{-|\beta|-1-\eps}\ri)
\qquad  (\beta\in \bN_{0}^3)
\Leq{smooth}
for some $0<\vep\leq 1$. Then
the M\o ller transformations $\Opm:\Pinp\ar s^{\pm}$ are $C^\infty$
diffeomorphisms and canonical transformations.

If $\pqs\equiv x_0$ with 
$q_0:=|\q_0| \geq \Rvir(E)$ and $\pm\LA\q_0,\p_0\RA \geq 0$,
then for multi-indices $\alpha,\beta\in\bN_0^3$ combined in
$\gamma:=(\alpha,\beta)$,
$\pa^\gamma_{x_0}:=\pa^\alpha_{p_0}\pa^\beta_{q_0}$
\beq
\pa^\gamma_{x_0}(\p^{\pm}(x_0)-\p_0)=
\cO\l(q_0^{-|\beta|-1}E^{-\eh(|\alpha|+1)}\ri)\ , 
\Leq{DO:p}
\beq
\pa^\gamma_{x_0}(\vec{L}^{\pm}(x_0)-\vec{L}(x_0))=
\cO\l(q_0^{-|\beta|-\eps}E^{-\eh(|\alpha|+1)})\ri),
\Leq{DO:L}
and for $(\P_0,\Q_0):=\Opm\pqs$ and $E>\Eth$
\beqn
\lefteqn{\hspace*{-10mm}
\pa^\gamma_{x_0}(\P_0-\p_0)=\cO\l(q_0^{-|\beta|-1-\vep}E^{-\eh(|\alpha|+1)}\ri)
\ , \
\pa^\gamma_{x_0}(\Q_0-\q_0)=\cO\l(q_0^{-|\beta|-\vep} E^{-|\alpha|/2-1}\ri).}
&&\NN\\
&&\label{DMoe:E}
\eeqn
\end{theorem}
\begin{remarks}\label{vep:one}
{\bf 1)}
The above decay condition (\ref{smooth}) is met by purely Coulombic potentials
(\ref{simple}) with $\eps=1$.\\
{\bf 2)}
It is more natural to rescale the momentum by setting $\vec{v}:=\p/\sqrt{2E}$
(and going to the rescaled angular momentum $\vec{L}/\sqrt{2E}$).
In terms of these variables all estimates in 
Thm.\ \ref{thm:smooth:moeller} show an energy dependence 
of the order $\cO(1/E)$.
\end{remarks}
{\bf Proof.} Note that all estimates of the theorem for $\gamma=0$ coincide 
with the ones of Thm.\ \ref{thm:both:moeller}.\\
{\bf 1)}
We estimate the derivatives of 
\[\q(t,x_0)=\q_0+t\p_0 - \int_0^t\int_0^s \nabla V(\q(\tau,x_0))d\tau \,ds.\]
We have for $g:=|\gamma|\geq 1$ 
\beqn
\lefteqn{\pa^\gamma_{x_0}
\q(t,x_0)=\pa^\gamma_{x_0}(\q_0+t\p_0) -\int_0^t\int_0^s 
\nabla \pa^\gamma_{x_0}V(\q(\tau,x_0))d\tau \,ds }&&\NN\\
&&\hspace*{-15mm}
=\pa^\gamma_{x_0}(\q_0+t\p_0) - \label{eq:schlange}\\
&&\hspace*{-15mm}\sum_{N=1}^g 
\sum_{\stackrel{\gamma^{(1)}+\ldots+\gamma^{(N)}=\gamma}{|\gamma^{(i)}|>0}}
\int_0^t\int_0^s 
D^N \nabla V(\q(\tau,x_0))
\l(\pa^{\gamma^{(1)}}_{x_0}\q(\tau,x_0),\ldots, 
\pa^{\gamma^{(N)}}_{x_0}\q(\tau,x_0)\ri)d\tau \,ds\NN
\eeqn
or
\beqn
\lefteqn{(\idty+ \cQ)(\pa^\gamma_{x_0}\q)(t) =\pa^\gamma_{x_0}(\q_0+t\p_0) - }&&
\label{oneminq}\\
&&\hspace*{-15mm}\sum_{N=2}^g 
\sum_{\stackrel{\gamma^{(1)}+\ldots+\gamma^{(N)}=\gamma}{|\gamma^{(i)}|>0}}
\int_0^t\int_0^s 
D^N \nabla V(\q(\tau,x_0))
\l(\pa^{\gamma^{(1)}}_{x_0}\q(\tau,x_0),\ldots, 
\pa^{\gamma^{(N)}}_{x_0}\q(\tau,x_0)\ri)d\tau \,ds\NN
\eeqn
with the linear operator $\cQ\equiv \cQ_{x_0}$ given by
\beq
\cQ(\vec{w})(t):=
\int_0^t\int_0^s  D\nabla V(\q(\tau,x_0))\vec{w}(\tau)d\tau \,ds
\qquad(t\geq 0).
\Leq{Q:op}
We note that on the r.h.s.\ of (\ref{oneminq}) only partial derivatives 
or order $|\gamma^{(i)}|<g$ appear, so that we can perform an 
induction in $g$, if we are able to invert the operator $\idty + \cQ$.

We assume $\vec{w}\in\hat{\cC}$ with
\beq
\hat{\cC} := \l\{ \vec{w}\in C\l([0,\infty),\bR^{3}\ri) \l| \,
\|\vec{w}\|_{\lambda}:=\sup_{t\geq0} \B{\vec{w}(t)}/\langle t\rangle_{\lambda} < 
\infty \ri.\ri\}
\Leq{hat:C}
for 
\[\langle t\rangle_{\lambda}:=\sqrt{1+(\lambda t)^2}\]
(note that $\hat{\cC}$
is independent of the choice of $\lambda>0$).

$\cQ$ maps  $\hat{\cC}$ into itself, 
and we want to prove that for all positive energies 
the operator norm of $\cQ$ is strictly smaller than one.

If we assume 
\beq
q_0\geq\Rvir(E)\qmbox{and}\LA \q_0,\p_0\RA\geq 0,
\Leq{ass:out}
then estimate (\ref{q:grows}) holds:
\beq
|\q(t)| \geq q_0\cdot \langle t\rangle_{\lambda} 
\qmbox{for all $t\geq 0$, with}\lambda:= \sqrt{E/2}/q_0.
\Leq{q:large}
(\ref{smooth}) implies that 
\beq
\pa^\beta_{q}V(\q)=\cO(|\q|^{-|\beta|-1}).
\Leq{pa:beta:V} 
So we find $c_N>0$ such that
\beq
\|D^N\nabla V(\q(\tau))\|\leq c_N (q_0\langle \tau\rangle_{\lambda})^{-N-2}
\qquad (N\in\bN).
\Leq{e1}
For estimating the norm of $\cQ$, we may restrict ourselves 
to $\vec{w}\in \hat{\cC}$ with 
\beq
\|\vec{w}\|_{\lambda} =1.
\Leq{e2}
Inserting (\ref{e1}) and (\ref{e2}) into (\ref{Q:op}), we get
\beqno\hspace*{-10mm}
|\cQ(\vec{w})(t)|&\leq& 
c_1 q_0^{-3} \int_0^t\int_0^s \langle \tau\rangle_{\lambda}^{-2}d\tau \,ds
\leq \frac{c_1 }{\sqrt{E/2} q_0^2}\int_0^t\int_0^\infty (1+u^2)^{-1}du \,ds\\
&=&\frac{c_1 \pi}{\sqrt{2E} q_0^2}\, t\leq \frac{c_1 \pi}{E q_0} \,
\langle t\rangle_{\lambda}.
\eeqno
Assuming (\ref{ass:out}), 
a suitable choice of the $E$-dependence of $\Rvir$ consistent 
with assumption (\ref{R:vir}) then gives 
\beq
\|\cQ_{x_0}\|_{\lambda}\leq \eh \qmbox{and} 
\|\cQ_{x_0}\|_{\lambda}=\cO(1/(q_0H(x_0))) 
\Leq{invert}
for the operator norm $\|\cQ_{x_0}\|_{\lambda}$ of 
$\cQ$ w.r.t.\ the norm $\|\cdot\|_{\lambda}$
in (\ref{hat:C}).

We return to (\ref{oneminq}) and estimate its r.h.s.
Using (\ref{q:large}) and (\ref{e1}), we get
\beqn\hspace*{-10mm}
\lefteqn{\l|\int_0^t\int_0^s D^N \nabla V(\q(\tau,x_0))
\l(\pa^{\gamma^{(1)}}_{x_0}\q(\tau,x_0),\ldots, 
\pa^{\gamma^{(N)}}_{x_0}\q(\tau,x_0)\ri)d\tau \,ds\ri|}&&\label{uff}\\
&\hspace*{-17mm}\leq& \hspace*{-10mm}c_N q_0^{-N-2}\int_0^t\int_0^s 
\langle \tau\rangle_{\lambda}^{-2}d\tau \,ds
\cdot
\prod_{i=1}^N\|\pa^{\gamma^{(i)}}_{x_0}\q(\cdot,x_0)\|_{\lambda} \leq
\CC \langle t\rangle_{\lambda} q_0^{-|\beta|}E^{-\eh |\alpha|-1},\NN
\eeqn
assuming 
\beq
\hspace*{-5mm}\|\pa^{\gamma'}_{x_0}\q(\cdot,x_0)\|_{\lambda}=
\cO\l(q_0^{-|\beta'|+\delta_{|\gamma'|,1}} 
        E^{-\eh |\alpha'|-1+\delta_{|\gamma'|,1}}\ri)
\qmbox{for}0<|\gamma'|<g.
\Leq{ind:ass}
Inserting (\ref{invert}) into (\ref{oneminq}), we see that (\ref{ind:ass})
holds for the start of the induction, i.e. $|\gamma'|=1$.
Also, (\ref{uff}) is consistent with (\ref{ind:ass}),
so that (\ref{ind:ass}) holds for all multi-indices
$\gamma'\in\bN_0^6$.
\\[2mm]
{\bf 2)}
Now similar to (\ref{oneminq}), the derivatives of the momentum meet the 
recursion
\beqn
\lefteqn{\pa^\gamma_{x_0}(\p(t,x_0)-\p_0) = -}&&\label{oneminp}\\
&&\hspace*{-15mm}\sum_{N=1}^g 
\sum_{\stackrel{\gamma^{(1)}+\ldots+\gamma^{(N)}=\gamma}{|\gamma^{(i)}|>0}}
\int_0^t
D^N \nabla V(\q(\tau,x_0))
\l(\pa^{\gamma^{(1)}}_{x_0}\q(\tau,x_0),\ldots, 
\pa^{\gamma^{(N)}}_{x_0}\q(\tau,x_0)\ri)d\tau. \NN
\eeqn
We insert (\ref{e1}) and (\ref{ind:ass}) into (\ref{oneminp}) and 
performing the time $t\ar\infty$ limit. Then by Lebesgue's 
Dominated Convergence Theorem $\pa^\gamma_{x_0}\p^+(x_0)$
exists, and estimate (\ref{DO:p}) is valid.  $\pa^\gamma_{x_0}\p^+$
As a locally uniform limit of continuous functions, $\pa^\gamma_{x_0}\p^+$ is 
continuous.
\\[2mm]
{\bf 3)}
Similar to the proof of (\ref{thm:both:moeller}) we use the ansatz
\beqno
\lefteqn{\pa^\gamma_{x_0}\l(\vec{L}(x(t,x_0))-\vec{L}(x_0)\ri)= 
\int_0^t\pa^\gamma_{x_0}\frac{d}{d\tau}\vec{L}(x(\tau,x_0))d\tau}&&\\ 
&=& \int_0^t \pa^\gamma_{x_0}\l(\nabla \Vsr(\q(\tau,x_0)) 
   \times \q(\tau,x_0)  \ri)d\tau \\
&=& \sum_{N=1}^g 
\sum_{\stackrel{\gamma^{(1)}+\ldots+\gamma^{(N)}=\gamma}{|\gamma^{(i)}|>0}}
\int_0^t D^N \l(\nabla \Vsr(\q(\tau,x_0)) \times \q(\tau,x_0)  \ri)\\
&&\hspace*{40mm}\l(\pa^{\gamma^{(1)}}_{x_0}\q(\tau,x_0),\ldots, 
\pa^{\gamma^{(N)}}_{x_0}\q(\tau,x_0)\ri)d\tau, 
\eeqno
in order to estimate the smoothness of the
difference between the actual and the asymptotic angular momentum.

Note that the short range potential $\Vsr=V-V_\infty$ appears in assumption 
(\ref{smooth}).
Thus for all $t\geq0$ the above expression is bounded above by 
\[c_N q_0^{-|\beta|-1-\vep} E^{-\eh|\alpha|}
\int_0^\infty \langle \tau\rangle_{\lambda}^{-1-\vep} d\tau\leq
\CC q_0^{-|\beta|-\vep} E^{-(|\alpha|+1)/2}.\]
So by dominated convergence $\pa^\gamma_{x_0}\vec{L}^+(x_0)$ exists, and 
is estimated by (\ref{DO:L}).\\[2mm]
{\bf 4)}
As in Thm.\ \ref{thm:both:moeller}
we now consider the Kepler hyperbola 
\[\l(\P(t;x_0),\Q(t;x_0)\ri) = \Pit(X_0)\qmbox{,}
X_0 \equiv (\P_0,\Q_0) =\Op(x_0)\]
which is positive asymptotic to 
$\Phi^t(x_0)=(\p(t,x_0),\q(t,x_0))$
and write 
\[\vec{r}(t)\equiv \vec{r}(t;x_0):=\q(t,x_0)-\Q(t;x_0).\]
Then $\vec{r}$ is the solution of the integral equation (\ref{integral:eq}).
Thus formally
\beqn
\pa^\gamma_{x_0}\vec{r}(t) &=& \int_{t}^{\infty} \int_{s}^{\infty} 
\ \pa^\gamma_{x_0}\nabla
\l(V_\infty(\Q(\tau;x_0)) - V(\q(\tau,x_0)) \ri)    d\tau\,  ds  \NN\\
&=& I_1^\gamma(t;x_0)+I_2^\gamma(t;x_0)
\label{D:integral:eq}
\eeqn
with the integrals
\[I_1^\gamma(t;x_0):=- \int_{t}^{\infty} \int_{s}^{\infty}
 \pa^\gamma_{x_0}\nabla \Vsr(\q(\tau,x_0))\,  d\tau\,  ds\]
and
\[I_2^\gamma(t;x_0):= \int_{t}^{\infty} \int_{s}^{\infty}
\pa^\gamma_{x_0}\nabla 
\l(V_\infty(\Q(\tau;x_0))-V_\infty(\q(\tau,x_0))\ri)\, 
d\tau\,  ds.\]
We estimate the first integral in a similar fashion as the r.h.s.\ of 
(\ref{eq:schlange}), but use the stronger estimate (\ref{smooth})
for $\Vsr$ instead of (\ref{pa:beta:V}) for $V$.
We thus obtain (with $g=|\gamma|$)
\beqno
\lefteqn{|I_1^\gamma(t;x_0)| 
\leq \sum_{N=1}^g 
\sum_{\stackrel{\gamma^{(1)}+\ldots+\gamma^{(N)}=\gamma}{|\gamma^{(i)}|>0}}}&&\\
&&\int_{t}^{\infty} \int_{s}^{\infty}
\l|D^N \nabla \Vsr(\q(\tau,x_0))
\l(\pa^{\gamma^{(1)}}_{x_0}\q(\tau,x_0),\ldots, 
\pa^{\gamma^{(N)}}_{x_0}\q(\tau,x_0)\ri)\ri|d\tau \,ds\\
&\leq&
\CC q_0^{-N-2-\vep}\int_{t}^{\infty} \int_{s}^{\infty}
\langle \tau\rangle_{\lambda}^{-2-\vep}d\tau \,ds
\cdot\prod_{i=1}^N\|\pa^{\gamma^{(i)}}_{x_0}\q(\cdot,x_0)\|_{\lambda}\\
&\leq&
\CC q_0^{-|\beta|-\vep}E^{-\eh |\alpha|-1}.
\eeqno
Unlike $I_1^\gamma$, $I_2^\gamma$ depends on $\vec{r}$. So we split it into
\[I_2^\gamma = \cP(\pa^\gamma_{x_0}\vec{r}) +I_3^\gamma+ I_4^\gamma\]
with the linear operator
\[\cP(\vec{w})(t):= -\int_{t}^{\infty} \int_{s}^{\infty}
D\nabla V_\infty(\Q(\tau;x_0))\vec{w}(\tau)d\tau \,ds\qquad
(t\geq 0)\]
\beqno
\lefteqn{\hspace*{-10mm}I_3^\gamma:=\sum_{N=1}^g 
\sum_{\stackrel{\gamma^{(1)}+\ldots+\gamma^{(N)}=\gamma}{|\gamma^{(i)}|>0}}
\int_{t}^{\infty} \int_{s}^{\infty} D^N \nabla 
\l[V_\infty(\Q(\tau;x_0))- V_\infty(\q(\tau,x_0))\ri]}&&\\
&&\hspace*{50mm} \l(\pa^{\gamma^{(1)}}_{x_0}\q(\tau,x_0),\ldots, 
\pa^{\gamma^{(N)}}_{x_0}\q(\tau,x_0)\ri)d\tau \,ds
\eeqno
and
\beqno
\lefteqn{I_4^\gamma:=\sum_{N=2}^g 
\sum_{\stackrel{\gamma^{(1)}+\ldots+\gamma^{(N)}=\gamma}{|\gamma^{(i)}|>0}}}&&\\
&&\int_{t}^{\infty} \int_{s}^{\infty} 
\l[D^N \nabla V_\infty(\Q(\tau;x_0))
\l(\pa^{\gamma^{(1)}}_{x_0}\Q(\tau;x_0),\ldots, 
\pa^{\gamma^{(N)}}_{x_0}\Q(\tau;x_0)\ri)\ri.\\
&&\hspace*{12mm}\l.-D^N \nabla 
V_\infty(\Q(\tau;x_0))\l(\pa^{\gamma^{(1)}}_{x_0}\q(\tau,x_0),\ldots, 
\pa^{\gamma^{(N)}}_{x_0}\q(\tau,x_0)\ri)\ri]d\tau \,ds .
\eeqno
On the space of bounded $\vec{w}$, $\cP\equiv\cP_{x_0}$ is bounded,
since for $\sup_t|\vec{w}(t)|=1$
\beqno
|\cP(\vec{w})(t)|&\leq& |\Zi|q_0^{-3} \int_{0}^{\infty} \int_{s}^{\infty}
\langle \tau\rangle_{\lambda}^{-3}d\tau \,ds\\
&=&\frac{2|\Zi|}{q_0E}\int_0^\infty \l(1-\frac{u}{\sqrt{1+u^2}}\ri)du=
\frac{2|\Zi|}{q_0E}.
\eeqno
Assuming (\ref{ass:out}), 
a suitable choice of the $E$-dependence of $\Rvir$ consistent 
with assumption (\ref{R:vir}) thus gives 
\[\|\cP_{x_0}\|\leq \eh \qmbox{and} 
\|\cP_{x_0}\|=\cO(1/(q_0H(x_0))), \]
so that we can invert $\idty-\cP$ for all energies $E>0$.

$I_3^\gamma$ and $I_4^\gamma$ only contain derivatives 
$\pa^{\gamma'}_{x_0}\vec{r}(\tau;x_0)$ with $|\gamma'|<g$.

Majorizing $V_\infty(\Q)-V_\infty(\q)$ by $2|DV_\infty(\q)\vec{r}|$,
the terms in $I_3^\gamma$ are estimated by 
\beqno
\lefteqn{\CC 
q_0^{-N-3}\int_{0}^{\infty} \int_{s}^{\infty} 
\langle \tau\rangle_{\lambda}^{-3}|\vec{r}(\tau;x_0)|d\tau \,ds
\prod_{i=1}^N\|\pa^{\gamma^{(i)}}_{x_0}\q(\cdot,x_0)\|_{\lambda}}&&\\
&=&\cO\l(q_0^{-|\beta|-1}E^{-\eh |\alpha|-1}\ri),
\eeqno
since $|\vec{r}(t,x_0)| < \Rmin$ (see (\ref{C:E})), and using (\ref{ind:ass}).

Finally, inserting $\Q=\q-\vec{r}$ into the first $N$--linear form
of $I_4^\gamma$ and expanding, we see that each term contains at least one 
factor of the form $\pa^{\gamma^{(i)}}_{x_0}\vec{r}(\cdot,x_0)$. So
\[I_4^\gamma=\cO\l(q_0^{-|\beta|-1}E^{-\eh |\alpha|-1}\ri),\] 
too.

Together with the above estimates for $I_1^\gamma$ and  $I_3^\gamma$ 
this finally shows that
\[\pa^{\gamma}_{x_0}\vec{r}(\tau;x_0)= 
(\idty-\cP)^{-1}(I_1^\gamma+I_3^\gamma+I_4^\gamma)=
\cO( q_0^{-|\beta|-\vep}E^{-\eh |\alpha|-1}), \]
which is equivalent to the second assertion in (\ref{DMoe:E}). 
By a parametrized version of the Banach Fixed Point Theorem
(see, {\em e.g.} \cite{DG}, Prop.\ A.2.2),
$x_0\mapsto \pa^{\gamma}_{x_0}\vec{r}(\cdot;x_0)$ is continuous.

The first assertion in (\ref{DMoe:E}) follows similarly, since
\[ \pa^\gamma_{x_0}(\P_0-\p_0) = -
\int_{t}^{\infty} \ \pa^\gamma_{x_0}\nabla
\l(V_\infty(\Q(\tau;x_0)) - V(\q(\tau,x_0)) \ri)    d\tau.\hspace*{10mm} \Box\]
%
%
\Section{The Flow between Near-Collisions} \label{sect:between}
%
We remind the reader of the relation (\ref{V:small}) for 
the virial radius $\Rvir(E)$. We assume that $\Eth>\Vmax$ so that we may
assume $\Rvir(E)\equiv\Rvir$ for all $E>\Eth$.
 
By the virial identity (\ref{virial}) a trajectory of energy $E$
leaving the interaction zone
$ \IZ = \{\q\in \Mu\mid \B{\q} \leq \Rvir(\Eth)\} $
cannot reenter it. 

We show first that a particle cannot stay inside the interaction zone 
for a long time without having close encounters with the nuclei.
Then we will control long trajectories within $\IZ$ by Poincar\'{e}
section techniques.

To quantify this, let
\[ \IZ(r):=\IZ\setminus \bigcup_{l=1}^{n} {\rm int}( B_l(r)),\]
$B_l(r)=\{\q\in \Mu\mid |\q-\s_l| \leq r \}$
being the ball of radius $r$ around the $l$th nucleus.

Instead of considering the restriction $\Phi^t\rstr_{\SuE}$ of the
flow generated by $H$, it is technically convenient to consider the 
flow  $\PutE$ generated by 
\[\Hh_E(\vec{v},\vec{x}):=\eh |\vec{v}|^2 + \frac{V(\vec{x})}{2E}
\qmbox{on}(\Hh_E)^{-1}(\eh).\]

We then have for 
$(\p(t),\q(t))\equiv\Phi^t\pqs$ and 
$(\vec{v}(s),\vec{x}(s))\equiv \Phi^s_E(\vec{v}_0,\vec{x}_0)$ with
\[(\vec{v}_0,\vec{x}_0):=(\p_0/\sqrt{2E},\q_0)\]
\beq
(\p(t),\q(t)) = \l( \sqrt{2E}\vec{v}(\sqrt{2E}t),\vec{x}(\sqrt{2E}t)\ri).
\Leq{umrechnen}
In the lemma below we use the standard Euclidean metric on 
$\bR^3_{\vec{v}}\times \bR^3_{\vec{x}}$.
\begin{lemma} \label{lem:C1}
For $E>\Eth$ the rescaled flow $\PutE$ in $\IZ(c_q)$ is $C^1$--near to 
the linear flow generated by the Hamiltonian function 
$\Hh_\infty\equiv \eh|\vec{v}|^2$ in the sense that
\beqn
\lefteqn{\sup_{s\in[0,T]} \l(
\l|\Phi^s_E(\vec{v}_0,\vec{x}_0)-(\vec{v}_0,\vec{x}_0+s\vec{v}_0)\ri| +
\l\| D\Phi^s_E(\vec{v}_0,\vec{x}_0) - 
\bsm \idty & 0 \\ s \idty & \idty\esm \ri\| \ri)}\NN\\
&=&\cO(1/E)
\label{C1:near}
\eeqn
if
\beq
x([0,T])\subset\IZ(c_q)
\Leq{in:iz}
Moreover, condition (\ref{in:iz}) is never satisfied if $T\geq3\Rvir$.
\end{lemma}
{\bf Proof.}
We {\em assume} that $T\leq 3\Rvir$, prove estimate (\ref{C1:near}) and then
the necessity of that assumption. The constants
\[L_0:=\sup_{\q\in\IZ(c_q)} |V(\q)|\quad,\quad 
  L_1:=\sup_{\q\in\IZ(c_q)} |\nabla V(\q)|\]
and
\[L_2:=\sup_{\q\in\IZ(c_q)} |D^2V(\q)|\]
are finite, since $\IZ(c_q)$ is compact. We bound the differences
\[Z_1(s):=\l(\vec{v}(s)-\vec{v}_0,\vec{x}(s)-(\vec{x}_0+s\vec{v}_0)\ri)\]
and
\[Z_2(s):= D\Phi^s_E(\vec{v}_0,\vec{x}_0) - 
\bsm \idty & 0 \\ s \idty & \idty\esm \]
in (\ref{C1:near}) using their integral equations
\[\vec{x}(s)-(\vec{x}_0+s\vec{v}_0) = 
- (2E)^{-1} \int_0^s \int_0^u \nabla V(\vec{x}(\tau))\, d\tau\, du\]
and
\[Z_2(s) = - (2E)^{-1} 
     \int_0^s \bem{cc} 0 & D^2V(\vec{x}(\tau)) \\ 0 & 0 \eem\, d\tau\]
and thus obtain
\[\sup_{s\in[0,T]} |Z_1(s)|\leq L_1\frac{T+T^2/2}{2E}\qmbox{,}
\sup_{s\in[0,T]} |Z_2(s)|\leq L_2\frac{T}{2E}.\]
This is indeed of order $\cO(1/E)$. 

Moreover, for $E>\Eth\geq \max(4L_0,L_1\Rvir)$
we have 
\[|\vec{v}_0|=\sqrt{1-V(\vec{x}_0)/E}\geq \frac{\sqrt{3}}{2}\] 
and
$|\vec{x}(s)-(\vec{x}_0+s\vec{v}_0)|\leq \Rvir/4$ for $s\leq T\leq 3\Rvir$
so that for $T:=3\Rvir$
\[|\vec{x}(T)-\vec{x}_0| \geq |\vec{v}_0|T-\Rvir/4 \geq 
\l( \frac{3\sqrt{3}}{2}-\frac{1}{4} \ri) \Rvir > 2\Rvir.\]
Thus for $T=3\Rvir$ we obtain a contradiction with our assumption 
(\ref{in:iz}), since the diameter of $\IZ$ equals $2\Rvir$ so that
$\vec{x}_0$ and $\vec{x}(T)$ cannot be both $\in\IZ$. \hfill $\Box$
%
\Section{The Single Scattering Process} \label{sect:si:sc}
%
We now consider motion inside the ball $B_l(c_q)$ near the $l$th singularity, 
using the KS estimates of Prop.\ \ref{propo:key}.
For initial conditions in the phase space region
\[\cD_l:=\pi(\Dol) = \{x\in P\mid |\q(x)-\s_l|\leq c_q, H(x)>\Eth\}\]
over the ball the exit times
\[T^\pm_l:\cD_l\ar\bR\qmbox{,}
  T^\pm_l(x_0):=\pm\inf\{t\geq 0\mid |\q(\pm t,x_0)-\s_l|=c_q\}\]
are well-defined and smooth. We compare them with the exit times
$T^\pm_{L,l}$ for the purely Keplerian motion in $B_l(c_q)$
generated by (\ref{pure:kepler}), and set
\[\Psul^\pm:\cD_l\ar\pa\cD_l, \quad x\mapsto \Phi(T^\pm_l(x),x)\]
and $\PsuLl^\pm(x):=\Phi_L(T^\pm_{L,l}(x),x)$.

Similar to (\ref{ToLl:0}) we define the {\em pericentric hypersurface}
\beqn
\Hul&:=&\pi(\Hol)\label{def:Hul}\\
&=&\l\{ x\in \cD_l\mid \q(x)=\s_l \qmbox{or} \LA\p(x),\q(x)-\s_l\RA=0\ri\}\NN
\eeqn
on which the {\em pericentric time} $T_l^0:\cD_l\ar\bR$
vanishes. If $Z_l>0$ we smoothly extend $T_l^0$ by demanding that 
\[\Psul^0(x):=\Phi(T_l^0(x),x)\in\Hul,\qquad (x\in\cD_l).\]
In analogy to (\ref{exit}) we introduce the diffeomorphism
\beq
\Psul:\pa\cD_l\ar\pa\cD_l ,\quad x\equiv (\p,\q)\mapsto 
\l\{\begin{array}{ll}
\Psul^+(x) & ,\LA \p,\q-\s_l\RA \leq 0\\
\Psul^-(x) & ,\LA \p,\q-\s_l\RA > 0\end{array} \ri.
\Leq{u:exit}
that permutes incoming and outgoing data
and its analog $\PsuLl$ for the Kepler flow.

We first consider the Kepler flow (with potential $\q\mapsto -Z_l/|\q-\s_l|$)
with incoming, resp.\ outgoing coordinates
\[\bem{c} \p^{\,-}\\ \q^{\,-}\eem \in \pa\cD_l\cap\SuE\qmbox{and}
\bem{c} \p^{\,+}\\ \q^{\,+}\eem := \PsuLl \bem{c} \p^{\,-}\\ \q^{\,-}\eem. \]
It is more convenient to work with the coordinates
\beq
\vec{v}^\pm:=\frac{\p^{\,\pm}}{\sqrt{2E}}\qmbox{and} 
\w^\pm:=\frac{\q^{\,\pm}-\s_l}{c_q}.
\Leq{vw:coord}

\begin{lemma} \label{lem:single:Kepler}
The Kepler transformation $\PsuLl$ is given by
\beq
\l\{\begin{array}{rcl}
\vec{v}^+(\vec{v}^-,\w^-)&=& \frac{(1-u^2(1+\beta))\vec{v}^-
-u\beta\w^-}{1-u^2}\\
\w^+(\vec{v}^-,\w^-)&=& \frac{2u(1-u^2(1+\beta/2))\vec{v}^- +
(1-u^2(1+\beta))\w^-}{1-u^2},\end{array}\ri.
\Leq{vw:p}
where
\beq
 \beta:= \frac{Z_l}{c_q E} \qmbox{and}
 u\equiv u(\vec{v}^-,\w^-) := -\frac{\LA \vec{v}^-,\w^-\RA}{1+\eh\beta} = 
-u(\vec{v}^+,\w^+).
\Leq{def:u}
\end{lemma}
{\bf Proof.} 
For simplicity of notation we assume $\s_l=\vec{0}$.
Since the KS transformation (\ref{def:KS}) transforms the scalar product in the
formula (\ref{PsoLl})
\[u = \frac{2\LA\tilde{P},Q\RA}{|\tilde{P}|^2+|Q|^2}\]
into 
\[ \LA\tilde{P},Q\RA = \frac{\tr(PQ^*)}{4\sqrt{2E}}=
-\frac{\tr(pq^*)}{2\sqrt{2E}}=-\frac{\LA\p,\q\RA}{\sqrt{2E}},\]
$|\tilde{P}|^2=c_q|\vec{v}^-|^2=c_q(1+\beta)$ and $|Q|^2=c_q$, 
we see that the definition (\ref{def:u}) of $u$ is consistent with 
the definition in (\ref{PsoLl}).

Eqs.\ (\ref{vw:p}) are obtained by inserting (\ref{PsoLl})
into the KS formula (\ref{def:KS}) for $(\p^{\,+},\q^{\,+})$, noticing that
\[-\tilde{P}^*I_3\tilde{P} = \frac{pqp}{2E}= 
\frac{|\p|^2\q-2\LA\p,\q\RA\p}{2E}\]
follows from the three-term-identity 
\[
YXY=2\LA X^*,Y \RA Y - \LA Y,Y\RA X^*\qquad (X,Y\in\bH) 
\]
of quaternions, applied to
$\p,\q\in\ImH$. \hfill$\Box$\\[2mm]
As we assumed that the NC-condition holds (no three singularities on a line), 
bounded orbits must be scattered by an 
angle $\sphericalangle(\vec{v}^+,\vec{v}^-)\leq \vartheta$
which must be at least of the order of $\amin>0$ (see (\ref{def:amin})).
This can only happen if the mismatch between the initial 
velocity $\vec{v}^-$ and position $\vec{w}^-$ is only of the order $\cO(1/E)$.
We first show this for Kepler scattering: 
\begin{lemma} \label{lem:large:dev}
For all energies $E>\Eth\geq 4 \Zmax/c_q$, angles 
$\theta\in [\pi\cdot\Eth/E,\pi]$, and initial conditions 
$(\p^{\,-},\q^{\,-})\in\pa\cD_l\cap\SuE$ and $u\geq 0$ with 
\beq
\l|\frac{\vec{v}^-}{|\vec{v}^-|} + \w^-\ri|\geq 
\frac{\pi\cdot \Eth}{\theta\cdot E},
\Leq{large:dev}
Kepler scattering is in the forward direction,
namely the total scattering angle is bounded by
$\sphericalangle(\vec{v}^+,\vec{v}^-)\leq \theta$.
\end{lemma}
{\bf Proof.}
Since the parameter $u\geq 0$ from (\ref{def:u}) equals
\[u(\vec{v}^-,\w^-)=-
\frac{\LA \vec{v}^-/|\vec{v}^-|,\w^-\RA\sqrt{1+\beta}}{1+\eh\beta} 
\leq -\LA \frac{\vec{v}^-}{|\vec{v}^-|},\w^-\RA\]
so that
\beq
1-u^2=(1+u)(1-u) \geq 1-u \geq 
\eh \l|\frac{\vec{v}^-}{|\vec{v}^-|} + \w^-\ri|^2,
\Leq{double:use} 
we get from formula (\ref{vw:p}) for $\vec{v}^+$ that
\beqno
|\vec{v}^+-\vec{v}^-| &=& \frac{u|\beta|}{1-u^2}|u\vec{v}^- + \w^-|\\
&\leq&\frac{u|\beta|}{1-u^2}\l(\l|\frac{\vec{v}^-}{|\vec{v}^-|} + \w^-\ri|
+|\vec{v}^-|\l((1-u)+\l|1-1/|\vec{v}^-|\ri|\ri) \ri),
\eeqno
with $\beta$ from (\ref{def:u}).

For $E>\Eth$ we have $|\beta|<|Z_l|/(c_q \Eth)\leq \ev$, so that
$|\vec{v}^-|=\sqrt{1+\beta}\in[\eh,2]$. 
Thus inserting the last two estimates of (\ref{double:use}) 
and (\ref{large:dev}) gives
\beqno
|\vec{v}^+-\vec{v}^-| &\leq& 
2|\beta|\l(\frac{\theta\cdot E}{\pi\cdot \Eth} +
1+\frac{\l|1-1/|\vec{v}^-|\ri|}{1-u^2} \ri)\\
&=& \frac{2|Z_l|}{c_q} \l(\frac{\theta}{\pi\cdot \Eth} + E^{-1} + 
E^{-1}\frac{\l|1-\frac{1}{\sqrt{1+\beta}} \ri|}{1-u^2}\ri)\\
&\leq&\frac{2|Z_l|}{c_q}
\l( \frac{\theta}{\pi\cdot \Eth}+\frac{\theta}{\pi\cdot \Eth}+ 
\frac{2|Z_l|}{c_q}\frac{\theta^2}{\pi^2\cdot \Eth^2} \ri)\\
&\leq&\frac{\theta}{2\pi}+\frac{\theta}{2\pi}+\frac{\theta^2}{(2\pi)^2}
\leq \frac{5}{4\pi}\theta
\eeqno
This proves the assertion, since 
\[\sphericalangle(\vec{v}^+,\vec{v}^-)\leq 
\frac{\pi}{2}\frac{|\vec{v}^+-\vec{v}^-|}{|\vec{v}^-|}\leq 
\frac{5}{8}(1+\beta)^{-1/2}\theta\leq \theta.\hspace*{4cm}\Box\]
Later on we will study the linearization of the flow $\Pt$, using Poincar\'{e} 
section techniques.
Therefore we now calculate the linearization of the Kepler transformation
(\ref{vw:p}), for 
tangent vectors $(\delta\p^{\,-},\delta\q^{\,-})$ in the four-dimensional 
subspace
\[T_{(\p^{\,-},\q^{\,-})}(\pa\cD_l\cap\SuE),\]
that is, for variations $(\delta\vec{v}^-,\delta\w^-)$ of 
(\ref{vw:coord}) meeting
\beq
\LA\vec{v}^-,\delta\vec{v}^-\RA=0\qmbox{,} \LA\w^-,\delta\w^-\RA=0.
\Leq{perp:pert}
\begin{lemma} \label{lem:lin:kepler}
The linearization of the Kepler Transformation (\ref{vw:p}) is given by
$$\delta\vec{v}^+ = \frac{(1-u^2)[((1-u^2)-u^2\beta)\delta\vec{v}^- - 
u\beta\delta\w^-]-\beta[2u\vec{v}^-+(1+u^2)\w^-]du}{(1-u^2)^2}$$
and
\beqno
\delta\w^+ &=& \frac{(2u(1-u^2)-u^3\beta)\delta\vec{v}^-
+((1-u^2)-u^2\beta)\delta\w^-}{1-u^2}\\
&&+\frac{\l[(2(1-u^2)^2+\beta u^2(u^2-3) )\vec{v}^- -2u\beta \w^-\ri]du}
{(1-u^2)^2}
\eeqno
with 
\[du := -\frac{
\LA\w^-,\delta\vec{v}^-\RA+\LA\vec{v}^-,\delta\w^-\RA}{1+\eh\beta}.\]
\end{lemma}
{\bf Proof.} By differentiation of (\ref{vw:p}).\hfill$\Box$\\[2mm]
As the next lemma shows, there are two regimes of Kepler scattering:
\begin{itemize}
\item {\bf hard scattering:}
If the 
scattering angle is larger than $c_1/\sqrt{E}$, then linearized scattering is 
basically 
a reflection combined with scaling. In particular, for an energy-independent
scattering angle the Liapunov exponent is of the approximate size $E$.
\item {\bf soft scattering:}
If, however, the scattering angle $\Delta\psi$
is smaller than $c_2/\sqrt{E}$, then
linearized scattering is a perturbation of free motion. 
\end{itemize}
\begin{lemma} \label{lem:lin}
${\bf 1)}$ For $c_1>0$ and $0\leq \delta\leq \eh$
we consider initial conditions $(\p^{\,-},\q^{\,-})\in\pa\cD_l\cap\SuE$
leading to a scattering angle
$\Delta\psi:=\sphericalangle(\vec{v}^+,\vec{v}^-)> c_1 E^{-\delta}$.
Then the linearized Kepler Transformation of Lemma \ref{lem:lin:kepler} 
is estimated by
\beq
\bem{c}\delta\vec{v}^+\\ \delta\w^+\eem = 
\frac{4c_q \sin^2(\eh\Delta\psi)}{-Z_l} E\cdot\bem{cc}R&R\\ R&R\eem
\bem{c}\delta\vec{v}^-\\ \delta\w^-\eem+
\l|\bem{c}\delta\vec{v}^-\\ \delta\w^-\eem\ri|\cdot\cO(E^0)
\Leq{est:dvwp}
where $R\equiv R_{\vec{v}^+-\vec{v}^-}\in O(3,\bR)$ is the reflection
by the plane perpendicular to the vector $\vec{v}^+-\vec{v}^-$.\\
${\bf 2)}$ 
If instead $\Delta\psi < c_2 E^{-\delta}$ with $\eh\leq \delta\leq1$, then
\beqn
\delta\vec{v}^+ &=& \delta\vec{v}^- + \hspace*{34mm}
\cO(c_2^2E^{1-2\delta})\cdot (|\delta\vec{v}^-|+|\delta\w^-|),\NN\\
\delta\vec{w}^+ &=& 2u \delta\vec{v}^- + \delta\vec{w}^- + 2\vec{v}^- du +
\cO(c_2^2E^{1-2\delta})\cdot (|\delta\vec{v}^-|+|\delta\w^-|).
\label{dev:from:free}
\eeqn
\end{lemma}
\begin{remark}
For $\delta\neq \eh$ the constants $c_1,c_2$ become irrelevant 
(for a threshold energy $\Eth\geq1)$.

Observe that the leading term in (\ref{est:dvwp}) scales like 
$c_1^2 E^{1-2\delta}\gg1$, whereas the error term of 
(\ref{dev:from:free}) scales like $c_2^2 E^{1-2\delta}\ll1$.
\end{remark}
{\bf Proof.}
{\bf 1)}
First we prove the identity
\beq
1-u^2=\frac{e^2}{(2/\beta+1)^2}
\Leq{emu}
with the eccentricity 
\beq
e=\sqrt{1+2E|\vec{L}_l|^2/Z_l^2}= 1/\sin(\eh\Delta\psi) 
\Leq{ecce}
of the Kepler hyperbola. By reinserting the $(\p,\q)$--coordinates with the
help of (\ref{vw:coord}), def.\ (\ref{def:u}) of $u$
acquires the form $u= -2\frac{\LA \p^-,\q^-\RA}{\sqrt{2E}c_q(2+\beta)}$
(assuming without loss of generality that $\s_l=\vec{0}$).
We now insert the relationship 
\[\LA \p^-,\q^-\RA^2= |\p^-|^2\,|\q^-|^2-\vec{L}_l^2=
2(E+Z_l/c_q)c_q^2-\vec{L}_l^2\]
into the expression for $1-u^2$ to deduce (\ref{emu}).

If $\Delta\psi> c E^{-\delta}$,
\beq
1-u^2=\cO(c^2E^{2(\delta-1)})
\Leq{emu2}
so that
\beq
\delta\vec{v}^+ = -\frac{2\beta(\vec{v}^-+\w^-)du}{(1-u^2)^2}
-\frac{\beta(\delta\vec{v}^-+\delta\w^-)}{1-u^2}+\cO(E^0)\cdot
(|\delta\vec{v}^-|+|\delta\w^-|)
\Leq{dvp}
and
\beq
\delta\w^+ = -\frac{2\beta(\vec{v}^-+\w^-)du}{(1-u^2)^2}
-\frac{\beta(\delta\vec{v}^-+\delta\w^-)}{1-u^2}+\cO(E^0)\cdot
(|\delta\vec{v}^-|+|\delta\w^-|).
\Leq{dwp}
Using (\ref{perp:pert}) we get
\[du= -\frac{\LA\w^-+\vec{v}^-,\delta\vec{v}^-+ \delta\w^-\RA}{1+\eh\beta}.\]
The vectors $\w^-$ and $\vec{v}^-$ in the last formula are nearly
anti-parallel.
From the expression (\ref{vw:p}) for $\vec{v}^+$ and (\ref{emu2}) we see that
$\vec{v}^+ = \vec{v}^- -\frac{\beta(\vec{v}^- +\w^-)}{1-u^2}+\cO(1/E)$ or
\[\vec{v}^- +\w^- = (1-u^2)\l(\frac{\vec{v}^- -\vec{v}^+}{\beta}+\cO(1)\ri)\]
so that
\beqno
\delta\vec{v}^+ &=& \frac{2(\vec{v}^- -\vec{v}^+ + \cO(1/E))\LA\vec{v}^- -
\vec{v}^+ + \cO(1/E),\delta\vec{v}^-+ \delta\w^-\RA}{\beta(1+\beta/2)}\\
& &-\frac{\beta(\delta\vec{v}^-+\delta\w^-)}{1-u^2}
+\cO(E^0)\cdot
(|\delta\vec{v}^-|+|\delta\w^-|)\\
&=& \frac{2(\vec{v}^- -\vec{v}^+)\LA\vec{v}^- -
\vec{v}^+ ,\delta\vec{v}^-+ \delta\w^-\RA}{\beta}
-\frac{\beta(\delta\vec{v}^-+\delta\w^-)}{1-u^2}\\ & &+\cO(E^0)\cdot
(|\delta\vec{v}^-|+|\delta\w^-|)\\
\eeqno
Extracting the square of the norm 
\beqno
|\vec{v}^- - \vec{v}^+|&=&
\l|\frac{\vec{v}^-}{|\vec{v}^-|} - \frac{\vec{v}^+}{|\vec{v}^+|}\ri|+\cO(1/E)\\
&=&\sqrt{2(1-\cos(\Delta\psi))}+\cO(1/E) = 2\sin(\eh\Delta\psi)+\cO(1/E),
\eeqno
we obtain the first estimate in (\ref{est:dvwp}), since
by (\ref{emu}) and (\ref{ecce})
\[\frac{\beta}{1-u^2} = \frac{4}{\beta e^2}+\cO(E^0)= 
\frac{4\sin^2(\eh\Delta\psi)}{\beta}+\cO(E^0),\]
and the reflection equals
\[R_{\vec{v}^+-\vec{v}^-}= 
\idty-2\frac{\l|\vec{v}^- - \vec{v}^+\RA\LA \vec{v}^- - \vec{v}^+\ri|}
{|\vec{v}^- - \vec{v}^+|^2}.\]
Estimate (\ref{est:dvwp}) for $\delta\w^+$ 
follows since the r.h.s.\ of (\ref{dwp}) and (\ref{dvp})
have the same form.\\[2mm]
{\bf 2)}
We write
\[\delta\vec{v}^+ = \delta\vec{v}^- + R_1 + R_2\qmbox{,} 
  \delta\vec{w}^+ = 2u \delta\vec{v}^- + \delta\vec{w}^- + 2\vec{v}^- du +
  R_1+R_3\] 
with rest terms
\[R_1:= -\frac{\beta}{1-u^2}\l[ \delta\vec{v}^- + \delta\vec{w}^- +
2\frac{(\vec{v}^- + \vec{w}^-)du}{1-u^2} \ri], \]
\[R_2:= \beta\l[ \delta\vec{v}^- + \frac{\delta\vec{w}^-}{1+u}+
\frac{\vec{v}^-du}{(1+u)^2} 
+ \frac{(\vec{v}^- + \vec{w}^-)du}{1-u^2} \ri]\]
and
\[R_3:= \beta \l[
\frac{1+u+u^2}{1+u}\delta\vec{v}^- - \delta\vec{w}^- +
\vec{v}^-du + \frac{\vec{w}^-du}{(1+u)^2} 
- \frac{(\vec{v}^- + \vec{w}^-)du}{1-u^2} \ri] \]
Now we assume that $\Delta\psi < c_2 E^{-\delta}$ with $\eh\leq \delta\leq1$.

So by the first part of the proof 
$R_1=\cO(c_2^2E^{1-2\delta})\cdot (|\delta\vec{v}^-|+|\delta\w^-|)$. 
By the same reasoning the (identical, up to sign) last term in $R_2$ and $R_3$
is of order $\cO(1/E)\cdot (|\delta\vec{v}^-|+|\delta\w^-|)$. 
The other terms in $R_2$ and $R_3$ are
of the same order, since they contain the multiplier $\beta$.
\hfill$\Box$\\[2mm]
We now use the results of Section \ref{sect:KS:Appl} to compare the true motion 
with the Kepler motion of the last lemmata. 
\begin{proposition} \label{prop:near:u}
For pericentric initial conditions $x_0\in\Hul$ with energy
$E:=H(x_0)>\Eth$
\beq
T^\pm_l(x_0)=T^\pm_{L,l}(x_0)+\cO(E^{-3/2}),
\Leq{T:diff}
and for $(\p^{\,\pm},\q^{\,\pm}):=\Psul^\pm(x_0)$ and 
$(\p_L^\pm,\q_L^\pm):=\PsuLl^\pm(x_0)$
\beq
\l( \p^{\,\pm}/\sqrt{2E},\q^{\,\pm} \ri) = 
\l( \p_L^\pm/\sqrt{2E},\q_L^\pm \ri)+\cO(1/E).
\Leq{C0:unten}
For $|\q_0-\s_l|=\cO(1/E)$
\beq
T^\pm_{L,l}(x_0) = \pm \l( \frac{c_q}{\sqrt{2E}} -
\frac{Z_l}{(2E)^{3/2}}\ln(E c_q/|Z_l|)\ri) +\cO(E^{-3/2}).
\Leq{T:l:est}
For arbitrary $x_0\in \Hul$ the r.h.s.\ of (\ref{T:l:est}) is an upper bound
for $|T^\pm_l(x_0)|$.
The diffeomorphism 
\[\Xi_l: \pa \cD_l\ar\pa\cD_{L,l},\quad x\equiv (\p,\q)\mapsto 
\l\{ \begin{array}{ll}\PsuLl^-\circ\Psul^0(x) & ,\LA \p,\q\RA \leq 0\\
\PsuLl^+\circ\Psul^0(x) & ,\LA\p,\q \RA > 0\end{array}  \ri.\] 
onto its image which conjugates the maps
($\Psul= (\Xi_l)^{-1}\circ\PsuLl\circ\Xi_l$)
is $C^0$-near to the identity in the sense that 
in the Euclidean norm for the $(\vec{v},\w)$-coordinates (\ref{vw:coord})
\beq
|\Xi_l(x)-x| = \cO(1/H(x)),\qquad (x\in\pa\cD_l),
\Leq{est:xiu}
and the solution $\PsuLl$ of the linear problem is 
$C^1$-near to $\Psul$ in the sense
\beqn
|\Psul(x)-\PsuLl\circ\Xi_l(x)| &=& \cO(H^{-1}(x))
\label{Pu:one}\\
 \l\|D\Psul(x) - D\PsuLl\circ\Xi_l(x) \ri\| &=& \cO(H^0(x)).
 \label{Pu:two}
\eeqn
\end{proposition}
\begin{remark}
In particular the exit times (\ref{T:l:est}) are
independent of $x_0$, up to $\cO(E^{-3/2})$, if $|\q_0-\s_l|=\cO(1/E)$.
\end{remark} 
{\bf Proof.} 
Estimate (\ref{C0:unten}) follows directly from (\ref{est:xio}).

Estimate (\ref{C0:near}) says that for KS initial conditions $X_0$ with 
$\pi(X_0)=x_0$ the relative error
\[\frac{|\PoLl(s,X_0)-\Pol(s,X_0)|}{|Q(s,X_0)|}=\cO(E^{-1}).\]
This translates into the estimates 
\beq
\frac{|\q(t,x_0)-\q_L(t,x_0)|}{|\q(t,x_0)|}=\cO(E^{-1})
\Leq{rel:q}
and 
\beq
\frac{|\p(t,x_0)-\p_L(t,x_0)|}{|\p(t,x_0)|}=\cO(E^{-1})\qquad (t\neq 0).
\Leq{rel:p}
We know that inside $\cD_l$ each trajectory $t\mapsto \Phi^t(x_0)$ 
passes the hypersurface $\Hul$ only at $x_0$, so that we can parametrize the 
curves $(0,T^+_l(x_0)]\ni t\mapsto \Phi^t(x_0)\equiv(\p(t),\q(t))$
by their radius $r(t):=|\q(t)-\s_l|$ (similarly for negative times).

The explicit formula for the exit time
\[T^+_l(x_0)=\frac{1}{\sqrt{2}}\int_{|\q_0|}^{c_q} \l(
E+\frac{Z_l}{r}-\frac{\hat{L}_l^2(\p(r),\q(r))}{2r^2} +W_l(\q(r))\ri)^{-\eh}dr\]
(with the smooth potential $W_l$ defined in (\ref{Wl}), and relative angular
momentum $\hat{L}_l$ in (\ref{def:hat:Ll}))
is obtained by integrating the inverse of the radial velocity 
$\LA\p,\q\RA/|\q|$. Similarly
\[T^+_{L,l}(x_0)=\frac{1}{\sqrt{2}}\int_{|\q_0|}^{c_q} \l(
E+\frac{Z_l}{r}-\frac{L_l^2(x_0)}{2r^2} \ri)^{-\eh}dr.\]
inserting (\ref{rel:q}) and (\ref{rel:p}) then shows the assertion
(\ref{T:diff}).

(\ref{T:l:est}) follows by evaluating the explicit formula
(\ref{eq:explicit:int}) for the exit time of the Kepler flow:
As $E\rmin^2+Z_l\rmin-\eh L_l^2=0$, the square root terms 
in (\ref{eq:explicit:int}) vanish for the minimal radius.
Further, by our assumption on $\q_0$ we have 
$\rmin\leq C|Z_l|/E$ for some $C>0$, so that 
\[L_l^2\leq |\q-\s_l|^2\p^{\,2}=2|\q-\s_l|(E|\q-\s_l|+Z_l)
\leq\frac{2Z_l^2 C(1+C)}{E}=\cO(1/E).\]

Est.\ (\ref{est:xiu}) follows from (\ref{est:xio}) and similarly
(\ref{Pu:one}), (\ref{Pu:two}) follow from (\ref{Po:one}) and (\ref{Po:two}),
since on $\pa\DoLl$ the KS transformation from coordinates $(\tilde{P},Q)$
to the $(\vec{v},\w)$-coordinates (\ref{vw:coord}) takes the simple
form
\[\bem{c}\vec{v}\\ \w \eem = \l.\bem{c}Q^*I_3\tilde{P}\\ Q^*I_3Q\eem\ri/ c_q.\]
\hfill $\Box$
%
\Section{Long Paths Within the Interaction Zone} \label{sect:long:paths}
%
Strong changes of directions only occur if the pericentric distance
of the orbit from a nucleus is of order $\cO(1/E)$.

To quantify this, we set for $0<\theta\leq\pi$
\beq
\Hul(\theta):= \l\{x\in\Hul\l| \
|\q(x)-\s_l|<\frac{|Z_l|}{H(x)\cdot\sin(\theta/4)}\ri.\ri\}
\Leq{Hul:theta}
(with the pericentric hypersurface $\Hul$ defined in (\ref{def:Hul})).
\begin{lemma} \label{lem:defl}
For $C_5>0$ large the following statements are true.
Consider pericentric initial conditions $x_0\in \Hul$
of energy $E:=H(x_0)>\Eth$ and $\theta\in[C_5/E,\pi]$.
\begin{itemize}
\item
If  $x_0 \not\in\Hul(\theta)$, 
then the directions $\hat{p}^\pm:=\hat{p}(\Psul^\pm(x_0))$ 
of the orbit through $x_0$ at the moment of entering, resp.\ 
exit from the ball $B_l(c_q)$ differ at most by the angle
\beq
 \Delta\psi = \sphericalangle \l( \hat{p}^-,\hat{p}^+ \ri) < \theta.
\Leq{winke:winke}
\item
Conversely, if $x_0\in \Hul(\theta)$, then 
\beq
\Delta\psi \geq \theta/4,
\Leq{0:0}
and positions and momenta are nearly
anti-parallel before scattering and parallel after scattering:
\beq
\sphericalangle\l( \q(\Psul^\pm(x_0))-\s_l , \pm\hat{p}^\pm \ri) =  
\cO((\theta E)^{-1}).
\Leq{equal:dir}
\end{itemize}
\end{lemma}
{\bf Proof.} 
Because of estimate (\ref{Pu:one}) of Prop.\ \ref{prop:near:u} 
we need only choose $C_5$ large and then show 
the following statements for the {\em Kepler} flow (which 
correspond to (\ref{winke:winke}),  (\ref{0:0}) resp.\ (\ref{equal:dir})):
\beq
\Delta\psi_L := \sphericalangle \l( \hat{p}_L^-,\hat{p}_L^+ \ri) \leq \eh\theta,
\Leq{small:K:angle}
\beq
\sin(\eh\Delta\psi_L) \geq 1.05\cdot\sin(\theta/8)
\Leq{0:1}
and
\beq
\sphericalangle\l( \q_L^\pm-\s_l , \pm\hat{p}_L^\pm \ri) =  
\cO((\theta E)^{-1}).
\Leq{gaehn}
with $\hat{p}_L^\pm:=\hat{p}(\PsuLl^\pm(x_0))$ and
$\q_L^\pm:=\q(\PsuLl^\pm(x_0))$.\\
$\bullet$
We use the explicit formulae of Lemma \ref{lem:single:Kepler} to compute
\beqno
\cos(\Delta\psi_L) &=& \frac{\LA \vec{v}^-,\vec{v}^+ \RA} {|\vec{v}^-|\, |\vec{v}^+|}= 
\frac{\LA \vec{v}^-,\vec{v}^+\RA}{1+\beta} = 1-\eh\frac{u^2\beta^2}{(1-u^2)(1+\beta)}\\
&=& 1-\eh\frac{u^2(2+\beta)^2}{e^2(1+\beta)},
\eeqno
using the substitution (\ref{emu}). Now the eccentricity $e$ of the Kepler
hyperbola equals 
\beq
e=|1+2E\rmin/Z_l|,
\Leq{nochneformel:fuer:e} 
as one sees by comparing the definition (\ref{rmin}) of 
the pericentric distance $\rmin$ with (\ref{ecce})
(formula (\ref{rmin}) is also valid if $Z_l<0$; 
then the energy is always positive). 
Thus 
\beqno
\sin(\eh \Delta\psi_L)&=&\sqrt{\eh(1-\cos(\Delta\psi_L))} = 
\frac{u(1+\eh\beta)}{\sqrt{1+\beta}}\frac{1}{|2E\rmin/Z_l + 1|}\\
&\leq& |2E\rmin/Z_l + 1|^{-1}\leq 
\frac{\sin(\ev\theta)}{2+\sin(\ev\theta)\sign(Z_l)} \leq  \sin(\ev\theta)
\eeqno
since $u(1+\eh\beta)/\sqrt{1+\beta}$ equals the cosine of the angle between
$\vec{v}^+$ and $\w^+$ and is thus smaller than one, see (\ref{def:u}).
Then (\ref{small:K:angle}) follows, since 
the sine function is monotone increasing on $[0,\pi/2]$.\\
$\bullet$ By (\ref{ecce}) and (\ref{nochneformel:fuer:e})
\beqno
\sin(\eh\Delta\psi_L)&=& 1/e= |1+2E\rmin/Z_l|^{-1}\geq
\frac{\sin(\theta/4)}{|1+\sign(Z_l)\sin(\theta/4)|}\\
&\geq& 
\frac{\sin(\theta/4)}{1+1/\sqrt{2}}
\geq\frac{\sin(\pi/4)}{2\sin(\pi/8)}\frac{2\sin(\theta/8)}{1+1/\sqrt{2}}
\eeqno
showing (\ref{0:1}).\\
$\bullet$ 
With (\ref{def:u}), (\ref{emu}) and (\ref{nochneformel:fuer:e})
\beqno
\lefteqn{\hspace*{-5mm}\sin\l(\sphericalangle\l( \q(\Psul^\pm(x_0))-\s_l ,
 \pm\hat{p}^\pm \ri)\ri)=
 \sin\l(\sphericalangle\l( \w^\pm ,\pm\vec{v}^\pm\ri)\ri)
=\sqrt{1-\frac{\LA\w^\pm ,\pm\vec{v}^\pm\RA^2} {\LA\w^\pm,\w^\pm\RA}}}\\
&=& \sqrt{1-u^2\frac{(1+\eh\beta)^2}{1+\beta}}= \sqrt{|1-u^2|}+\cO(1/E)\\
&=& \frac{e}{|2/\beta+1|}+\cO(1/E)= \l|\frac{1+2E\rmin/Z_l}{2/\beta+1}\ri|+\cO(1/E)\\ 
&=& \l|\frac{2\rmin/Z_l+1/E}{2c_q/Z_l+1/E}\ri|+\cO(1/E)=
\frac{\rmin}{c_q}+\cO(1/E) =\cO(1/(\theta E)), 
\eeqno
showing estimate (\ref{gaehn}).\hfill $\Box$\\[2mm]
Long trajectories in the interaction zone 
$\IZ$ must have close encounters with singularities of distance $\cO(1/E)$.
To show this, 
we now assume that the radius $c_q$ of the balls $B_k(c_q)$ is so small that
there is no straight line meeting more than two balls. To be concrete, we
assume
\beq
c_q\leq \ev\sin(\amin)\cdot\dmin.
\Leq{cq:in}
\begin{proposition} \label{propo:drei:R} %
For $\Eth$ large we consider a trajectory segment
\beq
[0,T]\ni t \mapsto \Pt(x_0)\equiv(\p(t),\q(t)) \in \cD
\Leq{tr:seg} 
with
\[\cD:=\{x\in P\mid H(x)>\Eth, \q(x)\in \IZ\},\]
which does not intersect any pericentric hypersurface $\Hul(\amin/2)$, 
$(l=1,\ldots,n)$.
Then, in configuration space, it does not reenter a ball $B_{k}(c_q)$ 
after leaving it,
and it does not intersect three or more such balls. 

Furthermore, the length of the time interval is bounded by 
$T<13\Rvir/\sqrt{2E}$.
\end{proposition}
{\bf Proof.} 
$\bullet$
There is no subsegment $[t_1,t_2]\ni t \mapsto \q(t)$ of the trajectory, 
lying in $\IZ(c_q)$,
which leaves and then reenters a given ball $B_k(c_q)$. 
Namely we may otherwise assume
\[\q([t_1,t_2])\subset \IZ(c_q)\qmbox{,}
|\q(t_1)-\s_k|= |\q(t_2)-\s_k|=c_q\qmbox{and} t_2>t_1.\] 
Then there would be a time
$t_0\in(t_1,t_2)$ with maximal distance $|\q(t_0)-\s_k|\geq c_q$ from $\s_k$:
\beq
\frac{d}{dt}|\q(t_0)-\s_k|^2=0\qmbox{and}\frac{d^2}{dt^2}|\q(t_0)-\s_k|^2\leq 0.
\Leq{second:der:u}
But 
\beqno
\lefteqn{\eh\frac{d^2}{dt^2} |\q(t_0)-\s_k|^2 = 
 |\p(t_0)|^2 -\LA\nabla V(\q(t_0)),\q(t_0)-\s_k\RA }\\
&\geq&  2(E-L_0)-\sup_{\q\in\IZ(c_q)} \LA\nabla V(\q),\q-\s_k\RA
\geq 2(E- L_0-L_1\Rvir)\\
\eeqno
with $L_0=\sup_{\q\in\IZ(c_q)} |V(\q)|$ and 
$L_1=\sup_{\q\in\IZ(c_q)} |\nabla V(\q)|$, 
so that (\ref{second:der:u}) does not hold for $E$ large.\\[2mm] 
$\bullet$
We now prove that there is no orbit segment (\ref{tr:seg}) whose
configuration space projection
intersects three balls $B_{k_0}(c_q)$, $B_{k_1}(c_q)$ and
$B_{k_2}(c_q)$ in succession. By going to a subsegment, we otherwise assume 
outgoing initial data 
\[\q(0)\in\pa B_{k_0}(c_q)\qmbox{,}\LA \p(0),\q(0)-\s_{k_0}\RA\geq 0,\]
ingoing final data 
\[\q(T)\in\pa B_{k_2}(c_q)\qmbox{,}\LA \p(T),\q(T)-\s_{k_2}\RA\leq 0,\]
and intermediate times 
\[0<t^-\leq t^+<T \qmbox{with} 
\q(t^\pm)\in\pa B_{k_1}(c_q) \qmbox{and}
\pm\LA \p(t^\pm),\q(t^\pm)-\s_{k_1}\RA\geq 0.\]

By the first part of the proof we know that $k_0\neq k_1\neq k_2$.
Then the angle $\gamma$ between the direction $\hat{s}^{k_0,k_1}$ 
of the axis through the centres and
$\hat{p}(t^-)$ is bounded by 
\beq
\gamma < \eh\amin.
\Leq{gamma:e}
Namely all straight lines intersecting $B_{k_0}(c_q)$ and $B_{k_1}(c_q)$ 
have a direction whose angle $\gamma'$ with $\hat{s}^{k_0,k_1}$ is bounded 
by
\[\sin(\gamma')\leq \frac{c_q}{\eh\dmin}\leq \eh\sin(\amin)<\sin(\eh\amin) ,\]
using (\ref{cq:in}), so that (\ref{gamma:e}) follows for $\Eth$ 
large from Lemma \ref{lem:C1}.

We now consider the unit vector 
$\hat{u}=\lambda_0 \hat{s}^{k_0,k_1}+ \lambda_2 \hat{s}^{k_2,k_1}$, 
$\lambda_i\geq 0$, which is perpendicular to $\s^{k_0,k_2}$.

Then by (\ref{gamma:e}) the angle between $\hat{p}(t^-)$ and $\hat{u}$ is bounded above
by 
\beqn
\sphericalangle(\hat{p}(t^-),\hat{u})&\leq&
\sphericalangle(\hat{s}^{k_0,k_1},\hat{u})+\gamma = 
\pi-\eh\pi-\alpha(k_1,k_0,k_2)+\gamma\NN\\
&\leq&\eh\pi- \amin + \gamma<\eh(\pi-\amin)
\label{angle:u}
\eeqn
(with $\cos(\alpha(i,j,k))=\LA\hat{s}^{j,i},\hat{s}^{j,k}\RA$).

On the other hand for some time $t\in[t^+,T]$ we have
\beq
\LA \hat{p}(t),\hat{u}\RA \leq
\LA \frac{\q(T)-\q(t^+)}{|\q(T)-\q(t^+)|},\hat{u}\RA <0,
\Leq{angle:pu}
because otherwise there is no trajectory between $\q(t^+)$ and
$\q(T)$. 

But by Lemma \ref{lem:C1} in the interval $t\in [t^+,T]$ 
the change of direction 
\[\sphericalangle(\hat{p}(t^+),\hat{p}(t)) = \cO(1/E),\]
so that with (\ref{angle:u}) and (\ref{angle:pu}) 
we must have a large scattering angle
\beq
\sphericalangle\l(\hat{p}(t^-),\hat{p}(t^+)\ri) \geq \eh\amin.
\Leq{contra}
within $B_{k_1}(c_q)$.

Our assumption that the orbit segment (\ref{tr:seg})
does not intersect the pericentric hypersurface 
${\cal H}_{k_1}(\amin/2)$ (see (\ref{Hul:theta}) implies by
Lemma \ref{lem:defl} that the deflection angle is bounded above by 
\[\sphericalangle\l(\hat{p}(t^-),\hat{p}(t^+)\ri)<\eh\amin\]
contradicting (\ref{contra}).

Thus a trajectory segment can be represented as the union of at most
two segments inside balls $B_{k_i}(c_q)$ and three segments
in $\IZ(c_q)$.\\[2mm] 
$\bullet$
By (\ref{T:diff}) and (\ref{T:l:est}) for large $\Eth$
the time interval spent by the trajectory inside the ball $B_{k_i}(c_q)$ has
length
$T^+_{k_i}-T^-_{k_i}\leq \frac{4c_q}{\sqrt{2E}} < \frac{2\Rvir}{\sqrt{2E}}$.
By Lemma \ref{lem:C1} with (\ref{umrechnen}) the 
time interval of the segments in $\IZ(c_q)$ are of lengths $\leq
3\Rvir/\sqrt{2E}$. 
So the total time $T$ is bounded above by
\[T<2\frac{2\Rvir}{\sqrt{2E}}+3\frac{3\Rvir}{\sqrt{2E}}=
13\frac{\Rvir}{\sqrt{2E}}.\hspace*{8cm} \Box\]
\bigskip

We now control the paths between successive close encounters with the centres
$\s_k$. 
For radii $c_q>0$ meeting (\ref{cq:in})
the balls $B_k(c_q)$, $k\inInd$ around the centres $\s_k$ do not intersect.
 
By the NC condition \ref{defi:noncol} the distance ${\rm dist}(\s_i, A^{k,l})$, $i\neq k,l$
between the centre and the {\em axis}
\[\Ax^{k,l}:=\{t\s_l+(1-t)\s_k\mid t\in [0,1]\},\qquad (k\neq l)\]
connecting $\s_k$ and $\s_l$ is bounded below by $\sin(\amin)\cdot\dmin$. 

\begin{figure}
\begin{center}
\begin{picture}(0,0)%
\epsfig{file=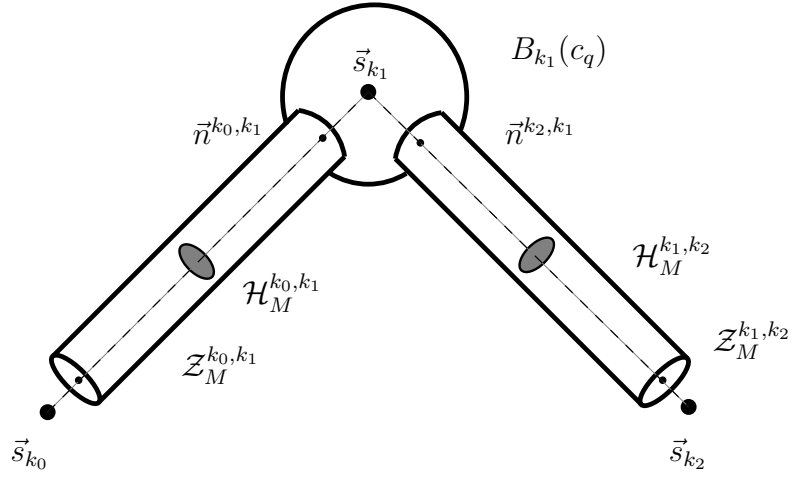}%
\end{picture}%
\setlength{\unitlength}{0.00083300in}%
\begingroup\makeatletter\ifx\SetFigFont\undefined
\def\x#1#2#3#4#5#6#7\relax{\def\x{#1#2#3#4#5#6}}%
\expandafter\x\fmtname xxxxxx\relax \def\y{splain}%
\ifx\x\y   
\gdef\SetFigFont#1#2#3{%
  \ifnum #1<17\tiny\else \ifnum #1<20\small\else
  \ifnum #1<24\normalsize\else \ifnum #1<29\large\else
  \ifnum #1<34\Large\else \ifnum #1<41\LARGE\else
     \huge\fi\fi\fi\fi\fi\fi
  \csname #3\endcsname}%
\else
\gdef\SetFigFont#1#2#3{\begingroup
  \count@#1\relax \ifnum 25<\count@\count@25\fi
  \def\x{\endgroup\@setsize\SetFigFont{#2pt}}%
  \expandafter\x
    \csname \romannumeral\the\count@ pt\expandafter\endcsname
    \csname @\romannumeral\the\count@ pt\endcsname
  \csname #3\endcsname}%
\fi
\fi\endgroup
\begin{picture}(4380,2916)(64,-2117)
\put( 64,-2083){\makebox(0,0)[lb]{\smash{\SetFigFont{12}{14.4}{rm}$\s_{k_0}$}}}
\put(4444,-1393){\makebox(0,0)[lb]{\smash{\SetFigFont{12}{14.4}{rm}$\Cyl_M^{k_1,k_
2}$}}}
\put(1509,-1095){\makebox(0,0)[lb]{\smash{\SetFigFont{12}{14.4}{rm}$\Poi_M^{k_0,k_
1}$}}}
\put(3182,427){\makebox(0,0)[lb]{\smash{\SetFigFont{12}{14.4}{rm}$B_{k_1}(c_q)$}}}
\put(1128,-1550){\makebox(0,0)[lb]{\smash{\SetFigFont{12}{14.4}{rm}$\Cyl_M^{k_0,k_
1}$}}}
\put(2194,352){\makebox(0,0)[lb]{\smash{\SetFigFont{12}{14.4}{rm}$\s_{k_1}$}}}
\put(3148,-102){\makebox(0,0)[lb]{\smash{\SetFigFont{12}{14.4}{rm}$\vec{n}^{k_2,k_
1}$}}}
\put(1200,-102){\makebox(0,0)[lb]{\smash{\SetFigFont{12}{14.4}{rm}$\vec{n}^{k_0,k_
1}$}}}
\put(3942,-865){\makebox(0,0)[lb]{\smash{\SetFigFont{12}{14.4}{rm}$\Poi_M^{k_1,k_2
}$}}}
\put(4170,-2083){\makebox(0,0)[lb]{\smash{\SetFigFont{12}{14.4}{rm}$\s_{k_2}$}}}
\end{picture}

\end{center}
\caption{Configuration space projections $\Poi_M^{k,l}$ of Poincar\'e sections} \label{fig:one}
\end{figure}

Thus the configuration space cylinders (see Figure \ref{fig:one})
$$\Cyl_M^{k,l} := \l\{\q\in\Mu \mid {\rm dist}(\q,\Ax^{k,l})\leq c_y,\, 
|\q-\s_k|\geq c_q \leq |\q-\s_l| \ri\}$$
have empty intersection with non-adjacent balls
($\Cyl_M^{k,l}\cap B_i(c_q)=\emptyset$ for $i\neq k,l$) if $c_y\leq c_q$.
\begin{itemize}
\item
The cylinders themselves are not mutually disjoint;\\ 
in particular
$\Cyl_M^{l,k}=\Cyl_M^{k,l}$. 
\item
However, under the NC condition
\[\Cyl_M^{k_1,l_1}\cap \Cyl_M^{k_2,l_2}=\emptyset\]
if the axes $\Ax^{k_1,l_1}\neq \Ax^{k_2,l_2}$ are parallel.
\item
For the choice 
\[c_y:=\eh\sin(\amin/2) c_q\]
of the radius we have in addition 
\[\Cyl_M^{k,l_1}\cap \Cyl_M^{k,l_2}=\emptyset\qquad (l_1\neq l_2).\]
\item
Finally there is a minimal nonzero angle between nonparallel axes.
More precisely, for this choice of $c_y$ a cylinder $\Cyl_M^{k_1,l_1}$
can only intersect a different ($\{k_1,l_1\}\neq\{k_2,l_2\}$) 
cylinder $\Cyl_M^{k_2,l_2}$ if the angle 
$\arccos(|\LA\hat{s}^{k_1,l_1},\hat{s}^{k_2,l_2}\RA|)$ between their axes
is larger than $\amin$. The proof of this fact is a nice exercise in conic
sections.
\end{itemize}
%
%
Thus for 
\[c_p:= \eh \min(\eh\amin, c_y/\dmax)\]
and $E>\Eth$ the phase space regions
$$\Cyl^{k,l} := \l\{ x\equiv(\p,\q)\in \cD  \l|\ \q\in\Cyl_M^{k,l},\
\l| \p/\sqrt{2 H(x)}-\hat{s}^{k,l}\ri| \leq c_p \ri.\ri\}$$
(with $\cD=\{x\in P\mid H(x)>\Eth, \q(x)\in \IZ\}$) 
do not intersect ($\Cyl^{k_1,l_1}\cap \Cyl^{k_2,l_2} =\emptyset$ for
$(k_1,l_1)\neq(k_2,l_2)$).

For $k\neq l$ we erect {\em Poincar\'{e} hypersurfaces} $\Poi^{k,l}$ 
near the {\em midpoint}
\[\Mid^{k,l}:= \eh(\s_k+\s_l)\in \Ax^{k,l}\] 
by setting for $e_0>0$
\[\Poi_M^{k,l}(E) := \l\{\q\in \Mu\mid \LA\q-\Mid^{k,l},\hat{s}^{k,l}\RA=0,
|\q-\Mid^{k,l}| < \ea c_p\cdot d^{k,l} e_0/E \ri\} ,  \]
so that $\Poi_M^{k,l}(E)\subset\Cyl_M^{k,l}$ for $\Eth$ large and $E>\Eth$, and 
\beqn
\hspace{-0.5cm}\Poi^{k,l}_E &:=&\l\{(\p,\q)\in\SuE \l|\  \q\in\Poi_M^{k,l}(E),\ri.\ri.\NN\\
& &\hspace{2.2cm}\l. \LA\p,\hat{s}^{k,l}\RA>0,\,
\l| \frac{\p}{\sqrt{2 E}}\times\hat{s}^{k,l}\ri| < \eh  c_p e_0/E \ri\},
\label{def:PoiE}
\eeqn
so that $\Poi^{k,l}_E\subset\Cyl^{k,l}$, too;
see Figure \ref{fig:one}.
The $\Poi^{k,l}_E$ are four-dimensional submanifolds without boundary of $\SuE$.
Sometimes we work with their closures 
\beq
\ov{\Poi}^{k,l}_E, 
\Leq{clos:Poi}
which are
diffeomorphic to products of two closed disks (and as such formally speaking 
no manifolds).

High energy orbits between near-collisions are nearly straight and thus 
move near some axis:
\begin{lemma} \label{lem:HK}
Let $\Eth$ and the constant $e_0>0$ in (\ref{def:PoiE}) be large and 
consider for $E>\Eth$ trajectory segments of the form
\[ [-T,T]\ni t\mapsto x(t) \equiv (\p(t),\q(t)) := \Phi^t(x_0)\]
with initial values $x_0\in \SuE$, starting and ending in the pericentric
hypersurfaces (\ref{Hul:theta}):
\[x(\pm T)\in\cH_{k_\pm}(\amin/2)\]
but not intersecting a pericentric hypersurfaces in between:
\beq
x((-T,T))\cap\Hul(\amin/2)=\emptyset\qquad(l=1,\ldots,n).
\Leq{ass}
Then the trajectory segment is contained in the phase space region
\beq
x([-T,T]) \subset \Cyl^{k_-,k_+}\cup \cD_{k_-} \cup \cD_{k_+},
\Leq{Hundeknochen}
and there is a unique $t_0\in(-T,T)$ with
\beq
\Phi^{t_0}(x_0)\in\Poi^{k_-,k_+}_E.
\Leq{in}
Conversely for $k_+\neq k_-$ the exit time 
$$T^{k_-,k_+,\pm}_E : \Poi^{k_-,k_+}_E \ar \bR^\pm\ ,\
T^{k_-,k_+,\pm}_E(x_0):=
\pm\inf\l\{t\geq 0\mid \Phi^{\pm t}(x_0)\in \pa \Cyl^{k_-,k_+}\ri\}
$$
from $\Cyl^{k_-,k_+}$ is of order 
\beq
T^{k_-,k_+,\pm}_E =\pm \frac{\eh d^{k_-,k_+} - c_q}{\sqrt{2E}} +\cO(E^{-3/2}),
\Leq{noch:ne:zeit}
and the exit points 
$(\p^{\,\pm},\q^{\,\pm}):=\Phi(T^{k_-,k_+,\pm}_E(x_0),x_0)$
are estimated by
\beq
(\p^{\,\pm}/\sqrt{2E},\q^{\,\pm})=
(\hat{s}^{k_-,k_+}, \vec{n\,}^{k_\mp,k_\pm}) + \cO(1/E),
\Leq{ray:family}
with the intersection point 
$\vec{n\,}^{k,l} := \s_l + c_q\hat{s}^{l,k}$
between the axis $\Ax^{k,l}$ and the sphere $\pa B_{l}(c_q)$, 
see Figure \ref{fig:one}.
\end{lemma} 
{\bf Proof. }\\
By straight line geometry and Lemma \ref{lem:C1} trajectories
$t\mapsto \q(t,x_0)$ with $x_0\in \Poi^{k_-,k_+}_E$ intersect 
$\pa \Cyl_M^{k_-,k_+}$ near $\vec{n\,}^{k_\mp,k_\pm}$, and
\[\p^{\,\pm}/\sqrt{2E}-\hat{s}^{k_-,k_+}=
(\p^{\,\pm}/\sqrt{2E}-\p_0)+(\p_0-\hat{s}^{k_-,k_+})=\cO(1/E),\] 
again using 
Lemma \ref{lem:C1} and the definition (\ref {def:PoiE}) of $\Poi^{k_-,k_+}_E$.
This shows (\ref{ray:family}).
Since $|\Mid^{k_-,k_+}-\vec{n\,}^{k_\mp,k_\pm}| = \eh d^{k_-,k_+} - c_q$,
(\ref{noch:ne:zeit}) follows from Lemma \ref{lem:C1}.

Concerning the first statement of the lemma, 
we know from (\ref{equal:dir}) of Lemma \ref{lem:defl} that for 
$t^\pm:= T^\pm_{k_\mp}(x(\mp T))$
\[(\p^{\,\pm},\q^{\,\pm}) := \Phi\l( t^\mp , x(\pm T) \ri) \in
(\q\rstr_{\SuE})^{-1}(\pa B_{k_\pm}(c_q))\]
have the property 
\beq
\q^{\,\pm}-\s_{k_\pm} = \pm \frac{c_q}{\sqrt{2E}}\p^{\,\pm} +\cO(1/E).
\Leq{one}
Furthermore, by the assumption (\ref{ass}) and the first statement of Prop.\ 
\ref{propo:drei:R} the trajectory lies in $\IZ(c_q)$ during the time interval 
$[-T+t^+,T-t^-]$. Thus 
(\ref{C1:near}) and (\ref{umrechnen}) imply that
\beq
\frac{\p^{\,+}}{\sqrt{2E}}=\frac{\p^{\,-}}{\sqrt{2E}}+\cO(1/E) 
\Leq{two}
and
\beq
\frac{\q^{\,+}-\q^{\,-}}{|\q^{\,+}-\q^{\,-}|}=
\frac{\p^{\,-}}{\sqrt{2E}}+\cO(1/E). 
\Leq{three}
Using (\ref{one}), (\ref{three}) and then (\ref{two}),
\beqno
\s_{k_+}-\s_{k_-} &\equiv& (\s_{k_+} - \q^{\,+}) + 
(q^{\,-} -\s_{k_-})+(\q^{\,+}-\q^{\,-})\\ 
&=& \frac{c_q}{\sqrt{2E}}(\p^{\,+} + 
\p^{\,-}) + \frac{|\q^{\,+}-\q^{\,-}|}{\sqrt{2E}}\p^{\,-}
+\cO(1/E)\\
&=& \frac{2c_q+|\q^{\,+}-\q^{\,-}|}{\sqrt{2E}}\p^{\,\pm} +\cO(1/E)
\eeqno
so that 
\[ \l|\p^{\,\pm}/\sqrt{2E}-\hat{s}^{k_-,k_+}\ri| = \cO(1/E).\]
A second application of (\ref{three}) shows that the second term in
(\ref{ray:family}), too is $\cO(1/E)$.

By Lemma \ref{lem:C1} we obtain (\ref{Hundeknochen}) and (\ref{in}).
\hfill $\Box$
%
\Section{The Poincar\'{e} Map} \label{sect:P:map}
%
For $E>\Eth$ consider the {\em Poincar\'{e} surfaces}
\beq
\Poi_E:= \bigcup_{\stackrel{k,l=1}{ k\neq l}}^{n} \Poi^{k,l}_E\quad,\quad 
\Poi^\pm_E := \Poi_E\cup \pa\cD_E^\pm
\Leq{all:Poincare}  
with $\Poi^{k,l}_E$ defined in (\ref{def:PoiE}) and
\beq
\cD_E:= \cD\cap\SuE,\qquad 
\pa\cD_E^\pm:= \l\{ (\p,\q)\in \pa\cD_E\mid  \pm\LA\p,\q\RA \geq 0\ri\}. 
\Leq{cDE}
Then the {\em return time} to the Poincar\'{e} surface ${\Poi}_E^+$
\[\TuR:\cD_E\ar [0,\infty)\cup \{\infty\},\] 
is defined by $\TuR(x) := 0$ for $x\in\pa\cD_E^+$ and 
\beq
\TuR(x) := \inf \l\{ t>0 \l| \Phi^t_E(x)\in {\Poi}_E^+ \ri. \ri\},\qquad
(x\in\cD_E\setminus \pa\cD_E^+).
\Leq{return:time} 
\begin{lemma}\label{lem:return}
For $E>\Eth$ the Poincar\'{e} return time $\TuR$ is finite and
\[\TuR =\cO(1/\sqrt{E}).\]
\end{lemma}
{\bf Proof.}
Let $x_0\in\cD_E\setminus \pa\cD_E^+$.
By Prop.\ \ref{propo:drei:R} for a time 
$t_1\in [0,13\Rvir/\sqrt{2E}]$ the point
$x_1:=\Phi^{t_1}(x_0)$
\begin{itemize}
\item
either exits $\cD_E$, that is
$x_1\in\pa\cD_E^+$
\item
or meets a pericentric surface near the $k$th nucleus:
$x_1\in\cH_k(\amin/2)$.
\end{itemize}
In the first case we are done. In the second case
we iterate the above argument and find
$t_2\in [t_1,t_1+13\Rvir/\sqrt{2E}]$ for which the point
$x_2:=\Phi^{t_2}(x_0)$
either exits $\cD_E$ ($x_2\in\pa\cD_E^+$)
or meets a pericentric surface near the $l\neq k$th nucleus:
\[x_2\in\cH_l(\amin/2).\]
In the relevant second case we know from (\ref{in}) of Lemma \ref{lem:HK}
that 
\[ x_3:=\Phi^{t_3}(x_0)\in \Poi_E^{k,l}\]
for some $t_3\in(t_1,t_2)$.
\hfill $\Box$\\[2mm]
We shall analyze the 
{\em Poincar\'{e} map} 
\beq
\Po: \Poi^-_E\ar \Poi^+_E ,\qquad
\Po(x) := \Phi(\TuR(x),x).
\Leq{Poincare:map}
When we work with the closures
\[\ov{\Poi}_E:= \bigcup_{\stackrel{k,l=1}{ k\neq l}}^{n} \ov{\Poi}^{k,l}_E\quad,\quad 
\ov{\Poi}^\pm_E := \ov{\Poi}_E\cup \pa\cD_E^\pm\]  
(see (\ref{clos:Poi})), we write $\ov{T}_E$, $\ov{{\cal P}}_E$ etc.
\begin{lemma} \label{lem:bij}
$\Po$ is a bijection, and its restriction to
$\Poi^-_E\setminus \Po^{-1}(\pa \Poi_E^+ )$ is smooth.
\end{lemma}
{\bf Proof.} \\
{\bf $\Po$ is one-to-one:} 
If $\Po(x_1)=\Po(x_2)$ and $\TuR(x_1)\leq \TuR(x_2)$, then
$x_1=\Phi^{t_2-t_1}(x_2)$ and
\begin{itemize}
\item
either $x_1\in \pa\cD_E^-$. Then $x_2=x_1$ since by the virial inequality
\[\Phi^t(\Poi^-_E) \cap\pa\cD_E^- = \emptyset\qmbox{for}t>0,\]
\item
or $x_1\in \Poi_E$. Then by definition of $\TuR$ using
the infimum,  $x_2=x_1$, too.
\end{itemize}
{\bf $\Po$ is onto:}
{\em Time reversal}
\beq
\TR:P\ar P, \qquad(\p,\q)\in T^*\Muh\mapsto (-\p,\q),\quad x\in
P\setminus T^*\Muh\mapsto x 
\Leq{time:reversal}
is a smooth anti-symplectic transformation, with
\[\TR\circ\Phi^t\circ {\TR}=\Phi^{-t}\qquad(t\in\bR)\]
and
\[\TR(\Poi^\pm_E)=\Poi^\mp_E.\] 
Thus
\[\Po^{-1}(x)=\TR\circ\Po\circ\TR(x).\]
{\bf Smoothness} of $\Po\rstr_{U}$ for 
$U:= \Poi^-_E\setminus \Po^{-1}(\pa \Poi_E^+)$ 
follows from transversality of the codimension one
$\pa$--mani\-fold $U\subset\SuE$ to the flow $\Phi^t\rstr_{\SuE}$, since 
$\Po\rstr_{U}$ maps to inner points of $\Poi^+_E$. \hfill $\Box$
\subsection*{Symbolic Dynamics}
In order to use symbolic dynamics, we introduce 
{\em symbol sequences} 
\[\uk =(k_i)_{i\in I}\in \cS^I\]
over the alphabet 
\[\cS := \{1,\ldots,n\},\]
where 
\beq
I\equiv I_l^r:=\{i\in\bZ\mid l \leq i \leq  r \}
\Leq{interval} 
for $l,r\in\bZ\cup\{\pm\infty\}$ is a 
(finite, half-infinite or bi-infinite) {\em interval}.

$\uk$ is called {\em admissible} if 
$k_i\neq k_{i+1}$ for all $\{i,i+1\}\subset I$.

For $E>\Eth$ and $(k_0,k_1)$ admissible we set
\beq
V_E(k_{0},k_{1}):= W_E(k_{0},k_{1}):= \Poi^{k_0,k_1}_E,
\Leq{VWE}
and for $(k_{-m},\ldots,k_0)$ admissible
\beq
W_E(k_{-m},\ldots,k_0):= W_E(k_{-1},k_0)\cap\Po(W_E(k_{-m},\ldots,k_{-1})),
\Leq{WEM}
resp.\ for $(k_{0},\ldots,k_m)$ admissible
\beq
V_E(k_{0},\ldots,k_m):= V_E(k_{0},k_1)\cap\Po^{-1}(V_E(k_{1},\ldots,k_{m})),
\Leq{VEM}
$m\geq 2$.
Then by Lemma \ref{lem:bij} the iterated maps
\beqno
\Po(k_{0},\ldots,k_{m})&:& V_E(k_{0},\ldots,k_{m})\ar W_E(k_{0},\ldots,k_{m})\\
\Po(k_{0},\ldots,k_{m}) &:=& \Po^{m-1}\rstr_{V_E(k_{0},\ldots,k_{m})} 
\eeqno
are diffeomorphisms (we will show in Prop.\ \ref{prop:bVb} 
that the sets $V_E(\uk)$ are non-empty).
Again, $\ov{W}_E(\uk)$ and $\ov{V}_E(\uk)$ denote the closures of these sets.

We decompose the Poincar\'{e} map $\Po(k_{-1},k_0,k_1)$ in the form
\beqn
\hspace{-10mm}\Po(k_{-1},k_0,k_1)(x)&=&
\Po^-(k_{-1},k_0)\circ \Po'(k_{-1},k_0,k_1)\circ\Po^+(k_0,k_1)(x)
\label{Poinc:decomp}\\
\hspace{-10mm}&\equiv& \Po^-\circ \Po'\circ\Po^+(x)\qquad \qquad (x\in V_E(k_{-1},k_0,k_1))
\NN
\eeqn
with the diffeomorphisms
\beqno
\Po^-:V_E(k_{-1},k_0)&\ar&
 V'_E(k_{-1},k_0):=\Po^-(V_E(k_{-1},k_0))\subset\pa \Cyl^{k_{-1},k_0}\\
x         &\mapsto& \Phi \l( T^{k_{-1},k_0,+}_E(x),x \ri),
\eeqno
\beqno
\l(\Po^+\ri)^{-1}\!:W_E(k_0,k_1)\!&\ar&
  W'_E(k_0,k_1):= \l(\Po^+\ri)^{-1}\!(W_E(k_0,k_1))\subset\pa \Cyl^{k_0,k_1}\\
            x         &\mapsto& \Phi \l(T^{k_0,k_1,-}_E(x),x \ri),
\eeqno
and
\[ \Po':V'_E(k_{-1},k_0,k_1) \ar  W'_E(k_{-1},k_0,k_1)\quad,\quad
x  \mapsto \Psul^+(x)\]
with $\Psul^+$ defined in (\ref{u:exit})) for
\beqno
 V'_E(k_{-1},k_0,k_1)&:=&\Po^-(V_E(k_{-1},k_0,k_1))\subset  V'_E(k_{-1},k_0),\\
 W'_E(k_{-1},k_0,k_1)&:=&(\Po^+)^{-1}(W_E(k_{-1},k_0,k_1))\subset 
 W'_E(k_0,k_1).
\eeqno
$\Po^\pm$ are considered as perturbations of Poincar\'e maps for the free flow, whereas
Prop.~\ref{propo:key} allows us to view $ \Po'$ as a perturbation of
Poincar\'e maps for the Kepler flow.
\subsection*{Adapted Coordinates}
We now define adapted coordinates in the Poincar\'e sections $\Poi_E$.
This will simplify to introduce invariant cone fields later on. 

We complement the unit vectors 
\[\hat{s}^\pm\equiv \hat{s}^\pm(k_{-1},k_0,k_1)
\qmbox{with}\hat{s}^\pm:=\hat{s}^{k_{\pm 1},k_0}\] 
pointing towards $\s_{k_0}$ 
by a unit vector
\beq
\hat{t}^\pm\equiv \hat{t}^\pm(k_{-1},k_0,k_1)\qmbox{with}\hat{t}^+=\hat{t}^-
\mbox{ perpendicular to } {\rm span}(\hat{s}^-,\hat{s}^+)
\Leq{def:t:unit}
(by the NC condition \ref{defi:noncol} the span
is one-dimensional iff $k_1=k_{-1}$), and set 
\beq
\hat{u}^\pm \equiv \hat{u}^\pm(k_{-1},k_0,k_1):= 
\hat{s}^\pm \times \hat{t}^\pm .
\Leq{def:u:unit}
Then $\hat{t}^\pm$ and $\hat{u}^\pm$, considered 
as elements of $T_\q \Mu$ for $\q\in \ov{\Poi}_M^{k_{\pm1},k_{0}}$, 
form an orthonormal basis
of $T_\q \ov{\Poi}_M^{k_{\pm1},k_{0}}$.
We introduce adapted coordinates on $V_E(k_{-1},k_0)$ and $W_E(k_0,k_1)$ by
mapping
\beqno
V_E(k_{-1},k_0) &\ar& \bR^2\times\bR^2\qmbox{,}
(\p,\q) \mapsto (\y^-,\z^-),\\
W_E(k_0,k_1) &\ar& \bR^2\times\bR^2\qmbox{,}
(\p,\q) \mapsto (\y^+,\z^+),
\eeqno
with
\beq
\y^\pm:=
\l.\bem{c} \LA \p ,\hat{t}^\pm \RA\\ \LA \p,\hat{u}^\pm \RA \eem \ri/
\sqrt{2E}\quad,\quad
\z^\pm:=        
\l.\bem{c} \LA \q - \Mid^{k_{\pm1},k_0},\hat{t}^\pm \RA\\ 
           \LA \q - \Mid^{k_{\pm1},k_0},\hat{u}^\pm \RA \eem \ri/ l^\pm 
\Leq{def:yz}            
and $l^\pm:=|\Mid^{k_{0,\pm1}}-\s_{k_0}|=\eh d^{k_{\pm1},k_0}$.
\begin{remark} \label{rem:rotation}
The coordinates $(\y^-,\z^-)$ on $V_E(k_{-1},k_0)$ depend through the vector
$\hat{t}^-\equiv \hat{t}^-(k_{-1},k_0,k_1)$ on the symbol $k_1$.
However, the coordinate systems for the symbols
$k_1^I$ and $k_1^{II}$ are related by
\[\bem{c}\y^{-,I}\\\z^{-,I}\eem = \bem{c}O\,\y^{-,II}\\O\,\z^{-,II}\eem,\]
where $O\in SO(2,\bR)$ is the rotation in the plane spanned by 
$\hat{t}^{-,I}$ and $\hat{u}^{-,I}$ which 
maps $\hat{t}^{-,I}$ onto $\hat{t}^{-,II}$ (and correspondingly 
$\hat{u}^{-,I}$ onto $\hat{u}^{-,II}$).
\end{remark}
Similarly we introduce coordinates in
\beqno
V'_E(k_{-1},k_0) &\ar& \bR^2\times\bR^2, \qquad
(\p,\q) \mapsto (\y'^{-},\z'^-)\\
W'_E(k_0,k_1) &\ar& \bR^2\times\bR^2, \qquad
(\p,\q) \mapsto (\y'^+,\z'^+)
\eeqno             
with 
\beq
 \y'^\pm:=
\l.\bem{c} \LA \p ,\hat{t}^\pm \RA\\ \LA \p,\hat{u}^\pm \RA \eem \ri/ \sqrt{2E}\quad,\quad
\z'^\pm:=        
\l.\bem{c} \LA \q - \vec{n}^{k_{\pm1},k_0},\hat{t}^\pm \RA\\ 
           \LA \q - \vec{n}^{k_{\pm1},k_0},\hat{u}^\pm \RA \eem \ri/ l'^\pm ,
\Leq{def:yz:prime}             
\[l'^-:=|\vec{n}^{k_{-1},k_0}-\s_{k_{-1}}|= d^{k_{-1},k_0}-c_q\qmbox{and}
l'^+:=|\vec{n}^{k_{1},k_0}-\s_{k_0}|= c_q.\]

The next lemma shows that arbitrary pairs $\z^-,\z^+$ of points in the 
configuration space projections of the Poin\-car\'{e} surfaces
are connected by a trajectory.
\begin{lemma} \label{ueberall}
Choose a large enough constant $e_0$ in Def.\ (\ref{def:PoiE})
of the Poincar\'{e} surfaces. Then
for $\Eth\equiv\Eth(e_0)$ large, all 
$E>\Eth$ and
\[\z^\pm\in\Poi_M^{k_{\pm 1},k_0}(E)\] 
there exists a $x^-\in V_E(k_{-1},k_0,k_1)$ 
of the form 
\[x^-=(\y^-,\z^-)\qmbox{with} x^+:=\Po(x^-)=(\y^+,\z^+).  \]
\end{lemma}
{\bf Proof.} \\
By (\ref{def:PoiE}) in the $(\y^-,\z^-)$--coordinates 
$V_E(k_{-1},k_{0})= \Poi^{k_{-1},k_0}_E$ 
is of the form
\beq
V_E(k_{-1},k_{0})=B_y\times B_z
\Leq{VEByBz}
with the two-disks $B_y$ in momentum space and $B_z$ in configuration space
\beq
B_y:=\{\y\in\bR^2\mid |\y| <\eh c_p e_0/E\}\supset
B_z:=\{\z\in\bR^2\mid |\z| < \ev c_p e_0/E\}.
\Leq{ByBz}
We thus consider initial conditions of the form
\[x^-=(\y^-,\z^-)\in B_y\times \{\z^-\}\subset V_E(k_{-1},k_0)\]
and their images
\[x'^-=(\y'^-,\z'^-):=\Po^-(x^-)\in V'_E(k_{-1},k_0).\]
The idea of the proof is that, as $B_y$ has twice the size of $B_z$,
the family of trajectories with initial conditions
in $B_y\times \{\z^-\}$ has an opening angle $\cO(1/E)$ that is large enough
to hit a full $\cO(1/E)$--neighbourhood of $c_{k_0}$. 
Then by near collision, the opening angle is 
nearly amplified to $2\pi$, if $e_0$ is chosen large enough.

The threshold energy $\Eth$ depends on the constant $e_0$. However, 
by Lemma \ref{lem:C1} and a compactness argument one sees that the error term $\cO(1/E)$
in the $C^1$-estimates
\beq
\l|\bem{c}\y'^-\\ \z'^-\eem -
\bem{cc}\idty & 0\\ (1-\frac{l^-}{l'^-})\idty &\frac{l^-}{l'^-}\idty\eem
\bem{c}\y^-\\ \z^-\eem \ri|\leq C/E
\Leq{nlin:min}
\beq
\l|\bem{c}\y^+\\ \z^+\eem - 
\bem{cc}\idty & 0\\ (1-\frac{l'^+}{l^+})\idty &\frac{l'^+}{l^+}\idty\eem
\bem{c}\y'^+\\ \z'^+\eem \ri|\leq C/E
\Leq{nlin:plus}
and an identical estimate for the linearizations 
\beq
\l|\bem{c}\delta\y'^-\\ \delta\z'^-\eem -
\bem{cc}\idty & 0\\ (1-\frac{l^-}{l'^-})\idty &\frac{l^-}{l'^-}\idty\eem
\bem{c}\delta\y^-\\ \delta \z^-\eem \ri|\leq 
\frac{C}{E}\cdot 
\l|\bem{c}\delta\y^-\\ \delta\z^-\eem\ri|
\Leq{lin:min}
\beq
\l|\bem{c}\delta \y^+\\ \delta \z^+\eem - 
\bem{cc}\idty & 0\\ (1-\frac{l'^+}{l^+})\idty &\frac{l'^+}{l^+}\idty\eem
\bem{c}\delta \y'^+\\ \delta \z'^+\eem \ri|\leq 
\frac{C}{E}\cdot 
\l|\bem{c}\delta \y'^+\\ \delta \z'^+\eem \ri|
\Leq{lin:plus}
of these maps
do {\em not} depend on $e_0$, that is, $C>0$ may be chosen fixed and then
$e_0$ chosen arbitrarily large. 

The relations
\beqno
\p &=& \sqrt{2E^\pm}\cdot \l( \mp\hat{s}^\pm + 
y_1^\pm \hat{t}^\pm + y_2^\pm \hat{u}^\pm 
+ \cO(E^{-2})\ri)\\
\q &=& \Mid^\pm + l^\pm \l(z_1^\pm \hat{t}^\pm + z_2^\pm \hat{u}^\pm \ri)
\eeqno
with $E^\pm:= E-V(\Mid^\pm)$ show that the angular momentum relative to
the $k_0$-st nucleus before and after scattering  equals
$$\vec{L}_{k_0}\equiv (\q-\s_{k_0})\times \p =
\pm\sqrt{2E^\pm}l^\pm \cdot 
\l( (y_1^\pm+z_1^\pm) \hat{u}^\pm -(y_2^\pm+z_2^\pm) \hat{t}^\pm 
+ \cO(E^{-2})\ri) .$$
Thus $|\vec{L}_{k_0}|=\cO(e_0/\sqrt{E})$, and 
\[|\vec{L}_{k_0}|\geq \ev l^\pm c_p e_0/\sqrt{E} \qmbox{for}
(\y^-,\z^-)\qmbox{with} \y^-\in \pa \ov{B_y}.\]

Estimate (\ref{nlin:min}) shows that similar statements hold true for
 $\vec{L}_{k_0}\circ \Po^-$.

Finally (changing to $(\p,\q)$--coordinates for a moment)
by Prop.\ \ref{prop:near:u} we know that the image 
\[x'^+\equiv (\p\,'^+,\q\,'^+):=\psi_l^+(x'^-)\in W'_E(k_{-1},k_0,k_1)\]
of $x'^-\in V'_E(k_{-1},k_0,k_1)$ is near to the image 
\[x_L'^+\equiv (\p_L'^+,\q_L'^+):=
\psi_{L,l}^+(x_L'^-)\in W'_E(k_{-1},k_0,k_1)\] 
of a point $x_L'^-=\Xi(x'^-)$ 
in a $\cO(1/E)$--neighbourhood of $x'^-$ in the sense that
\beq |\p\,'^+-\p_L'^+|/\sqrt{2E}+|\q\,'^+-\q_L'^+|=\cO(1/E).
\Leq{nn}

For the Kepler exit map $\psi_{L,l}^+$ the formula 
\[e\equiv e(x_L'^-)=\sqrt{1+2E\vec{L}^2/Z^2}\]
for the eccentricity of the Kepler hyperbola parametrized by
\[r(\vv)=\frac{\vec{L}^2}{Z(1+e\cos(\vv))}\]
shows that the maximal total scattering angle $2\arccos(1/e)\in[0,\pi)$
can be increased to $\pi-\vep$ (uniformly in $E$) by increasing $e_0$, so that
the diffeomorphism 
\[B_y\times \{\z^-\}\ni x^-\mapsto \q_L'^+(\Xi\circ\Po^-(x^-)) \in \pa B_{k_0}(c_q) \]
onto its image covers the sphere $\pa B_{k_0}(c_q)$ around $\s_{k_0}$ 
except an $\vep$--neigh\-bour\-hood of the forward
direction $\s_{k_0}+c_q\cdot\hat{s}^-\in\pa B_{k_0}(c_q)$.

The same holds true for the diffeomorphism 
\[B_y\times \{\z^-\}\ni x^-\mapsto \q'^+(\Po^-(x^-)) \in \pa B_{k_0}(c_q) \]
onto its image. So by the NC condition \ref{defi:noncol} an
$\vep$--neighbourhood
$U\subset \pa B_{k_0}(c_q)$ of the point
$\vec{n\,}^{k_1,k_0} = \s_{k_0} - c_q\hat{s}^{+}$
is contained in that image for $e_0$ large. 
By Def.\ (\ref{def:yz:prime}) of the $\z'^+$--coordinates this corresponds
to a $\vep$--neighbourhood of $\z'^+=0$. Estimate (\ref{ray:family}) 
then implies that on this neighbourhood
$$ \y'^+( \Po'\circ\Po^-(x^-)) = \z'^+( \Po'\circ\Po^-(x^-)) + \cO(1/E)
\qquad (x^-\in B_y\times \{\z^-\}).$$
Thus we conclude from  (\ref{nlin:plus}) that the intersection
\beq
\Po(B_y\times \{\z^-\})\equiv
\Po^+\circ\Po'\circ\Po^-(B_y\times \{\z^-\})\cap B_y\times \{\z^+\}
\subset W(k_0,k_1)
\Leq{inters}
is non-empty.\hfill $\Box$
\begin{remark}
In fact the point $x^-\in V_E(k_{-1},k_0,k_1)$ is unique, as will
follow from Prop.\ \ref{propo:cf} below.
\end{remark}
\begin{lemma} \label{ret:time:est}
For all admissible $(k_{-1},k_0,k_1)\in \cS^3$ the Poincar\'{e} return 
time equals
$$T_E\rstr_{V_E(k_{-1},k_0,k_1)}=
\frac{d^{k_{-1},k_0}+d^{k_0,k_1}}{2\sqrt{2E}}
-\frac{Z_{k_0}}{(2E)^{3/2}}\ln\!\l(\!E\frac{d^{k_{-1},k_0}+d^{k_0,k_1}}{2|Z_{k_0}|}\!\ri)
+\cO(E^{-3/2}).$$
\end{lemma}
{\bf Proof.}
\[T_E(x)=T^{k_{-1},k_0,+}_E(x) + 
(T^+_{k_0}(y)-T^-_{k_0}(y)) - T^{k_{0},k_1,-}_E(\Po(x))\]
with $y:=\psi^0_{k_0}\circ\Phi(T^{k_{-1},k_0,+}(x),x)$.
The estimate thus follows from (\ref{noch:ne:zeit}), (\ref{T:diff}) and 
(\ref{T:l:est}).
\hfill $\Box$
%
\Section{Existence of an Invariant Cone Field} \label{sect:ICF}
%
\begin{definition}
A {\em cone at} $x\in\Poi_E$ 
is the image of the {\em standard cone} 
\[\{(\vec{u},\vec{v})\in\bR^2\times\bR^2\mid
|\vec{u}|\geq|\vec{v}|\}\]
w.r.t.\ an invertible linear map 
\[ \bR^2\times\bR^2\ar T_x\Poi_E.\]
A {\em cone field} $\cC$ in $U\subset\Poi_E$ associates to every
$x\in U$ a cone $\cC(x)$ at $x$.\\
$\cC$ is called {\em invariant} if for every $x\in U$ with $\Po(x)\in U$
\[T_x\Po(\cC(x))\subset \cC(\Po(x)),\]
and {\em strictly invariant} if
\[T_x\Po(\cC(x))\subset {\rm Int}(\cC(\Po(x)))\cup\{0\}.\]
\end{definition}
Alternatively one may think of cones as subsets of $\bR P^3$.

In the next proposition we use the Poincar\'{e} section
coordinates $(\y^\pm,\z^\pm)$ defined in (\ref{def:yz}).
\begin{proposition} \label{propo:cf}
For $C>0$ large and $E>\Eth$ the cone field $\cC$ in 
$$U := \bigcup_{\stackrel{(k_{-1},k_0,k_1)\in \cS^3}{\rm
admissible}}V_E(k_{-1},k_0,k_1)$$
\beq
\cC(x) := 
\l\{(\delta\y,\delta\z)\in T_{(\y,\z)} V_E(k_{-1},k_0,k_1)\mid
|\delta\y-\delta\z|\leq \frac{C}{E} |\delta\y+\delta\z|\ri\}
\Leq{cone}
for $x\equiv (\y,\z)\in U$, is strictly $\Po$--invariant, and the linearized
Poincar\'{e} map equals
\beq
T_x\Po = f(k_{-1},k_0,k_1)E\cdot\bem{cc}
\idty&\idty\\\idty&\idty\eem+\cO(E^0)\qquad (x\in V_E(k_{-1},k_0,k_1)),
\Leq{Tx}
with
\beq
f(k_{-1},k_0,k_1):=
\frac{2d^{k_{-1},k_0}\cos^2(\eh\alpha(k_{-1},k_0,k_1))}{-Z_{k_0}}.
\Leq{def:f}
\end{proposition}
\begin{remarks}
{\bf 1)}
The linearized Poincar\'{e} map (\ref{Tx}) is invertible, although
the matrix $\bsm\idty&\idty\\\idty&\idty\esm$ is not. \\
{\bf 2)}
Although the coordinates $(\y^+,\z^+)$ on 
$W_E(k_{-1},k_0,k_1)\subset\Poi^{k_0,k_1}_E$
do not coincide on their common domain 
$W_E(k_{-1},k_0,k_1)\cap V_E(k_{0},k_1,k_2)$ 
with the coordinates 
$(\y^-,\z^-)$ on $V_E(k_{0},k_1,k_2)\subset\Poi^{k_0,k_1}_E$, 
they are related by
\[\bem{c}\y^+\\\z^+\eem = \bem{c}O\,\y^-\\O\,\z^-\eem,\]
where $O\in O(2,\bR)$ is the reflection in the plane perpendicular to 
$\hat{s}^{k_{0},k_1}$, transforming the unit vector
$\hat{t}^+(k_{-1},k_0,k_1)$ into $\hat{t}^-(k_{0},k_1,k_2)$ and the unit vector
$\hat{u}^+(k_{-1},k_0,k_1)$ into $\hat{u}^-(k_{0},k_1,k_2)$, see
Def.\ (\ref{def:t:unit}) and (\ref{def:u:unit}).
Thus the above cone field $\cC$ is invariant under that transformation.
\end{remarks}
{\bf Proof.}\\
By (\ref{lin:min}) 
\beqn
\bem{c}\delta\y'^-\\ \delta\z'^- \eem &:=&
(T_{x^-}\Po^-)\bem{c}\delta\y^-\\ \delta\z^- \eem\label{chain1}\\ 
&=& 
\bem{cc}\idty & 0\\ (1-\frac{l^-}{l'^-})\idty &\frac{l^-}{l'^-}\idty\eem
\bem{c}\delta\y^-\\ \delta\z^- \eem +
\l|\bem{c}\delta\y^-\\ \delta \z^-\eem \ri|\cdot\cO(1/E).\NN
\eeqn
From (\ref{vw:coord}) and (\ref{def:yz:prime}) we deduce the transformation
formulae
\beq
\bem{c}\LA \delta\vec{v}^-,\hat{t}^-\RA \\ 
         \LA \delta\vec{v}^-,\hat{u}^-\RA\eem = \delta\y'^-
\qmbox{and}\bem{c}\LA \delta\w^-,\hat{t}^-\RA \\ 
         \LA \delta\w^-,\hat{u}^-\RA\eem = 
         (l'^-/c_q)\delta\z'^-
\Leq{chain2}         
from the $(\delta\y'^-,\delta\z'^-)$--coordinates to the 
$(\delta\vec{v}^-,\delta\w^-)$--coordinates. 

After near-collision we have
\beq
\bem{c}\LA \delta\vec{v}^+,\hat{t}^+\RA \\ 
         \LA \delta\vec{v}^+,\hat{u}^+\RA\eem = \delta\y'^+
\qmbox{and}\bem{c}\LA \delta\w^+,\hat{t}^+\RA \\ 
         \LA \delta\w^+,\hat{u}^+\RA\eem = \delta\z'^+ .   
\Leq{chain3}  
The reflection $R_{\hat{s}^--\hat{s}^+}$ by the plane perpendicular to the
vector $\hat{s}^--\hat{s}^+$ transforms the unit vectors as follows:
\[\hat{s}^+=-R_{\hat{s}^+-\hat{s}^-}(\hat{s}^-)\qmbox{,}
\hat{t}^+ = +R_{\hat{s}^+-\hat{s}^-}(\hat{t}^-)\qmbox{and}
\hat{u}^+ = +R_{\hat{s}^+-\hat{s}^-}(\hat{u}^-).\]
Moreover, up to an error term $\cO(1/E)$, we may substitute 
$R_{\hat{s}^+-\hat{s}^-}$ for the reflection
$R_{\vec{v}^+-\vec{v}^-}$ in (\ref{est:dvwp}), since by estimate 
(\ref{ray:family}) the differences 
$|\vec{v}^\pm-\hat{s}^\pm|=\cO(1/E)$ (and since $|\hat{s}^+-\hat{s}^-|>0$).
In particular 
\beq
\Delta\psi = \pi - \alpha(k_{-1},k_0,k_1)+\cO(1/E). 
\Leq{chain4}

Thus putting together (\ref{chain1}), (\ref{chain2}), (\ref{est:dvwp}),
(\ref{Pu:two}),
(\ref{chain3}) and (\ref{chain4}), we obtain
\beqn
\hspace{-9mm}
\bem{c}\delta\y'^+\\ \delta\z'^+ \eem &=&
\frac{4c_q\sin^2(\eh\Delta\psi) E}{-Z_{k_0}}\cdot\NN\\
& &\bem{cc}\idty&\idty\\\idty&\idty\eem 
\bem{cc}\idty&0 \\ 0&(l'^-/c_q)\idty\eem 
\bem{cc}\idty & 0\\ (1-\frac{l^-}{l'^-})\idty &\frac{l^-}{l'^-}\idty\eem
\bem{c}\delta\y^-\\ \delta\z^- \eem \NN\\
& &\hspace{6cm}+\l|\bem{c}\delta\y^-\\ \delta \z^-\eem \ri|\cdot\cO(E^0)\NN\\
&=& \frac{4l^-\cos^2(\eh\alpha(k_{-1},k_0,k_1)) E}{-Z_{k_0}}
\bem{cc}\idty&\idty\\ \idty&\idty\eem \label{chainn}
\bem{c}\delta\y^-\\ \delta\z^- \eem \\
& &\hspace{6cm}
+\l|\bem{c}\delta\y^-\\ \delta \z^-\eem \ri|\cdot\cO(E^0),\NN
\eeqn
as $c_q+l'^--l^-=l^-$.
Inserting (\ref{chainn}) in (\ref{lin:plus}) yields
\beqn
\hspace{-9mm}
\bem{c}\delta\y^+\\ \delta\z^+ \eem &=&
\bem{cc}\idty & 0\\ (1-\frac{c_q}{l^+})\idty &\frac{c_q}{l^+}\idty\eem
\bem{c}\y'^+\\ \z'^+\eem +
\l|\bem{c}\delta\y'^+\\ \delta \z'^+\eem \ri|\cdot\cO(1/E)\NN\\
&=& f(k_{-1},k_0,k_1) E
\bem{cc}\idty&\idty\\ \idty&\idty\eem 
\bem{c}\delta\y^-\\ \delta\z^- \eem +
\l|\bem{c}\delta\y^-\\ \delta \z^-\eem \ri|\cdot\cO(E^0),\NN
\label{chain}
\eeqn
as $l^-= \eh d^{k_{-1},k_0}$, proving (\ref{Tx}).

Strict invariance of the cone field $\cC$
then follows from (\ref{Tx}).
\hfill $\Box$\\[2mm]
We now describe the domains $V_E(\uk)$ and images
$W_E(\uk)$ of the iterated Poincar\'{e}
map with more precision. So let
\beq
V_E(\uk)(\z):= V_E(\uk)\cap 
\l( B_y\times\{\z\}\ri)\qquad (\z\in B_z)
\Leq{VEz}
and similarly 
\[W_E(\uk)(\z):= W_E(\uk)\cap 
\l( B_y\times\{\z\}\ri)\qquad (\z\in B_z)\]
consist of the points with the same configuration space coordinate $\z$.
We also consider them as subsets of $B_y$, forgetting the fixed coordinate
$\z$. 
It turns out that they are diffeomorphic to two-dimensional disks.
We call
\[\diam_y(V_E(\uk)) := \sup_{\z\in B_z}\diam(V_E(\uk)(\z)),\]
resp.\
\[\diam_y(W_E(\uk)) := \sup_{\z\in B_z}\diam(W_E(\uk)(\z))\]
the {\em $y$--diameter} of $V_E(\uk)$, resp.\ $W_E(\uk)$ (which is measured
with the Euclidean metric in the $\y$ coordinates).
\begin{corollary}
For all $\z^-\in B_z$ the map
\beq
V_E(k_{-1},k_0,k_1)(\z^-)\ni\y^-\mapsto \z^+(\y^-,\z^-)\in B_z
\Leq{mapA}
is a diffeomorphism. Thus the domain can be represented as
\[V_E(k_{-1},k_0,k_1)(\z^-) = \{\y^-(\z^-,\z^+)\mid \z^+\in B_z\},\]
with 
\beq
D_1\y^-(\z^-,\z^+)=-\idty+\cO(E^{-1}),
\Leq{diffB1}
\beq
D_2\y^-(\z^-,\z^+)=\frac{\idty}{f(k_{-1},k_0,k_1)E} +\cO(E^{-2}),
\Leq{diffB}
and 
\beq
\y^-(\z^-,\z^+)=-\z^-+\cO(E^{-1}),
\Leq{ymmzm}
the error terms being independent of the parameter $e_0$ in (\ref{def:PoiE}).\\
In particular the $y$--diameter of $V_E(k_{-1},k_0,k_1)$ is of order 
$\cO(e_0/E^{2})$.
\end{corollary}
{\bf Proof.}
From formula (\ref{Tx}) for the derivative of the Poincar\'{e} map
we see (observing that $f(k_{-1},k_0,k_1)\neq 0$) that (\ref{mapA}) is a local
diffeomorphism, and in fact a diffeomorphism onto its image.

Thus by the definitions (\ref{VEM}) and (\ref{WEM}) of domain and image 
(\ref{mapA}) is a diffeo. 
Estimates  (\ref{diffB1}) and (\ref{diffB}) follow by inverting (\ref{Tx}).
Defs.\ (\ref{def:PoiE}), (\ref{VWE}) and (\ref{def:yz})
show that the $y$--diameter of $W_E(k_{0},k_1)$ is of order $\cO(e_0/E)$.
Thus (\ref{diffB}) implies that the $y$--diameter of $V_E(k_{-1},k_0,k_1)$
is smaller by one order. 

By definition (\ref{def:PoiE}), the domain $\Poi^{k,l}_E$ has size $\cO(e_0/E)$,
so that directly we get only $\y^-(\z^-,\z^+)=-\z^-+\cO(e_0/E)$
instead of (\ref{ymmzm}). However, we then see 
from (\ref{diffB1}) that we may enlarge the size parameter $e_0$ without
enlarging the error, and thus prove (\ref{diffB1}).
\hfill $\Box$
\begin{proposition} \label{prop:bVb}
There exist $C_6>1$ and $\delta E>0$, such that for all $E>\Eth$, $m\geq 1$, 
$\uk\in\Adm_0^m$ and $(\y^+_\uk,\z^+_\uk):=
\ov{\cal P}_E(\uk)(\y^-,\z^-)\subset\ov{\Poi}^{k_{m-1},k_m}_E$
the following holds true:
\begin{enumerate}
\item
For $V_E(\uk)(\z^-)$ defined in (\ref{VEz}) the maps
\beq
V_E(\uk)(\z^-)\ni \y^-\mapsto \z^+_\uk(\y^-,\z^-)\in B_z\qquad (\z^-\in B_z)
\Leq{mapAA}
are diffeomorphisms.
\item
\beq
D_1\y^-_\uk(\z^-,\z^+)= -\idty +\cO(E^{-1})\qquad (\z^\pm\in\ov{B}_z)
\Leq{uniform:Lip}
uniformly in $m$.    
\item
For all $\z^+\in \ov{B}_z$   
the vector field on $\ov{B}_z$ given by $\z^-\mapsto \y^-_\uk(\z^-,\z^+)$
points inside the boundary:
\beq
\y^-_\uk(\z^-,\z^+)\cdot \z^- < 0\qquad (\z^-\in\pa\ov{B}_z).
\Leq{ineq:disk}
\item
$V_E(\uk)(\z^-)$ contains a disk of radius 
\beq
C_6^{-1}E^{-1}\prod_{i=1}^{m-1}
\l(|f(k_{i-1},k_i,k_{i+1})|\cdot (E+\delta E)\ri)^{-1}, 
\Leq{contains}
and is contained in a disk of radius 
\beq
C_6 E^{-1}\prod_{i=1}^{m-1}
\l(|f(k_{i-1},k_i,k_{i+1})|\cdot (E-\delta E)\ri)^{-1}. 
\Leq{is:contained}
\item
For $\uk,\ul\in \Adm_0^{m+1}$ with $k_i=l_i$ for $i=0,\ldots,m$ but
$k_{m+1}\neq l_{m+1}$
\beq \hspace{-10mm}
\dist\l(V_E(\uk)),V_E(\ul))\ri)\geq 
C_6^{-1}E^{-1}\prod_{i=1}^{m-1}
\l(|f(k_{i-1},k_i,k_{i+1})|\cdot (E+\delta E)\ri)^{-1}.
\Leq{dist:large}
\end{enumerate}
\end{proposition}
\begin{remark}
Using time reversal, one obtains similar estimates for $W_E(\uk)$.
\end{remark}
{\bf Proof.}
\begin{itemize}
\item $m=1$: \\
Then $V_E(k_{0},k_{1})=B_y\times B_z $ by (\ref{VEByBz}), and the radius
of these two disks is proportional tp $1/E$, see def.\ (\ref{ByBz}).
So 
\[V_E(k_{0},k_{1})(\z)= 
V_E(k_{0},k_{1})\cap \l( B_y\times\{\z\}\ri)=B_y\times\{\z\},\] 
showing (\ref{contains}) and (\ref{is:contained}).

If $\uk=(k_0,k_1,k_2)$ and $\ul=(k_0,k_1,l_2)$ with $l_2\neq k_2$,
then $V_E(\uk)$ and $V_E(\ul)$ are both contained in the Poincar\'e
section $V_E(k_0,k_1)= \Poi^{k_0,k_1}_E$, but are mapped to different  
sections in the next iteration. These have a distance from each other 
that is bounded below by an energy-independent constant.
On the other hand by (\ref{est:dvwp}), scattering by the nucleus $k_1$ 
only leads to an expansion of order $E$, which conversely implies 
\[\dist\l(V_E(k_0,k_1,k_2),V_E(k_0,k_1,l_2)\ri)\geq C_6^{-1}E^{-1}.\]
\item $m\mapsto m+1$: \\
By def.\ (\ref{VEM}) of $V_E(\uk)$
\[V_E(k_{0},\ldots,k_m)(\z)= 
\l( B_y\times\{\z\}\ri)\cap\Po^{-1}(V_E(k_{1},\ldots,k_{m})).\]
Thus the induction step for (\ref{contains}) and (\ref{is:contained})
is provided by 
\[D_2\y(\z,\z^+)=\frac{\idty}{f(k_0,k_1,k_2)E} +\cO(E^{-2})\qquad (\z\in B_z),\]
{\em i.e.} formula (\ref{diffB}). 

The additional control of the $\z$-dependence of the sets $V_E(\uk)$
needed for the distance estimate (\ref{dist:large}) is guaranteed by 
the coupling (\ref{ymmzm}) between the $\y$- and the $\z$-coordinate.
\item
One application of the linearized Poincar\'{e} map (\ref{Tx})
maps the cone field with cones
\[\l\{(\delta\y,\delta\z)\in T_{(\y,\z)} V_E(k_{-1},k_0,k_1)\mid
\delta\y\cdot \delta\z \geq 0\ri\}\]
into the invariant cone field (\ref{cone}).
In particular this applies to vectors of the form $(\delta\z,0)$.

The image of the invariant cone field $\cC$ under time reversal
(using $\TR$ defined in (\ref{time:reversal}))
meets the equation
\[ |\delta\y+\delta\z|\leq \frac{C}{E} |\delta\y-\delta\z|.\]
This then implies the uniform estimate (\ref{uniform:Lip}).
\item 
By (\ref{ymmzm})
\[\diam_y(V_E(k_0,k_1,k_2)(\vec{0}) =\cO(1/E)\quad\mbox{independent of }e_0. \]
Thus $\y^-(0,\z^+)=\cO(1/E)$ independent of $e_0$, too. So the estimate
(\ref{uniform:Lip}) of the first partial derivative gives. 
$$\y^-_\uk(\z^-,\z^+)=\y^-_\uk(\vec{0},\z^+)+\int_0^1 D_1\y^-_\uk(t\z^-,\z^+)\,dt
=-\z^-+\cO(1/E).$$
Now the modulus of $\z^+\in \pa\ov{B}_z$ equals 
$\ev c_p e_0/E$ (see def.\ (\ref{ByBz})), or 
\[\y^-_\uk(\z^-,\z^+)\cdot \z=
\ev c_p e_0(-\ev c_p e_0+\cO(1))/E^2\qquad (\z^-\in\pa\ov{B}_z).\]
Enlarging the parameter $e_0$ if necessary gives (\ref{ineq:disk}).
\item 
To prove the iteration step for the first statement, we first show that 
that map (\ref{mapAA}) is onto: 
\beq
\mbox{For all }\z,\z^-\in B_z\mbox{ there exists
an }\y^-\mbox{ with }\z^+_\uk(\y^-,\z^-)=\z.
\Leq{zu:beweisen}

Namely for $\z^+\in B_z$, $\uk^{I}:=(k_0,k_1,k_2)$ and 
$\uk^{II}:=(k_1,\ldots,k_m)$ we consider the vector field 
\[\vec{F}:B_z\ar\bR^2\qmbox{,} \vec{F}(\z):=
\y^-_{\uk^{II}}(\z,\z^+)-\y^+_{\uk^{I}}(\z^-,\z).\]
If $\vec{F}$ has a zero, there exists an initial point 
$(\y^-,\z^-)\in V_E(\uk)(\z^-)$
meeting (\ref{zu:beweisen}). However, the vector
field $\z\mapsto\y^-_{\uk^{II}}(\z,\z^+)$ points inside the boundary of 
$\ov{B}_z$, whereas by time reversal
\[\y^+_{\uk^{I}}(\z^-,\z)\cdot \z > 0\qquad (\z\in\pa\ov{B}_z).\]
Thus $\vec{F}$ points inwards, too, so that the degree of the vector
 field $\vec{F}$ on the disk
is non-zero, which in turn implies that $\vec{F}(\z)=\vec{0}$ for some $\z$
(see, {\em e.g.}, Hirsch \cite{Hi}, Chapter 5).

Inspection of the derivative (\ref{Tx}) of the Poincar\'{e} map
shows that $\y^-\mapsto \z^+_\uk(\y^-,\z^-)$ is injective, smooth 
and smoothly invertible.
\hfill $\Box$
\end{itemize}
\bigskip
\begin{lemma} \label{lem:y:diam}
For $E>\Eth$, $l<0<r$ and $\uk=(k_l,\ldots,k_r)$ admissible 
\beqno
\lefteqn{\diam(W_E(k_l,\ldots,k_1)\cap V_E(k_0,\ldots,k_r)) \leq}\\
& &4 \l(\diam_y(W_E(k_l,\ldots,k_1)) + \diam_y(V_E(k_0,\ldots,k_r))\ri).
\eeqno
\end{lemma}
{\bf Proof.}
Set $\ul:=(k_l,\ldots,k_1)$ and $\um:=(k_0,\ldots,k_r)$.
Let $A^I\equiv(\y^I,\z^I)$ and $A^{II}\equiv(\y^{II},\z^{II})$
be two points in $W_E(\ul)\cap V_E(\um)$ and 
\[\z^{I,+}:=\z_\um^+(A^I)\in B_z\qmbox{,}\z^{II,-}:= \z_\ul^-(A^{II})\]
the $\z$--components of $\Po(\um)(A^I)$ resp.\ $\Po(\ul)^{-1}(A^{II})$.
The disks 
\[D^I:= \l\{x\in V_E(\um)\cap W_E(\ul)\mid
\z_\um^+(x)= \z^{I,+}\ri\}\]
and
\[D^{II}:= \l\{x\in V_E(\um)\cap W_E(\ul)\mid
\z_\ul^-(x)= \z^{II,-}\ri\}\]
intersect in a (unique) point
$A^{III}\equiv(\y^{III},\z^{III})$ of $\Poi^{k_0,k_1}_E$.

As $\|A^{II}-A^{I}\|\leq\|A^{II}-A^{III}\|+\|A^{III}-A^{I}\|$, 
the lemma follows from the estimates 
\[\diam(D^I)\leq 4\diam_y(W_E(\ul))\qmbox{and}
\diam(D^{II})\leq 4\diam_y(V_E(\um)),\]
and by time reversal symmetry (\ref{time:reversal}) it suffices to prove the 
first one.

$D^I$ is mapped diffeomorphically onto $W_E(\ul)(\z^{I})$ by
\[x\mapsto (\y',\z^{I})
\qmbox{with the unique} \y'\in B_y \mbox{ meeting } 
\z_\ul^-(\y',\z^{I})=\z_\ul^-(x)\]
(uniqueness of $\y'$ follows from Prop.\ \ref{prop:bVb}.1).
By the estimate (\ref{Tx}) on the linearized Poincar\'{e} map
this map increases distances by 
a factor less than 2.
\hfill $\Box$
%
\Section{Symbolic Dynamics} \label{sect:symbol}
%
In this section we analyze the set 
\[\buE = b\cap\SuE\]
of bounded $\Phi^t$--orbits of energy $E$.

To that aim we equip the alphabet $\cS=\{1,\ldots,n\}$
with the discrete topology, and for an interval (see (\ref{interval}))
$I\equiv I_l^r\subset\bZ$ the space
$\cS^I$ with the product topology.
Finally, we introduce the topological subspace
\beq
\Adm_l^r := \{\uk\in\cS^I\mid \uk\,\,{\rm admissible}\}
\Leq{adm:space}
of admissible sequences (that is, $k_{i+1}\neq k_i$), and use 
the abbreviations
\[\Adm:=\Adm_{-\infty}^\infty\qmbox{,}
\Adm^+:=\Adm_{0}^\infty \qmbox{and} \Adm^-:=\Adm_{-\infty}^1.\]

The space $\Adm$ of bi-infinite admissible sequences 
is empty for $n=1$, consists of two points for $n=2$,
and is a Cantor set for $n\geq3$. From now on we assume $n\geq 2$.

The {\em shift} 
\[\sigma:\Adm\ar\Adm\quad, \quad\sigma(\uk)_i := k_{i+1} \quad(i\in\bZ)\]
is a homeomorphism on $\Adm$.
It is well-known that the topology on $\Adm$ is
generated by the metric
\beq
d(\uk,{\underline l}) := \sum_{i\in\bZ}2^{-|i|}\cdot (1-\delta_{k_i,l_i}), \qquad
(\uk,{\underline l}\in \Adm).
\Leq{Sigma:metric:d}

For an admissible sequence 
$\uk^{+} = (k_{0},k_{1},\ldots)\in\Adm^{+}$
we define
\beq
V_E(\uk^{+}) := \bigcap_{m\in\bN} V_E(k_{0},\ldots,k_{m})\subset V_E(k_0,k_1).  
\Leq{eq:stable:c}
Similarly, for an admissible sequence 
$\uk^{-}=(\ldots,k_{0},k_{1})\in\Adm^{-}$ we define
\beq
W_E(\uk^{-}) := \bigcap_{m\in\bN} W_E(k_{-m},\ldots,k_{1})\subset W_E(k_0,k_1).  
\Leq{eq:unstable:c}
\begin{lemma} \label{lem:lipschitz}
For $E>\Eth$ and $\uk\in\Adm^+$ the sets $V_E(\uk^{+})$ and $W_E(\uk^{-})$ 
are the graphs of functions
\[v_E(\uk^{+}):B_z\ar B_y\qmbox{resp.} w_E(\uk^{-}):B_z\ar B_y\]
meeting the Lipschitz estimates
\beq
\l|v_E(\uk^{+})(\z_1)-v_E(\uk^{+})(\z_2)-(\z_2-\z_1)\ri| \leq 
C\frac{|\z_1-\z_2|}{E} \qquad (\z_1,\z_2\in B_z)
\Leq{vk:plus}
resp.
$$\l|w_E(\uk^{-})(\z_1)-w_E(\uk^{-})(\z_2)-(\z_1-\z_2)\ri| \leq 
C\frac{|\z_1-\z_2|}{E} \qquad (\z_1,\z_2\in B_z) $$
for some $C\equiv C(\Eth)$.
\end{lemma} 
{\bf Proof.}
The sets in (\ref{eq:stable:c}) are nested: 
$V_E(k_{0},\ldots,k_{m_{2}})\subset V_E(k_{0},\ldots,k_{m_{1}})$ for
$m_{2}\geq m_{1}$. 
By estimate (\ref{is:contained}), the
$y$-diameter of these sets goes to zero as $m\ar \infty$.

We set $\ul:=(k_{0},\ldots,k_{m})$. For $\z^+\in\pa\ov{B}_z$,
the $\z$-dependence of the curves 
$\{(\y^-_\ul(\z,\z^+),\z)\mid\z\in B_z\}$ in the boundary of $V_E(\ul)$
is controlled by the $m$--uniform estimate (\ref{uniform:Lip}), showing
(\ref{vk:plus}).
Finally, the statements concerning  $W_E(\uk^{-})$ follow by
time reversal (\ref{time:reversal}),  since 
\[W_E(\uk^{-})=\TR (V_E(k_1,k_0,k_{-1},\ldots)).\hspace*{6cm} \Box\]
\bigskip
Let 
\beq
\Lambda_E^{+} := \bigcup_{\uk^{+}\in \Adm^{+}} V_E(\uk^{+})\qmbox{,} 
\Lambda_E^{-} := \bigcup_{\uk^{-}\in \Adm^{-}} W_E(\uk^{-})
\Leq{eq:Lambda:pm}
and 
\[\Lambda_E := \Lambda_E^{+}\cap \Lambda_E^{-}.\]

By restriction, we associate to $\uk\in\Adm$ half-infinite admissible sequences
$\uk^{\pm}\in\Adm^{\pm}$.
Then we define a map $\cF_E: \Adm\ar \Lambda_E$ by
\beq
\cF_E(\uk) := V_E(\uk^{+})\cap W_E(\uk^{-})\qquad(k\in\Adm).
\Leq{cFE} 
Note that in view of Prop.\ \ref{prop:bVb} the disks
$V_E(\uk^{+})$ and $W_E(\uk^{-})$ in $\Poi^{k_{0},k_1}_E$ intersect. 
The Lipschitz estimates of Lemma \ref{lem:lipschitz} imply that
their intersection consists of precisely one point, which
we identify with an element of $\Lambda_E$.

On $\Adm$ we introduced in (\ref{Sigma:metric:d}) the metric $d$.

On the Poincar\'{e} surfaces 
$\Poi_E = \bigcup_{k\neq l} \Poi^{k,l}_E$ 
we use the metric 
\beqno
\dist&:&\Poi_E\times\Poi_E\ar \bR\\
\dist(x,x')&:=&
\l\{ \begin{array}{cl} \sqrt{(\y-\y')^2+(\z-\z')^2} &,\ x,x'\in \Poi^{k,l}_E\\
1 &,\  {\rm otherwise} \end{array} \ri.
\eeqno
based on the $(\y,\z)$--coordinates (\ref{def:yz}) of $x$ and $x'$
(Remark \ref{rem:rotation} showing that $\dist$ is well-defined).
\begin{lemma} \label{shift:lemma} %
There exist $\alpha > 0$ such that for $E>\Eth$, 
$\cF_E$ is an $(\alpha\cdot\ln E)$--H\"older continuous homeomorphism, that is, 
\beq 
\dist(\cF_E(\uk),\cF_E(\ul))\leq C(E)\cdot d^{\alpha\ln E}(\uk,\ul), 
\qquad (\uk,\ul\in\Adm),
\Leq{est:Hoelder}
for some function $C>0$ of the energy,  conjugating the shift
with the restricted Poincar\'{e} map $\Po^{\Lambda} := {\Po}\rstr_{\Lambda_E}$:
\beq
\cF_E\circ \sigma = \Po^{\Lambda}\circ \cF_E.
\Leq{conjugate:map}
\end{lemma} 
{\bf Proof.}
Since $C$ may depend on $E$, we can assume without loss of generality 
that the central blocks of $\uk$ and $\ul$ coincide, {\em i.e.}
$(k_{-1},k_0,k_1)=(l_{-1},l_0,l_1)$. Let $I^r_l\subset\bZ$ be the maximal
interval containing $0$ on which $\uk$ and $\ul$ coincide.
Then by Lemma \ref{lem:y:diam} and (\ref{is:contained})
\beqno
\lefteqn{\dist(\cF_E(\uk),\cF_E(\ul))}\\
&\leq& \diam(W_E(k_l,\ldots,k_1)\cap V_E(k_0,\ldots,k_r)) \\
&\leq &4 \l(\diam_y(W_E(k_l,\ldots,k_1)) + \diam_y(V_E(k_0,\ldots,k_r))\ri)\\
&\leq& 8\frac{C}{E} \l[(\eh f_{\min} E)^{-|l|}+(\eh f_{\min} E)^{1-r}\ri]
\leq 4 f_{\min}^2C(\eh f_{\min} E)^{-\min(1-l,r+1)}\\
&=& 4 f_{\min}^2C\cdot 2^{-\l(\min(1-l,r+1)\frac{\ln(f_{\min} E/2)}{\ln 2}\ri)}
\leq4 f_{\min}^2C\cdot d^{\alpha\ln E}(\uk,\ul)
\eeqno 
with $f_{\min}:=\min_{i\neq j\neq k}|f(i,j,k)|$,
setting $\alpha$ slightly smaller than $1/\ln(2)$.
\hfill $\Box$\\[2mm]
For $n\geq 2$ we denote $\TuR\circ\cF_E:\Adm\ar \bR^{+}$, 
with the return time $\TuR$ defined in 
(\ref{return:time}) by 
$\TuR$, too. Being defined by composition of a smooth map with a
H\"older continuous map, $\TuR$ is H\"older.

The continuous flow on $\buE$ is modelled as follows. 
\begin{definition}\label{defi:suspension}
Given a {\em roof function} $r\in C^0(\Adm,\bR^+)$, we set
\[\Adm_r:= \Adm\times\bR/\sim\]
where $\sim$ is the equivalence relation defined by
\[(\uk,t+r(\uk))\sim(\sigma(\uk),t) \qquad((\uk,t)\in\Adm\times\bR).\]
Then the $r$--{\em suspension flow} is given by
\beq
\sigma_r^t:\Adm_E\ar \Adm_E, \quad [(\uk,s)]\mapsto[(\uk,s+t)]
\qquad (t\in\bR). 
\Leq{def:susp:flow}
\end{definition}
In the interesting case $r=\TuR$ we abbreviate $\Adm_E:=\Adm_{\TuR}$
and $\sigma_E^t:=\sigma_{\TuR}^t$.
\begin{definition} \label{defi:index} %
The set $\CT$ of {\em collision times} of a trajectory $c:I\ar P$
is given by
\[\CT(c) := \l\{t\in I\mid \q\circ c(t)\in\{\s_1,\ldots,\s_n\} \ri\}. \]
The {\em Morse index} of a hyperbolic $T$-periodic trajectory $c:[0,T)\ar P$
is given by
\[\Ind(c):=\sum_{t\in [0,T)\setminus \CT(c)} 
\dim\l(E_u(c(t))\cap{\rm Vert}_{c(t)}\ri) + |\CT(c)|,\]
where $E_{s/u}(c(t))\subset T_{c(t)}P$ denotes the weak (un)stable subspace at 
$c(t)$
and the {\em vertical subspace} 
${\rm Vert}_x\subset T_xP$ at $x\in T^*\Muh\subset P$
is the one annihilated by the linearized configuration space projection
$T_x\q$.
\end{definition}
\begin{remarks}
{\bf 1)} 
The weak (un)stable subspaces $E_{s/u}(c(t))$
are the direct sums of the (2-dim.)
strong (un)stable subspaces and the neutrally stable flow direction 
(see \cite{KH}).
Similar to ${\rm Vert}_{c(t)}$, they are thus 3-dim.\ Lagrangian subspaces
of $T_{c(t)}P$.\\
{\bf 2)} 
Actually $\Ind(c)$ is always finite. We already know from our study of
the (near)-collision process that $\CT(c)$ does not have accumulation points.
On the other hand, the set of times $t$ where the (un)stable subspace 
$E_{s/u}(c(t))$ turns vertical is finite, too, since the kinetic energy term
$\eh\p^{\,2}$ in the Hamiltonian function is a positive quadratic form, see
Duistermaat \cite{Du}.\\
{\bf 3)} 
The additional term $|\CT(c)|$ in the definition of the Morse index is
chosen in the only way that makes that definition invariant under small
perturbations of the flow.
\end{remarks}
For $0<\theta\leq\pi$ we introduced in (\ref{Hul:theta}) 
the pericentric hypersurfaces
\[\Hul(\theta)= \l\{x\in\Hul\l| \
|\q(x)-\s_l|<\frac{|Z_l|}{H(x)\cdot\sin(\theta/4)}\ri.\ri\}
\qquad(l=1,\ldots,n)\]
(with $\Hul$ defined in (\ref{def:Hul})) near the $l$th nucleus.
So the angle parameter $\theta$ fixes the precise meaning of the term
`near-collision'.
\begin{definition} \label{defi:visit} %
The set $\NCT_{\!\!\theta}(c)$ of {\em $\theta$-near-collision times}
of a trajectory $c:I\ar P$ is given by
\[\NCT_{\!\!\theta}(c) := 
\l\{ t\in I\mid c(t)\in \cup_{n=1}^n \Hul(\theta)\ri\}.\]
We say that the trajectory $c$\hspace{2mm}  
{\em  $\theta$--visits the nuclei 
$\uk\in\Adm_l^r$ in succession}
if $\NCT_{\!\!\theta}(c)\neq \emptyset$,
\[r = \l|\NCT_{\!\!\theta}(c)\cap[0,\infty) \ri|  \qmbox{,} 
l = 1- \l|\NCT_{\!\!\theta}(c)\cap(-\infty,0)\ri| \]
and 
\[c(t_i)\in \cH_{k_i}(\theta)\qquad (i\in I_l^r)\] 
for the enumeration
$\NCT_{\!\!\theta}(c)=
\{t_i\mid  i\in I_l^r,\ t_i<t_{i+1},\ t_{0} < 0 \leq t_1\}$
of near-collision times.
\end{definition}
\begin{remarks}
{\bf 1)} $\uk$ is really well-defined, since the hypersurfaces 
$\Hul(\theta)\subset \Hul$ 
do not intersect, and since the flow $\Phi^t$ is transversal to 
$\Hul$ so that $\NCT_{\!\!\theta}(c)$ is discrete.\\
{\bf 2)} As long as we analyze bounded motion, the natural angle will be
$\theta=\amin$, so in that case we do not write it explicitly.

Lateron, in the analysis of the scattering process, 
$\theta$ will be a parameter, since scattering with small angles
shows non-universal features.
\end{remarks} 
\begin{theorem} \label{thm:homeo}
\begin{enumerate}
\item
For an NC configuration and $E>\Eth$ the map 
\[\Ho^T_E:\Adm_E\ar \buE,\qquad [(\uk,s)]\mapsto \Phi^s(\cF_E(\uk))\]
is a H\"{o}lder continuous  homeomorphism conjugating the suspension 
flow with the flow on the set $\buE$ of energy $E$ bound states:
\[\Phi^t\circ\Ho^T_E = \Ho^T_E\circ\sigma_E^t\qquad (t\in\bR).\]
Thus for 
\begin{itemize}
\item
$n=1$ there are no bounded orbits, 
\item
$n=2$ that set consists of one closed orbit: $\buE\cong S^1$
\item
$n\geq3\ $ $\ \buE$ is locally homeomorphic to the product of a Cantor 
set and an interval.
\end{itemize}
\item
All bounded orbits are hyperbolic. 
\item
$\buE$ has measure zero w.r.t.\ Liouville
measure $\lambda_E$ on $\SuE$.
\item
If a $T$-periodic trajectory $c:[0,T)\ar\SuE$ 
visits the nuclei $(k_1,\ldots,k_m)\in \Adm_1^m$ in succession, then
its Morse index equals $\Ind(c)=m$.
\end{enumerate}
\end{theorem}
{\bf Proof.}
{\bf 1)}
By the virial identity (\ref{virial}) all trajectories that leave the 
interaction zone $\IZ$ go to spatial infinity and thus do not belong to
$\buE$. Hence 
\[\buE\subset \cD_E\qmbox{with}\cD_E\subset\SuE\mbox{ defined in (\ref{cDE})}.\]
Lemma \ref{lem:return} tells us that for $\buE$
the Poincar\'{e} return time $\TuE$ is uniformly bounded by $\cO(1/\sqrt{E})$,
so that a fortiori
\[\Phi(\TuE(x),x)\in \Poi_E\qquad(x\in\buE).\]
Thus 
\[\buE=\Phi(\bR,\Lambda_E),\]
showing that $\Ho^T_E:\Adm_E\ar \buE$ is a bijection. 
Using Lemma \ref{shift:lemma}, we see that $\Ho^T_E$ is a 
H\"{o}lder continuous homeomorphism, since the flow $\Phi\rstr\SuE$ 
on the energy shell can locally be straightened out
(see, e.g., \cite{AM}, Thm.\ 2.1.9), and since the Poincar\'{e}
section $\Poi_E$ is transversal to the flow.

Now for $n=1$ the space $\Adm$ of admissible sequences is empty,
and for $n=2$ 
\[\Adm=\{(\ldots12121\ldots), (\ldots21212\ldots)\}
\qmbox{so that}\Adm_E\cong S^1. \]
Finally, for $n\geq 3$ the sequence space $\Adm$ is a Cantor set
(a non-void compact totally disconnected set without isolated points, see,
{\em e.g.} Katok and Hasselblatt \cite{KH}, A 1.)\\
{\bf 2)}
The existence of the strictly invariant cone field $\cC$ 
(defined in (\ref{cone})) on $\Lambda_E$ implies the hyperbolicity of the
flow. The expanding (as well as the contracting) subspace of $T_x\SuE$, 
$x\in\Lambda_E$, is two-dimensional, since it equals 
\[\bigcap_{m\in\bN} D\Po^m\cC(\Po^{-m}(x)).\]
{\bf 3)}
The property of $\buE\subset\SuE$ to have measure zero is defined 
{\em without reference to a measure}
(since in every local chart of the smooth manifold $\SuE$ the image of
Liouville measure $\Lambda_E$ is continuous w.r.t.\ Lebesgue measure, 
see Hirsch \cite{Hi}, Chapt. 3.1).
In fact, applying the Straightening Out Theorem (\cite{AM}, Thm.\ 2.1.9)
it suffices to show that $\Lambda_E$ has measure zero.
This follows for $E>\Eth$ from estimate (\ref{est:Hoelder}).\\
{\bf 4)}
Setting $k_0:=k_m$, we denote by $s_i$ the intersection times
with the Poincar\'{e} surfaces, that is, 
$c(s_i)\in \Poi^{k_{i-1},k_i}_E$ ($i=1,\ldots,m$).

We assume w.l.o.g.\ that the $T$--periodic trajectory $c:[0,T)\ar\SuE$
begins and ends in the Poincar\'{e} surface $\Poi^{k_0,k_1}_E$, 
so that we have $s_i<t_i < s_{i+1}$ and 
$s_1=0$. Additionally we set $s_{m+1}:=T$.
Then we prove that $\Ind(c)=m$ by showing that 
$$\sum_{t\in [s_i,s_{i+1})\setminus \CT(c)}\hspace*{-4mm} 
\dim\l(E_u(c(t))\cap{\rm Vert}_{c(t)}\ri) + |\CT(c)\cap [s_i,s_{i+1})|=1
\qquad(i=1,\ldots,m),$$
which means that every near-collision adds one to the Morse index.
This claim is equivalent to
\beqn 
\lefteqn{\hspace*{-8mm}\sum_{t\in [s_i,s_{i+1})\setminus \{t_i\}} 
\dim\l(E_u(c(t))\cap{\rm Vert}_{c(t)}\ri)=
\l\{\begin{array}{cl}0 & ,\q(t_i)=\s_{k_i}\\ 1 & ,{\rm otherwise}
\end{array}\ri.
\quad(i=1,\ldots,m),}&&\NN\\
&&
\label{maslov}
\eeqn
where we used the coordinates $c(t)\equiv(\p(t),\q(t))$.

First we treat the case $Z_{k_i}<0$ of a repelling Coulomb singularity.
There we know that the trajectory does not touch the 
singularity ($\q(t_i)\neq \s_{k_i}$).

The simplest case is the one where along a neighbourhood
of $\q([s_i,s_{i+1}])$ the potential $V$ coincides with 
the Keplerian potential $-Z_l/|\q-\s_{k_i}|$, and the angular momentum 
$\hat{L}_{k_i}=(\q(t)-\s_{k_i})\times \p(t)$ of
the trajectory relative to the nucleus vanishes.

Then $\p(t_i)=\vec{0}$, which gives a contribution of one in (\ref{maslov}),
since there the flow direction, given by the Hamiltonian vector field 
\[X_H(c(t_i))=(-\nabla V(\q(t_i)),\p(t_i))\in{\rm Vert}_{c(t_i)},\]
is vertical.

However in that case there are no further contributions to the index.
This may be seen as follows.
The intersection of the (three-dimensional) 
weak unstable subspace at $c(s_i)$  with the tangent 
space to the local Poincar\'{e} section $\Poi_E^{k_i,k_{i+1}}$
is a two-dim.\ subspace which 
lies inside the local cone field $\cC(c(s_i))$ defined in (\ref{cone}).
Thus we have $\delta\q\cdot \delta\p\geq 0$ for a variation vector
\beq
(\delta\p,\delta\q)\in E_u(c(s_i))\cap T_{c(s_i)}\Poi_E^{k_i,k_{i+1}}
\subset T_{c(s_i)}P,
\Leq{dp:dq}
and this property is preserved by the forward flow, since 
\[\frac{d}{dt} ( \delta\q(t)\cdot \delta\p(t))=
-\delta\q(t)\cdot D^2V(\q(t))\delta\q(t)+\delta\p(t)\cdot\delta\p(t)\geq 0.\] 
For general repelling potentials and general trajectories $c$ 
(\ref{maslov}) equals one, too, by a continuity argument based on 
estimate (\ref{Tx}).

Now we treat the case $Z_{k_i}>0$, again starting with a
Keplerian potential $-Z_l/|\q-\s_{k_i}|$.
and vanishing angular momentum 
$\hat{L}_{k_i}$. 
Here $\q(t_i)=\s_{k_i}$, so that we have to show that $E_u(c(t))$
does not turn vertical for $t\in [s_i,s_{i+1})\setminus \{t_i\}$.
However, the configuration space trajectories $\q([s_i,t_i])$ and 
$\q([t_i,s_{i+1}])$ are straight lines, and the infinitesimal two-parameter 
family of diverging Kepler hyperbolae corresponding to initial conditions (\ref{dp:dq})
do not have a conjugate point. A more formal way of seeing this is
to use the Jacobi metric $g_E(\q)= (1-V(\q)/E) g(\q)$ 
discussed in (A.\ref{def:jacobi}) of the Appendix.
Then the energy $E$ solution curves correspond to geodesics $c$ in that metric, 
and the linearization of the flow $\Phi_t$ corresponds to the Jacobi equation 
\[\nabla^2\vec{Y}(t)+R_{\dot{c}(t)}\vec{Y}(t)=\vec{0}, \]
see, e.g., \cite{Kli2}.
For a variational vector field $\vec{Y}$ orthonormal to the geodesic velocity
vector $\dot{c}(t)$ 
the self-adjoint curvature operator $R_{\dot{c}(t)}$ has an 
$\vec{Y}$--expectation value equal to the sectional curvature of the 
plane spanned by these two vectors (\cite{Kli2}, Prop.\ 1.11.3).
The sectional curvature of $g_E$ has been calculated in (A.\ref{sect:curv})
of the Appendix.

For the geodesic under consideration the relevant plane contains the direction
$\q-\s_{k_i}$ of the singularity and is thus seen to be 
negative definite (setting $q_3=0$ in (A.\ref{sect:curv})). 

However, negative sectional curvature, together with initial conditions
(\ref{dp:dq}), lead to absence of conjugate points. This proves (\ref{maslov})
for the collision orbit.

In the case of non-vanishing angular momentum $\hat{L}_{k_i}$.
the sectional curvature (A.\ref{sect:curv})) is still negative in 
the plane perpendicular to $\hat{L}_{k_i}$.
However, near the singularity (A.\ref{sect:curv})) becomes positive for 
plane containing the vector $\hat{L}_{k_i}$, leading to a conjugate point.
So (\ref{maslov}) holds for that case, too.

Like in the case $Z_{k_i}<0$, a continuity argument based on 
estimate (\ref{Tx}) shows assertion 4.\ for general potentials $V$
and large $E$.
\hfill $\Box$
%
\Section{Fractal Dimension} \label{sect:fractal}
%
In this section we estimate the fractional dimension of the set 
$\buE$ of  energy $E$ bound states, for $E$ large.
This quantity, being of interest in its own right, governs the
measure of those scattering orbits which have a large time delay. 
\subsection*{Dimensions: Definitions and Elementary Properties} 
Besides the well-known dimension $\dimH$ introduced by Hausdorff and
Besicovitch, there exist several other definitions of the fractional
dimension. Of those we will only consider upper
box-counting dimension $\dimB$, since 
most dimensions take values between $\dimH$ and 
$\dimB$, see Falconer \cite{Fa}.
\begin{definition} 
Let $(X,d)$ be a separable metric space and $U\subset X$, $U\neq\emptyset$.
The {\em diameter} $\diam(U)$ {\em of} $U$ is given by
\[\diam(U) := \sup\{d(x,y)\mid x,y \in U\}\qmbox{and} \diam(\emptyset):=0.\]
For $E\subset X$, $s\geq 0$ and $\delta > 0$ let
\[{\cal H}^{s}_{\delta}(E) := \inf \l\{\sum_{i=1}^{\infty}
(\diam(U_{i}))^{s}\mid E\subset \cup_{i\in\bN} U_{i},
\diam(U_{i})\leq \delta \ri\}\]
(setting $0^0:=1$ except for the case $\diam(\emptyset)^0:=0$).
The {\em Hausdorff $s$-dimensional outer measure of} $E$, ${\cal
H}^{s}(E)$, is then defined by
\beq
{\cal H}^{s}(E) := \sup_{\delta > 0} {\cal H}^{s}_{\delta}(E) =
\lim_{\delta \searrow 0} {\cal H}^{s}_{\delta}(E).
\Leq{def:Hs}
The {\em Hausdorff dimension} \label{Hausdorff:dimension}
of $E$ is given by 
\beq
\dimH (E) := \sup \{s\in\bR \mid {\cal H}^{s}(E) = \infty \}
= \sup \{s\in\bR \mid {\cal H}^{s}(E) >0 \}.
\Leq{def:dimH}
\end{definition}
\bigskip
The important property of ${\cH}^{s}$ is that it defines a measure
on, say, the Borel sets.
Therefore, Hausdorff dimension is not only {\em monotone}, that is
\begin{equation}
E_{1}\subset E_{2}\Rightarrow \dimH(E_{1}) \leq \dimH(E_{2}),
\label{monotone}
\end{equation}
but also {\em countably stable}:
\beq
\dimH(\cup_{i\in\bN} E_{i}) = \sup_{i\in\bN}\ \dimH(E_{i}).
\Leq{s:stable}
\begin{definition} \label{def:dimB}
Let $(X,d)$ be a manifold $X$ with metric $d$ and 
$E\subset X$ be a non-empty bounded subset. 
The (upper) {\em box-counting dimension} (or Minkowski dimension) 
$\dimB(E)$ is given by
\[\dimB(E) := \limsup_{\vep\searrow 0} \frac{\ln(N_\vep(E))}{-\ln(\vep)},\]
where $N_\vep(E)$ is the minimal number of balls of radius $\vep$
needed to cover $E$.
\end{definition}
\bigskip
$\dimB$ is monotone (see (\ref{monotone})), and 
\beq
\dimB(E_{1}\cup E_{2})=\max(\dimB(E_{1}),\dimB(E_{2})),
\Leq{dimB:max}
but it is not countably stable, since $\dimB$ is invariant under closure.
Furthermore, for all (bounded) $E$,
\beq
\dimH(E)\leq \dimB(E).
\Leq{HsB}
\begin{example}
The triadic Cantor set $E\subset \bR$ has
\[\dimH(E) = \dimB(E) = \frac{\ln 2}{\ln 3} = 0.6309\cdots\, .\]
On the other hand, the set $\bQ\subset\bR$ of rational numbers has
dimensions
\[\dimH(\bQ) = 0\qmbox{,} \dimB(\bQ\cap[0,1])=1.\]
\end{example}
We have 
\beq
\dimH(f(E))\leq\dimH(E)\qmbox{and} \dimB(f(E))\leq\dimB(E)
\Leq{map:in}
if the map $f$ is Lipschitz. 
In general both dimensions are {\em not} additive w.r.t.\ cartesian products,
but the inequalities
\beq
\dimH(E_1)+\dimH(E_2)\leq \dimH(E_1\times E_2)\leq \dimH(E_1)+\dimB(E_2)
\Leq{prod:in}
and
\beq
\dimB(E_1\times E_2)\leq \dimB(E_1)+\dimB(E_2)
\Leq{B:prod:in}
hold true, see Mattila \cite{Ma} and Tricot \cite{Tr}.
\subsection*{Lower Estimate for the Hausdorff Dimension}
By (\ref{HsB}) we need a lower estimate for the Hausdorff dimension 
and an upper estimate for the box counting dimension.

Whereas {\em any} choice of coverings by $\vep$-balls leads to
an upper estimate for $\dimB$,
a lower estimate for $\dimH$ is provided by the 
{\em mass distribution principle}:\\[2mm]
{\bf Proposition (see \cite{Fa}).}
Let $\mu$ be a probability measure on $(X,d)$ and $\mu(E)>0$.
Suppose that for some $s\geq0$, $C>0$ and $\delta>0$
\[\mu(U)\leq C\cdot\diam(U)^s\qquad(\diam(U)\leq\delta),\]
Then $\dimH(E)\geq s$.\\[2mm]
{\bf Proof.}
If $E\subset \cup_{i\in \bN} U_i$, then
\[\mu(E)\leq \mu(\cup_{i\in \bN} U_i)\leq \sum_{i=1}^\infty\mu(U_i)\leq 
C\sum_{i=1}^\infty\diam(U_i)^s\]
so that $0<\mu(E)\leq C{\cal H}^{s}_{\delta}(E)$. Then the statement
follows from the second expression for $\dimH$
in (\ref{def:dimH}).\hfill$\Box$
\subsection*{Thermodynamic Formalism for Dimension Estimates}
We base ourselves on the size estimate Prop.\ \ref{prop:bVb} for
$V_E(\uk)$ in terms of the {\em finite} 
geometric data encoded in $f$ (see (\ref{def:f})).

Hence we consider the parameter-depen\-dent 
$(n(n-1)\times n(n-1))$--matrix $M(s)$ with double-indices in
$\{(i,j)\in\cS\times\cS\mid i\neq j \}$ and 
entries
\[M(s)_{i,j;k,l} :=\l\{\begin{array}{cl}
|f(i,k,l)|^{-s} & ,i\neq j,k\neq l\mbox{ and }j=k\\
0 & ,\mbox{otherwise}
\end{array}\ri. . \]
For $s\in\bR$ this is a matrix with non-negative entries, and for 
$n\geq 3$ centres
all entries of $(M(s))^m$ are strictly positive iff $m\geq 3$.
Thus by the Perron-Frobenius (PF) Theorem $M(s)$ has a unique eigenvalue 
$\lambda_{\rm max}(s)$
of largest modulus, which is positive and of multiplicity one, and 
the corresponding eigenvector $v(s)$ can be chosen to have 
strictly positive entries. 

In the case $n=2$ the two eigenvalues of $M(s)$ are given by 
$\pm\lambda_{\rm max}(s)=\pm|f(1,2,1)|^{-s}$.
\begin{lemma}
For all $E>\Eth$ there is a unique solution $d(E):=s$ of the equation
\beq
\lambda_{\rm max}(s)=E^s,
\Leq{ev:eq}
and $d(E)=0$ for $n=2$, whereas for $n\geq 3$
\beq
d(E)=\frac{\ln(n-1)}{\ln(E)}\cdot\l(1+\cO(\frac{1}{\ln(E)})\ri).
\Leq{dE:est}
Finally
\beq
d'(E) =  -\frac{\ln(n-1)}{E(\ln(E))^2}\l(1+ \cO(\frac{1}{\ln(E)})\ri).
\Leq{dpE:est}
\end{lemma}
{\bf Proof.} \\
We conjugate $M(s)/E^s$ with the 
$(n(n-1)\times n(n-1))$--matrix $D(s)$ given by
\[D(s)_{i,j;k,l} :=\l\{\begin{array}{cl}
(d^{i,j})^{s/2} & ,i=k\neq j=l\\
0 & ,\mbox{otherwise}
\end{array}\ri. . \]
Then $\tilde{M}(s,E) := D(s)M(s)D(s)^{-1}/E^s$ has the non-negative entries
\beq
\tilde{M}(s,E)_{i,j;k,l} =\l\{\begin{array}{cl}
(\tilde{f}(i,k,l)E)^{-s} & 
,i\neq j,k\neq l\mbox{ and }j=k\\
0 & ,\mbox{otherwise}
\end{array}\ri.
\Leq{def:tM}
with  
\[\tilde{f}(i,k,l):=\sqrt{d^{k,l}/d^{i,k}}|f(i,k,l)|>0\qquad(i\neq k\neq l).\]
The right Perron-Frobenius eigenvector of
$\tilde{M}(s,E)$ equals 
\[\tilde{v}(s):=D(s)v(s)\qmbox{with eigenvalue} 
\tilde{\lambda}_{\rm max}(s,E):=\lambda_{\rm max}(s)/E^s,\]
and (\ref{ev:eq}) corresponds to the implicit equation 
\beq
\tilde{\lambda}_{\rm max}(d(E),E)=1.
\Leq{imp:eq}
Although $\tilde{f}(l,k,i)=\tilde{f}(i,k,l)$ by def.\ (\ref{def:f}) 
of $f$, for $n>2$ the matrix $\tilde{M}(s,E)$ is non-symmetric
and even non-normal, leading to slightly more complicated estimates.

We denote the left PF eigenvector by $\tilde{w}(s)$, again assuming positivity
of its entries. Then
for $E$ large $\tilde{\lambda}_{\rm max}(s,E)$ is strictly
decreasing in $s$, since 
\beq
D_1 \tilde{\lambda}_{\rm max}(s,E) = 
\frac{\LA\tilde{w}(s),D_1\tilde{M}(s,E)\,\tilde{v}(s)\RA}
{\LA\tilde{w}(s),\tilde{v}(s)\RA}<0,
\Leq{decrease}
$D_i$ denoting the derivative w.r.t.\ the $i$th argument. 

Inequality (\ref{decrease}) follows from the fact that all components of 
$\tilde{v}(s)$ and $\tilde{w}(s)$ are positive, and for $E$ large
all entries of $D_1\tilde{M}(s,E)$
are non-positive, and some are negative. So the l.h.s.\ of the equation
$\tilde{\lambda}_{\rm max}(s,E)=1$ is strictly decreasing in $s$. Furthermore,
\[\tilde{\lambda}_{\rm max}(0,E)=n-1\geq 1\qmbox{and} 
\lim_{s\ar\infty}\tilde{\lambda}_{\rm max}(s,E)=0.\]
This implies that (\ref{ev:eq}) has a unique solution 
$d(E)=0$ for $n=2$ and $d(E)>0$ for $n\geq 3$.

More precisely we observe that the
$n-1$ non-vanishing entries in each row $(i,j)$ of $\tilde{M}(s,E)$
are of the form $(\tilde{f}E)^{-s}$.
Setting 
\[\tilde{f}_{\rm min}:=\min_{i\neq k\neq l}\tilde{f}(i,k,l)\qmbox{and}
\tilde{f}_{\rm max}:=\max_{i\neq k\neq l}\tilde{f}(i,k,l),\]
we note that by a consideration of the largest and the smallest 
components of the eigenvalue equation
$\tilde{\lambda}_{\rm max}(s,E)\tilde{v}(s)=\tilde{M}(s,E)\tilde{v}(s)$
\beq
\frac{\ln(n-1)}{\ln(\tilde{f}_{\rm max}E)}\leq d(E)\leq 
\frac{\ln(n-1)}{\ln(\tilde{f}_{\rm min}E)}.
\Leq{minimax}
This shows the estimate (\ref{dE:est}).
The implicit equation (\ref{imp:eq}) for $d(E)$ shows that 
\beq
d'(E)=\frac{-D_2\tilde{\lambda}_{\rm max}(d(E),E)}
              {D_1\tilde{\lambda}_{\rm max}(d(E),E)}=
  \frac{d(E)}{E\cdot D_1\tilde{\lambda}_{\rm max}(d(E),E)}.            
\Leq{eq:dpE}  
So we have to estimate 
\beq
D_1 \tilde{\lambda}_{\rm max}(d(E),E) = 
\l.\frac{\LA\tilde{w}(s),D_1\tilde{M}(s,E)\,\tilde{v}(s)\RA}
{\LA\tilde{w}(s),\tilde{v}(s)\RA}\ri|_{s=d(E)} 
\Leq{Deins}
more precisely in order to show (\ref{dpE:est}).
The matrix 
\[\cM:= \lim_{E\ar\infty} \tilde{M}(d(E),E)\] 
has the form 
\[\cM_{i,j;k,l} =\l\{\begin{array}{cl}
1/(n-1) & ,i\neq j,k\neq l\mbox{ and }j=k\\ 0 & ,\mbox{otherwise}
\end{array}\ri.,\]
since by (\ref{minimax}) the quotients 
\beqn
\tilde{M}(d(E),E)_{i,k;k,l}/\tilde{M}(d(E),E)_{i',k';k',l'}&=&
\l(\tilde{f}(i',k',l')/\tilde{f}(i,k,l)\ri)^{d(E)}\NN \\
&=& 1+\cO(1/\ln(E))
\label{eq:quotient}
\eeqn
of the non-zero coefficients converge to one, and since the PF eigenvalue 
of $\tilde{M}(d(E),E)$ equals one. We now use the algebraic relation
\[((n-1)\cM)^2+ (n-1)\cM= \cF \qmbox{with} \cF_ {i,j;k,l}:=1\qquad 
(i\neq j,k\neq l),\]
which is approximately met for finite energies in the sense that by
(\ref{eq:quotient})
\beqn
R(E)&:=&\l((n-1)\tilde{M}(d(E),E)\ri)^2+ (n-1)\tilde{M}(d(E),E)\NN\\
&=&\cF+\cO(1/\ln(E)). 
\label{eq:RE}
\eeqn
By (\ref{imp:eq}) the PF eigenvectors fulfil the equations
$$R(E)\tilde{v}(d(E))= n(n-1)\tilde{v}(d(E))\qmbox{,} 
R(E)^t\tilde{w}(d(E))= n(n-1)\tilde{w}(d(E))$$
so that in view of estimate (\ref{eq:RE}) its components are nearly equal:
\beq
\tilde{v}(d(E))_{i,j}=
\frac{\sum_{k\neq l} \tilde{v}(d(E))_{k,l}}{n(n-1)}
\cdot \l(1+\cO(1/\ln(E))\ri) \qquad (i\neq j)
\Leq{nearly}
and similarly for $\tilde{w}$.

Finally, for $i\neq j,k\neq l\mbox{ and }j=k$
\[D_1\tilde{M}(d(E),E)_{i,j;k,l} =
-\ln(|\tilde{f}(i,k,l)E|)/(n-1)\cdot(1+\cO(1/\ln(E))). 
\]
Inserting that estimate and (\ref{nearly}) into (\ref{Deins}) we obtain
\[D_1 \tilde{\lambda}_{\rm max}(d(E),E) = -\ln(E)+\cO(1).\]
Putting that estimate and (\ref{minimax}) into (\ref{eq:dpE})
gives (\ref{dpE:est}).
\hfill $\Box$
\subsection*{Dimensions of the set $\buE$ of Bound States}
We now estimate $\dimB(\buE)$ and $\dimH(\buE)$ by the solution
$d(E)$ of the matrix eigenvalue problem (\ref{ev:eq}).
\begin{theorem} \label{thm:fractal}
For $E>\Eth$ and $n\geq 3$ the Hausdorff dimension $\dimH$ and
the upper box-counting dimension $\dimB$ of the energy $E$
bound states $\buE$ meet the estimates
\beqn
\lefteqn{\hspace{-8mm}1+2d(E)\cdot\l(1-\cO((E\ln E)^{-1})\ri)
\leq\dimH(\buE)\leq}\NN\\
&\leq& \dimB(\buE)
\leq 1+2d(E)\cdot\l(1+\cO((E\ln E)^{-1})\ri)
\label{dim:e}
\eeqn
with the solution $d(E)$ of (\ref{ev:eq}).
In particular they meet the rough estimate
\beq
\dimH(\buE)=1+\frac{2\ln(n-1)}{\ln(E)}+\cO\l((\ln E)^{-2}\ri)=\dimB(\buE)
\Leq{rough:e}
For $n=2$ centres $\dimH(\buE)=\dimB(\buE)=1$.
\end{theorem}
{\bf Proof.} 
The second inequality in (\ref{dim:e}) is the abstract inequality (\ref{HsB}).
Estimate (\ref{rough:e}) follows from (\ref{dim:e}) 
by inserting (\ref{dE:est}).

So it remains to prove the lower bound for $\dimH$ and then
the upper bound for $\dimB$.
For $n=2$ centres $\buE$ consists of one closed orbit, whose dimensions equal
one. So we assume from now on $n\geq3$.
\\[2mm]
{\bf 1)}
Using the constant $\delta E$ from Prop.\ \ref{prop:bVb} and setting
$E_L:=E+\delta E$ for $E>\Eth$, 
estimate (\ref{dpE:est}) shows that 
\[d(E)=d_L\cdot (1-\cO(1/(E\ln E)))
\qmbox{for} d_L:=d(E_L),\]
so that the first inequality in (\ref{dim:e}) follows from an estimate
\beq
\dimH(\buE)\geq 1+2d_L \qquad(E>\Eth).
\Leq{dLb}
We show that this follows from 
\beq
\dimH(\Lambda_E)\geq 2d_L \qquad(E>\Eth).
\Leq{dLL}
First we find $\vep>0$ such that the flow $\Phi:\bR\times P\ar P$, restricted
to $(-\vep,\vep)\times \Poi_E$, is a diffeomorphism onto its image.
Thus
\[U_\vep:=\Phi((-\vep,\vep)\times \Lambda_E)\subset\buE.\] 
satisfies 
\beqn
\dimH(U_\vep) &=& \dimH((-\vep,\vep)\times \Lambda_E)\NN\\
&\leq&\dimH((-\vep,\vep))+ \dimH(\Lambda_E)= 1+\dimH(\Lambda_E)
\label{cc}
\eeqn
using (\ref{map:in}) and (\ref{prod:in}).
By finiteness of the Poincar\'{e} return time $\TuR$ (Lemma \ref{lem:return}),
$\buE$ can be covered by finitely many time translates 
\beq
\buE= \bigcup_{j=0}^{j_{\max}}\Phi(\vep j,U_\vep)\qmbox{with}
j_{\max}:=[\vep^{-1}\sup_{x} \TuR(x)]
\Leq{finite:cover}
so that 
\[\dimH(\buE)=\max_j \dimH(\Phi(\vep j,U_\vep))=\dimH(U_\vep)\]
follows from stability of $\dimH$ and (\ref{map:in}).
So we are reduced to show (\ref{dLL}), by employing the mass distribution
principle. This will be based on a 
$\Po^{\Lambda}$--invariant probability measure $\mu_E$ on $\Lambda_E$, which is the image 
\[\mu_E := \cF_E\, \mu_{\Adm,E}\]
w.r.t.\ (\ref{cFE}) of a measure $\mu_{\Adm,E}$ on $\Adm$. $\mu_{\Adm,E}$
is defined through its values
\beq
\mu_{\Adm,E}(Z(\uk)) := 
\tilde{w}(d_L)_{k_l,k_{l+1}}\cdot 
\prod_{i=l+1}^{r-1} \l(\tilde{f}(k_{i-1},k_i,k_{i+1})E_L\ri)^{-d_L}
\quad(\uk\in \Adm_l^r)
\Leq{mu:cyl}
on the cylinder sets
\[Z(\uk):= \l\{ \uk'\in\Adm\mid k'_i=k_i\, \forall i\in\{l,\ldots,r\}\ri\}
\qquad (l<r\in\bZ,\, \uk\in\Adm_l^r),\]
with the left Perron-Frobenius eigenvector $\tilde{w}(d_L)$ of
$\tilde{M}(d_L,E_L)$. 

Using the $l^1$--normalization 
\beq
\sum_{i\neq k} \tilde{w}(s)_{i,k}=1
\Leq{tv:norm}
of $\tilde{w}(s)$  
and the relations (\ref{def:tM})  and (\ref{imp:eq}) we see that
the definitions (\ref{mu:cyl}) are compatible and define a
$\sigma$--invariant Borel
probability measure $\mu_{\Adm,E}$. By the conjugacy 
(\ref{conjugate:map}) between the 
shift $\sigma$ and the restricted Poincar\'e map
the image measure $\mu_E$ is then indeed 
$\Po^{\Lambda}$--invariant.

We now claim that the mass distribution principle applies.
Namely we have for some $C\equiv C(E)>0$ and $\delta\equiv \delta(E)>0$ 
\beq
\mu_E(U)\leq C\cdot \diam(U)^{2d_L}
\Leq{mass:distr}
for all measurable $U\subset\Lambda_E$ of small diameter $\diam(U)<\delta$.
Instead of general such $U$ we first consider balls 
$B(\delta')\subset\Poi_E$ of small radius $\delta'>0$, 
which are centered at a point
$x\in\Lambda_E$. The symbol sequence $\uk:= \cF_E^{-1}(x)\in\Adm$
of this point projects to the half-infinite sequences 
$\uk^\pm\in\Adm^\pm$ (see (\ref{cFE})). There are unique 
integers $l<0<r$ with
\beq
\qquad
C_6^{-1}E_L^{-1}\prod_{i=l+1}^{0}(\tilde{f}(k_{i-1},k_i,k_{i+1})E_L)^{-1}
>\delta'
\Leq{in:lgd}
\beq
\qmbox{but}
C_6^{-1}E_L^{-1}\prod_{i=l}^{0}(\tilde{f}(k_{i-1},k_i,k_{i+1})E_L)^{-1}
\leq\delta'.
\Leq{in:lld}
respectively 
\beq
C_6^{-1}E_L^{-1}\prod_{i=1}^{r-1}(\tilde{f}(k_{i-1},k_i,k_{i+1})E_L)^{-1}
>\delta'
\Leq{in:rgd}
\beq
\qmbox{but}C_6^{-1}E_L^{-1}\prod_{i=1}^{r}(\tilde{f}(k_{i-1},k_i,k_{i+1})E_L)^{-1}
\leq\delta',
\Leq{in:rld}
By the lower estimates (\ref{in:lgd}) and (\ref{in:rgd}) and 
(\ref{dist:large}) of Prop.\ \ref{prop:bVb}
\[V_E(k_0',\ldots,k_r')\cap B(\delta')\neq \emptyset
\qmbox{only if} (k_0',\ldots,k_r')=(k_0,\ldots,k_r).\]
and similarly for $W_E(k_{l}',\ldots,k_1')$.
However, this implies that 
\[B(\delta')\cap \Lambda_E \subset \cF_E(Z(k_l,\ldots,k_r)).\]
Thus by def.\ (\ref{mu:cyl}) of $\mu_{\Adm,E}$ and the upper estimates
(\ref{in:lld}) and (\ref{in:rld})
\beqno
\mu_E(B(\delta'))&\leq& \mu_{\Adm,E}(Z(k_l,\ldots,k_r))\\
&=&
\tilde{w}(d_L)_{k_l,k_{l+1}}\cdot 
\prod_{i=l+1}^{r-1} \l(\tilde{f}(k_{i-1},k_i,k_{i+1})E_L\ri)^{-d_L}\\
&\leq& c_\alpha \cdot(2\delta')^{2d_L}=c_\alpha \cdot\diam(B(\delta'))^{2d_L},
\eeqno 
the constant
\[c_\alpha:= \l(\max_{a\neq b}\tilde{w}(d_L)_{a,b}\ri)\cdot
\l(C_6 E_L^2\tilde{f}_{\rm max}\ri)^{2d_L}\] 
being independent of $\delta'$ and $\uk$ 
(remembering that $\tilde{f}_{\rm max}=\max_{i\neq k\neq l}\tilde{f}(i,k,l)$).
That is, the balls $B(\delta')$ satisfy the mass distribution estimate
(\ref{mass:distr}).

General sets $U\subset\Poi_E$ of diameter $\eh\delta'$ and
$U\cap\Lambda_E\neq\emptyset$ are subsets of such balls of radius 
$\delta'$ centered at a point in $U\cap\Lambda_E$.
This proves the mass distribution principle (\ref{mass:distr}) 
in general.

We see from (\ref{mass:distr}) that
for any cover $\Lambda_E\subset \cup_i U_i$ with $\diam(U_i)<\delta$,
\[1\leq\mu_E(\Lambda_E)\leq \sum_i \mu_E(U_i)\leq C\sum_i \diam(U_i)^{2d_L}\]
which implies ${\cal H}^{2d_L}(\Lambda_E)\geq 1/C$ and thus 
$\dimH(\Lambda_E)\geq 2d_L$.
\\[2mm]
{\bf 2)}
Similar to the first case, to obtain the upper bound in 
(\ref{dim:e}) for $\dimB(\buE)$, it suffices to show that for 
$\delta E>0$ from Prop.\ \ref{prop:bVb}
\beq
\dimB(\buE)\leq 1+2d_U
\Leq{d:U}
with 
\beq
d_U:=d(E-2\delta E)\qmbox{and}E_U:=E-\delta E.
\Leq{dU:EU}
Also (\ref{d:U}) follows from the estimate 
\beq
\dimB(\Lambda_E)\leq 2d_U,
\Leq{dimU:L}
for then $\dimB(U_\vep)\leq 1+2d_U$ (for the set 
$U_\vep=\Phi((-\vep,\vep)\times \Lambda_E)$) is a consequence of
(\ref{map:in}) and (\ref{B:prod:in}).
Relation (\ref{dimB:max}) may then be used to determine
the box dimension of the covering (\ref{finite:cover}).

We now claim that for $\vep>0$ small we can cover $\Lambda_E$ in the form 
\beq
\Lambda_E \subset \bigcup_{\uk\in I_\vep} \cF_E(Z(\uk))
\Leq{le:cover}
with index set 
\beqno
I_\vep&:=& \l\{ \uk\in \bigcup_{l=-1}^{-\infty}
\bigcup_{r=1}^\infty \Adm_{l}^{r}\,\ri| \\
& &\qquad
\frac{\vep}{16}\geq C_6E^{-1} \prod_{i=l+1}^0 (|f(k_{i-1},k_i,k_{i+1})|\cdot E_U)^{-1}
\geq \frac{\vep}{16 f_{\max} E_U}, \\
& &\l.\qquad
\frac{\vep}{16}\geq C_6E^{-1} \prod_{i=1}^{r-1}(|f(k_{i-1},k_i,k_{i+1})|\cdot E_U)^{-1}
\geq \frac{\vep}{16 f_{\max} E_U}\ri\}.
\eeqno
\begin{itemize}
\item
To prove (\ref{le:cover}), we choose an arbitrary bi-infinite sequence 
$\ul\in\Adm$, and show that there is a $\uk\in I_\vep$ with $\ul\in Z(\uk)$.

Since we have chosen a small $\vep$, the right inequalities in the definition 
of $I_\vep$ are met for $l=-1$ and $r=1$. We find $\uk$ by setting $k_i:=l_i$
and choosing $r$ such that the right inequality would be violated for $r+1$.
Then the left inequality holds. Similarly we choose the minimal
possible $l<0$.
\item
Next we claim that for $\uk\in I_\vep$ the subset $\cF_E(Z(\uk))$
of the Poincar\'{e} section is contained in a ball of radius $\vep$.
This holds true  since 
\[\cF_E(Z(\uk))=W_E(k_l,\ldots,k_1)\cap V_E(k_0,\ldots,k_r),\]
see (\ref{cFE}), so that using Lemma \ref{lem:y:diam} and (\ref{is:contained})
\beqno
\lefteqn{\diam(\cF_E(Z(\uk)))}\\
&\leq& 
4\l(\diam_y(W_E(k_l,\ldots,k_1)) + \diam_y(V_E(k_0,\ldots,k_r))\ri)\\
&\leq& 8\l( 
C_6E^{-1} \prod_{i=l+1}^0 (|f(k_{i-1},k_i,k_{i+1})|\cdot E_U)^{-1} \ri.\\
& &\quad
\l. + C_6E^{-1} \prod_{i=1}^{r-1}(|f(k_{i-1},k_i,k_{i+1})|\cdot E_U)^{-1} \ri)
\ \leq \vep.
\eeqno
\item
Finally, denoting by $b := 1/(16 f_{\max} E_U)$ the constant appearing in the 
definition of $I_\vep$,
the number $N_\vep(\Lambda_E)$ of $\vep$-balls needed to cover
$\Lambda_E$ is bounded by 
\beqno
\lefteqn{\hspace*{-10mm}N_\vep(\Lambda_E)\leq |I_\vep|}\\
&\hspace*{-25mm}\leq& \hspace*{-15mm}
\l|\l\{ \uk\in \bigcup_{l=-1}^{-\infty}
\bigcup_{r=1}^\infty \Adm_{l}^{r}\,\ri|\
 C_6E^{-1} \prod_{i=l+1}^0 (|f(k_{i-1},k_i,k_{i+1})|\cdot E_U)^{-1}
\geq b\vep,\ri.  \\
& & \l.\l.\hspace*{20mm}
 C_6E^{-1} \prod_{i=1}^{r-1}(|f(k_{i-1},k_i,k_{i+1})|\cdot E_U)^{-1}
\geq b\vep\ri\}\ri|\\
&\hspace*{-25mm}\leq& \hspace*{-15mm}
\l|\l\{\l. \uk\in \bigcup_{l=-1}^{-\infty}
\bigcup_{r=1}^\infty \Adm_{l}^{r}\,\ri|\
 C_6^2E^{-2} \prod_{i=l+1}^{r-1} (|f(k_{i-1},k_i,k_{i+1})|\cdot E_U)^{-1}
\geq (b\vep)^2\ri\}\ri|  \\
&\leq& (\vep bE/C_6)^{-2d_U}
\sum_{l=-1}^{-\infty} \sum_{r=1}^\infty \sum_{\uk\in\Adm_l^r}
\prod_{i=l+1}^{r-1} \l(\tilde{f}(k_{i-1},k_i,k_{i+1}) E_U\ri)^{-d_U}\\
&=& (\vep bE/C_6)^{-2d_U}\sum_{m=1}^\infty m\cdot 
\LA\vidty, \tilde{M}(d_U,E_U)^m\vidty\RA,
\eeqno
with the matrix $\tilde{M}(s,E)$ from (\ref{def:tM}), and 
$\vidty\in\bR^{n(n-1)}$ denoting the vector whose components equal 1.
Now for the arguments (\ref{dU:EU})
the Perron-Frobenius eigenvalue $\tilde{\lambda}_{\rm max}(s,E)$
of $\tilde{M}(s,E)$ is smaller than one:
\beqn
\tilde{\lambda}_{\rm max}(d_U,E_U)&=&\frac{\lambda_{\rm max}(d_U)}{E_U^{d_U}}
=\frac{\lambda_{\rm max}(d(E-2\delta E))}{E_U^{d_U}}\NN\\
&<&\frac{\lambda_{\rm max}(d(E-2\delta E))}{(E-2\delta E)^{d_U}}=1.
\label{ineq:PF}
\eeqn
In view of the near-constancy (\ref{nearly})
of the left PF eigenvector $\tilde{v}\equiv\tilde{v}(d(E_U))$ of
$\tilde{M}\equiv\tilde{M}(d_U,E_U)$
we use the inequality
\[\LA\vidty,  \tilde{M}^m \vidty\RA\leq  2n(n-1)
\frac{\LA\tilde{v},  \tilde{M}^m \tilde{v}\RA}{\LA\tilde{v},  \tilde{v}\RA}
\qquad(m\in\bN), \]
to exchange the vector $\vidty$ by $\tilde{v}$:
\beqno
N_\vep(\Lambda_E)
&\leq &(\vep bE/C_6)^{-2d_U}2n(n-1)
\sum_{m=1}^\infty m\cdot 
\frac{\LA\tilde{v},  \tilde{M}(d_U,E_U)^m \tilde{v}\RA}
     {\LA\tilde{v},  \tilde{v}\RA} \\
&=&
(\vep bE/C_6)^{-2d_U}2n(n-1)
\sum_{m=1}^\infty m\cdot \tilde{\lambda}_{\rm max}(d_U,E_U)^m\\
&\leq& C'\vep^{-2d_U}
\eeqno
with finite $C'$, using (\ref{ineq:PF}). 
Inserting this into the Def.\ \ref{def:dimB} 
of the box-counting dimension shows (\ref{dimU:L}).
\hfill $\Box$
\end{itemize}
%
\Section{Topological Entropy} \label{sect:top:ent}
%
We shall now determine the topological entropy of the
flow $\PutE=\Pt\rstr_{\SuE}$ on the energy shell $\SuE$.

The estimate of $\htop(\Phi_E^1)$ is based on Prop.\ \ref{prop:restriction},
stating that the topological entropy of the flow is
determined by its restriction to the set $\buE$ of bound states. 
It is then relatively easy to compute that topological entropy
using symbolic dynamics. 

First we formally introduce the notion of topological entropy of $T:X\ar X$,
following the definition of Bowen. That definition is in a way more
general than others in allowing for non-compact metric spaces $(X,d)$
(see Walters \cite{Wa}).
\begin{definition} \label{top:entropy}
Let $(X,d)$ be a metric space and $T:X\ar X$ a uniformly continuous
map. Then for $m\in \bN$, $\vep>0$, a subset $F\subset X$ is said to
$(m,\vep)$-{\em span} a compact $K\subset X$ if 
\[\forall x\in K\, \exists y\in F : d_{m}(x,y)\leq \vep, \]
with the metric $d_{m}(x,y) := \max_{0\leq i\leq m-1}
d(T^{i}x,T^{i}y)$.

Let $r_{m}(\vep,T,K)$ be the smallest cardinality of an
$(m,\vep)$-spanning set $F$ of $K$,
\[h_{r}(\vep,T,K) := \limsup_{m\ar\infty} \frac{1}{m} \ln(r_{m}(\vep,T,K)),\]
and $\htop(T,K) := \lim_{\vep\ar 0} h_{r}(\vep,T,K)$.

Then the {\em topological entropy of} $T$ is
\[\htop(T) := \sup \l\{ \htop(T,K)\mid K\subset X \mbox{ compact} \ri\}.\]
\label{def:htop}
\end{definition}
Note that $h_{r}(\vep,T,K)$ is monotonically increasing in $\vep$, so
that its limit $\htop(T,K)$ exists (it can be $\infty$).

If $d'$ is a second metric on $X$ {\em uniformly equivalent to} $d$, that
is, if the maps
\[{\rm Id}:(X,d)\ar(X,d')\qmbox{and}{\rm Id}:(X,d')\ar(X,d) \]
are both uniformly continuous, then the topological entropies
$\htop(T,d)$ and $\htop(T,d')$ coincide, see \cite{Wa}.
Therefore, for compact
metrisable spaces $X$ the topological entropy only depends on the
topology generated by the metric.

We are to estimate the topological entropy for the time-one flow $T :=
\pptE{1}$ on the energy shell $X := \SuE$, 
and for any metric on $\SoE$ uniformly equivalent outside the
compact region projecting to the interaction zone to the metric
induced by the Euclidean metric of $T^*(\bR^3-\IZ)$. 
$T$ is uniformly continuous since the motion is
asymptotically free.

In our case we have 
\[\htop(\pptE{-1}) = \htop(\pptE{+1}),\]
although in general $\htop(\pptE{-1},K) \neq
\htop(\pptE{+1},K)$, since we have a symmetry w.r.t.\ time
reversal.

Prop.\ \ref{prop:restriction} below will show the importance of the
set $\buE$ of energy $E$ bound states.
\begin{definition} \label{defi:non:wandering} 
The {\em non-wandering set}\, 
$\Omega(\Phi)$ of a continuous flow
$\Pt:X\ar X$ is given by 
\[\Omega(\Phi) := \l\{ x\in X\l| \forall\mbox{ neighb. } U\ni x,
\forall T>0: U\cap \l( \bigcup_{t\geq T}\Pt(U) \ri) \neq
\emptyset \ri.\ri\}.\] 
\end{definition}
For non-compact spaces $X$ the topological entropy is in general larger than
the topological entropy of the restriction to the non-wandering set.
Nevertheless, in our case we have

\begin{proposition} \label{prop:restriction} %
For $E>\Eth$, the non-wandering set equals
\[\Omega(\Phi_E) = \buE,\] 
and  the topological entropy of the flow on the energy shell is
determined by the bound states, that is, $\htop(\Phi_E^{1}) =0$ for $n=1$ 
centre and
\beq
\htop(\Phi_E^{1}) = \htop(\Phi_E^{1}\rstr_{\buE})
\Leq{restriction}
for $n\geq 2$.
\end{proposition}
{\bf Proof.}
 The proof parallels the one of Lemma 7.4 of \cite{KK}.
In particular $\htop(\Phi_E^{1}) =0$ for $n=1$, since then
by Thm.\ \ref{thm:homeo} $\buE=\emptyset$.
\hfill $\Box$
\begin{theorem} \label{thm:top:ent}
{\bf 1)} 
For $E>\Eth$ and $n=1$ or $n=2$ centres the topological entropy
equals $\htop(\Phi_E^{1})=0$.\\
{\bf 2)} For $n\geq 3$ centres
\beq
\htop(\Phi_E^{1})=\htopi\cdot\sqrt{2E}\cdot
\l(1+\frac{\ln(E)}{E}C_\htop + \cO(1/E)\ri).
\Leq{htop:est}
Here $\htopi$
is the unique solution (with eigenvector $\vec{v}$) 
of the largest eigenvalue problem
$\lambda(s)=1$ for the one parameter family of $n\times n$ 
Perron-Frobenius matrices $\cM(s)$,
\beq
\cM(s)_{i,j}:= \l\{
\begin{array}{cl} \exp(-s d^{i,j}) & i\neq j\\ 0 & i=j\end{array}
\ri. ,
\Leq{matrix:M}
\[C_\htop := \frac{\LA \vec{v},\vec{Z}\RA\LA\vidty,
\vec{v}\RA}{2\LA \vec{v},D\vec{v}\RA}\qmbox{with}
\vec{Z}:=(Z_1,\ldots,Z_n)\,,\,\vidty:=(1,\ldots,1)\]
and the matrix $D$ of distances $d^{i,j}$.

{\bf 3)} 
Let $N_E(T)$ be the number of closed $\Phi_E^t$--orbits of period
smaller than $T$. Then there exists a ($E$-independent)
constant $C>1$ so that for $T$ large 
\beq
C^{-1} \frac{\exp(\htop T)}{\htop T} \leq \NuE(T) \leq
C \frac{\exp(\htop T)}{\htop T},\qquad(E>\Eth)
\Leq{schachtel}
with $\htop\equiv \htop(\Phi^{1}_{E})$.
\end{theorem}
\begin{remark}
For configurations which form an equilateral triangle
or a tetrahedron, the trivial estimate
\[\htopi\in \l[ \frac{\ln(n-1)}{\dmax}, \frac{\ln(n-1)}{\dmin} \ri]\]
for the PF eigenvalue problem of the matrices $\cM(s)$ is sharp.
\end{remark}
{\bf Proof.}
{\bf 1)} 
For $n=1$, by Prop.\ \ref{prop:restriction}
the topological entropy is zero for $E$ large.
By Prop.\ \ref{prop:restriction}, for $n\geq 2$
we need only estimate $\htop(\Phi_E^{1}\rstr_{\buE})$.

For $n=2$, we know from Thm.\ \ref{thm:homeo} 
that $\buE$ consists of only
one closed orbit. Therefore, by Thm.\ 7.14 of \cite{Wa}, the topological 
entropy $\htop(\Phi_E^{1})=\htop(\Phi_E^{1}\rstr_{\buE})=0$,
too.\\[2mm]
{\bf 2)} 
The interesting case left is $n\geq3$.
Since $\htop$ is a conjugacy invariant (Thm.\ 7.2 of \cite{Wa}),
we have by Thm.\ \ref{thm:homeo}
\beq
\htop(\Phi_E^{1}\rstr_{\buE})=\htop(\sigma_E^1)
\Leq{conjugacy}
for the time-one $\TuR$--suspension flow $\sigma_E^1$ on $\Adm_E$ (see
(\ref{def:susp:flow})).

Generally, given a roof function $r:\Adm\ar\bR^+$ 
(see Def.\ \ref{defi:suspension}),
any $\sigma$-invariant probability measure $\mu$ on $\Adm$
induces a probability measure $\mu_r$
on $\Adm_r$ which is invariant under the $r$-suspension flow
$\sigma^t_r$, namely the one obtained by normalization
of the measure $\mu\times dx$ on $\Adm\times \bR$, restricted to
the fundamental domain
$\{(\uk,t)\mid 0\leq t<r(\uk)\}\cong \Adm_r$.
Conversely, any $\sigma^t_r$-invariant probability measure
can be obtained in that way.

Furthermore, the topological entropy is the supremum of the
$KS$-entropies over all {\em ergodic} measures, and the same is 
generally true
for the {\em topological pressure} 
\beqno
P(T,f) &:=& \sup \l\{\l.h_{\mu}(T)+\int_X f\, d\mu\ \ri|
 \ \mu \ T-\mbox{invariant prob.\
measure}\ri\}\\
&=& \sup \l\{ \l.
h_{\mu}(T)+\int_X f\, d\mu \ \ri|
 \ \mu \ \mbox{ergodic w.r.t.\ }T\ri\}
\eeqno
of a continuous map $T:X\ar X$ of a compact metrisable space $X$ and 
$f\in C(X,\bR)$, see Cor.\ 9.10.1 of \cite{Wa}.

Finally, any $\sigma^t_r$-ergodic measure $\mu_r$ on 
$\Adm_r$ comes from a $\sigma$-ergodic measure $\mu$ on $\Adm$.
 
So we can apply Abramov's formula
\beq
h_{\mu_r}(\sigma_r^{1}) = 
\l(\int_{\Adm} r\, d\mu \ri)^{-1} h_{\mu}(\sigma),
\Leq{true}
(which is valid for any ergodic measure $\mu$ on $\Adm$, see
Theorem 2.1 of Chap.\ 3 in Sinai \cite{Sin}), to obtain
the formula
\beq
\htop(\sigma^1_r)=\sup \l\{ \l.
h_{\mu}(\sigma)\l/\int_{\Adm} r\, d\mu \ \ri.\ri|
 \ \mu \ \mbox{ergodic}\ri\}.
\Leq{htop:r}
Eq.\ (\ref{htop:r}) implies the scaling behaviour
\beq
\htop(\sigma^1_{\lambda
r})=\lambda^{-1}\htop(\sigma^1_r)\qquad(\lambda\in\bR^+)
\Leq{scaling}
and the inequality 
\beq
\htop(\sigma^1_{r_1})\geq\htop(\sigma^1_{r_2})\qquad (r_1\leq r_2).
\Leq{htop:ineq}

Using Lemma \ref{ret:time:est} for estimating the Poincar\'{e} 
return time $\TuR$, we obtain for $C$ large 
\[\l(1-C/E\ri)\cdot r(E)\leq \sqrt{2E}\cdot \TuR \leq 
\l(1+C/E\ri) \cdot r(E)
\qquad(E>\Eth),\]
with 
$$r(E)\in C(\Adm,\bR^+)\qmbox{,}r(E)(\uk):=
\eh\l(d^{k_{-1},k_0}+d^{k_0,k_1}-Z_{k_0}\frac{\ln E}{E}\ri)
\quad(\uk\in\Adm).$$
So by applying (\ref{scaling}) and (\ref{htop:ineq}), we get
\beq
\htop(\sigma^1_E)= \htop(\sigma^1_{r(E)})\cdot\sqrt{2E}\cdot
\l(1+\cO\l(1/E\ri)\ri).
\Leq{htop:reduction}
The claim (\ref{htop:est}) follows from (\ref{htop:reduction}) and
(\ref{conjugacy}), if we can prove
\beq
\htop(\sigma^1_{r(E)})= 
\htopi\cdot\l(1+\frac{\ln(E)}{E}C_\htop + \cO(1/E)\ri).
\Leq{htop:htop}
This is done by reformulating it as a question about the topological
pressure $P(\sigma,-s\cdot r(E))$ for the shift $\sigma$. 
We know that for $s=0$
\[P(\sigma,0)=\htop(\sigma)=\ln(n-1)>0,\]
since $\ln(n-1)$ is the largest eigenvalue of the transition matrix
of the shift, and since $n\geq 3$. On the other hand, for $E\geq\Eth$
the values $r(E)(\uk)$
are larger than $\dmin/2$, so that the roof function 
$r\equiv r(E)\in C(\Adm,\bR^+)$
has mean $\int_X r\, d\mu\geq \dmin/2$ w.r.t.\ to any probability measure 
$\mu$. Thus $P(\sigma,-s\cdot r)\leq 0$ for $s\geq2\ln(n-1)/\dmin$, and
there exists an $s>0$ with
\beq
P(\sigma,-s\cdot r)=0.
\Leq{pressure:zero}
By convexity of the map $f\mapsto P(\sigma,f)$ (Thm.\ 9.7(v) of 
\cite{Wa}) this $s$ is unique.

Moreover the roof function $r$ is locally
constant, so that there exists a unique equilibrium state $\mu$ for 
(\ref{pressure:zero})  (see \cite{Ru}), {\em i.e.} 
\[h_\mu(\sigma)= s\int_X r\, d\mu\qmbox{and}
h_\nu(\sigma)< s\int_X r\, d\nu\qmbox{for}\nu\neq\mu.\]
Thus $\htop(\sigma^1_r)=s$, using (\ref{htop:r}) (and
ergodicity of $\mu$, which follows from uniqueness, see Thm.\ 9.13 of
\cite{Wa}).
Since $\sigma:\Adm\ar\Adm$ is an expansive homeomorphism, by Thm.\ 9.6 
of \cite{Wa} we can use the formula
\[P(\sigma,f)=\lim_{m\ar\infty}\frac{1}{m} \ln(p_m(\sigma,f,\alpha)),\]
with the generating partition
\[\alpha=\{A_1,\ldots,A_n\}\qmbox{with atoms} A_l:=\{\uk\in\Adm\mid k_0=l\} \]
and
$$p_m(\sigma,f,\alpha):=
\inf\l.\l\{\sum_{B\in\beta} \sup_{\uk\in
B}\,\exp(S_m f(\uk))\,\ri|\, \beta\mbox{ finite subcover of }
\bigvee_{i=0}^{m-1}\sigma^{-i}\alpha\ri\}.
$$
Here $S_m f:=\sum_{i=0}^{m-1}f\circ \sigma^i$.
Now for $\delta:=\ln(E)/E$
\beqn
\lefteqn{\exp(-s \cdot S_m r(E)(\uk))=}\NN\\
& &\hspace{-1.5cm}
\exp\l(-\frac{s}{2}
\l(d^{k_{-1},k_0}-d^{k_{m-1},k_{m}}+(Z_{k_m}-Z_{k_0})\eh\delta\ri)\ri)
\cdot\prod_{i=0}^{m-1} \cM(s,\delta)_{k_i,k_{i+1}},\label{eq:Smr:M}
\eeqn
with the $n\times n$--matrix $\cM(s,\delta)$ given by
\[\cM(s,\delta)_{i,j}:= \l\{
\begin{array}{cl} 
\exp\l(-s \l(d^{i,j}-\frac{\delta}{4}(Z_{i}+Z_{j})\ri)\ri)&,i\neq j\\ 
0 &, i=j\end{array}\ri. \]

So the function is constant on the atoms of the partition
$\beta:=\bigvee_{i=-1}^{m}T^{-i}\alpha$. 
The first factor in (\ref{eq:Smr:M}) is bounded from below and from
above, uniformly in $m$.
The PF property of 
the symmetric matrix $\cM(s)$  then implies 
that $P(\sigma,-s\cdot r(E))$ is equal to the logarithm
of the largest eigenvalue $\lambda(s,\delta)$ of $\cM(s,\delta)$. 

We now do perturbation theory around $\cM(s,0)=\cM(s)$, with 
$\cM(s)$ defined in (\ref{matrix:M}),
and obtain from the condition $\lambda(s(\delta),\delta)=1$
\[s'(\delta)=
-\frac{D_2\lambda(s(\delta),\delta)}{D_1\lambda(s(\delta),\delta)}.\]
Since the PF eigenvalue is isolated, 
the derivatives exist, and since $\cM(s,\delta)$ is symmetric, they 
equal at $\delta=0$
\[D_1\lambda(s(0),0)=\LA\vec{v},D_1\cM(s(0),0)\,\vec{v}\RA=
-\LA\vec{v},D\vec{v}\RA\]
and
\[D_2\lambda(s(0),0)=\LA\vec{v},D_2\cM(s(0),0)\,\vec{v}\RA=
\eh\htopi\LA\vec{v},\vec{Z}\RA\LA\vidty,\vec{v}\RA.\]
So
\[s(\delta)=s(0)+s'(0)\delta +\cO(\delta^2) =
\htopi\cdot(1+C_\htop\delta+\cO(\delta^2)),\]
proving (\ref{htop:htop}).\\[2mm]
{\bf 3)} 
Finally, for fixed energy $E>\Eth$ the estimate 
(\ref{schachtel}) for the number of periodic
orbits follows from Thm.\ 2 of the article \cite{PP}
by Parry and Pollicott (or alternatively by 
arguments based on the Renewal Theorem):
If the Axiom A flow $\sigma_E^t$ is topologically weak-mixing, 
then
\[N_E(T)\sim \frac{\exp(\htop T)}{\htop T}.\]
Otherwise all periods are known to
be integral multiples of some $T_0>$, and the formula
\beq
N_E(T)\sim \frac{\htop T_0}{1-\exp(-\htop T_0)}
\frac{\exp(\htop [\frac{T}{T_0}]T_0)}{\htop T}
\Leq{NE:sim}
follows directly from Thm.\ 2 of \cite{PP}.
Here $T_0$ is the largest positive eigenperiod, which is bounded by
\[T_0\leq\inf_{\uk\in\Adm} T_E(\uk) =\cO(1/\sqrt{E}).\]
As we just proved that $\htop=\cO(\sqrt{E})$,
$\htop\cdot T_0$ has an $E$-independent upper bound, and we may thus 
use an $E$-independent $C>1$ in (\ref{schachtel}) to bound the r.h.s.\
of (\ref{NE:sim}) from above and from below.
\hfill
$\Box$
\begin{remarks}
{\bf 1)}
Observe that, up to a small error term, $\htop(\Phi_E^{1})$
is independent of the charges $Z_l$ of the centres.
Effectively an attracting potential speeds up the particle a bit, 
which then leads
to the $\vec{Z}$--dependent correction term in (\ref{htop:est}).
\\[2mm]
{\bf 2)}
It is natural to ask whether one may improve (\ref{schachtel}) to show
\beq
\NuE(T) \sim \frac{\exp(\htop T)}{\htop T}.
\Leq{question}
However, 
the return time estimate of Lemma \ref{ret:time:est} for $T_E(x)$ is, 
up to the relative
order $\cO(1/E)$, independent of the point $x\in V_E(k_{-1},k_0,k_1)$. 

So we cannot exclude by that estimate that all return times are equal to a
constant $T_R$
for a symmetric configuration (e.g., for 
an equilateral triangle or tetrahedron with equal charges $Z_l$).
But by formula (\ref{NE:sim}) in that case $\NuE(T)$ is not asymptotic 
to {\em any} smooth function, so that (\ref{question}) does not hold.
In this sense our statement is optimal, given the, already quite precise,
return time estimate. 

Any eventual improvement of Lemma \ref{ret:time:est} must be 
dependent of the additional smooth component $W$ of the 
potential $V$, and thus be complicated.

Although the iterated Poincar\'{e} {\em map} is certainly mixing 
for the measure 
of maximal entropy, by the above argument we cannot decide whether or 
not the {\em flow} is mixing.
\end{remarks}
%
\Section{Characterization of the Scattering Orbits} \label{sect:sc}
%
Up to now we were mainly concerned with the bound states
$\buE\subset\SuE$. However, the topological entropy
analyzed in the last section is an example for a quantity which, though
a priori depending on the dynamics on the {\em whole} 
energy shell $\SuE$, is determined by that subset of 
Liouville measure zero.

When we now consider the scattering states $s_E\subset\SuE$,
again their description will be based on
symbolic dynamics of the {\em bound} states.

Our concrete question will be to enumerate the scattering orbits
with given energy $E>0$ and asymptotic directions $\htheta^\pm\in S^2$, 
that is, the subset
\[(H,\hat{p}^-,\hat{p}^+)^{-1}(E,\htheta^-,\htheta^+)\]
of phase space (noticing that the asymptotic
$\hat{p}^\pm$ defined in (\ref{as:dir}) are constant on orbits).

Before we come to the case of several centres, we first consider 
the simple
\begin{example}{\bf: Keplerian motion}.\\
Without loss of generality $V(\q)=-\frac{Z}{|\q|}$. 
There are no undeflected orbits, {i.e.} only asymptotic directions
\[(\htheta^-,\htheta^+)\in (S^2\times S^2)\setminus \Delta,
\qmbox{with diagonal}
\Delta:=\{(\htheta,\htheta)\mid\htheta\in S^2\}\]
occur.
The pericentric time $T$ is a smooth function on the positive energy
part 
\[P_+ := \{ x\in P \mid H(x)>0 \} \]
of phase space, see
(\ref{peric:time}) and the definition of the manifold $P$.
In fact in that case the map
\[(H,\hat{p}^-,\hat{p}^+,T)\rstr_{P_+}:P_+\ar 
\bR^+\times\l((S^2\times S^2)\setminus \Delta\ri)\times\bR\]
is a diffeomorphism. In particular for $E>0$ and asymptotic
directions $\htheta^-\neq\htheta^+$
there is exactly one Kepler hyperbola.
\end{example} 
If we allow for an additional smooth potential $W$, we need
to exclude a whole neighbourhood of the diagonal 
$\Delta\subset S^2\times S^2$, since even for large energies
forward scattering may be dominated by $W$ and not by the 
Coulomb potential.

Already in the case of bounded orbits we introduced the NC condition
of Def.\ \ref{defi:noncol} which excluded collinear configurations
of nuclei. Similarly we need to exclude certain asymptotic directions 
if we want to obtain statements which are independent of the precise form 
of the potential $V$. 

Thus for $\vartheta\in (0,\pi]$ and
\[\thmin:S^2\ar[0,\pi]\qquad \thmin(\htheta):=\min_{1\leq i\neq k\leq n}
\sphericalangle(\htheta,\hat{s}^{i,k})\] 
(letting $\thmin(\htheta):=\pi$ for $n=1$ nucleus)
we restrict our interest to the asymptotic directions in
$$\AD(\vartheta):=\l\{(\htheta^-,\htheta^+)\in S^2\times S^2\,\l|\,
\min\l( \sphericalangle(\htheta^-,\htheta^+),
\thmin(\htheta^-),\thmin(\htheta^+)\ri)
> \vartheta\ri.\ri\},$$
and  the ($\Phi^t$-invariant) subset 
\[s_E(\vartheta):= \{x\in
s_E\mid(\hat{p}^-(x),\hat{p}^+(x))\in\AD(\vartheta) \}\]
of scattering states,
thus excluding near-forward scattering
and scattering from or to any direction near an axis 
through two nuclei.

For $n\geq2$ centres we have bound states which influence the
scattering trajectories in their vicinity. As the high energy bound
states can be described by symbolic dynamics, we introduce
the set 
\[\cW:= \bigcup_{k\in\bN} \Adm_1^k\]
of (non-empty) {\em admissible words} to enumerate 
the scattering states, see (\ref{adm:space}).

Every orbit within $s_E(\vartheta)$ will be uniquely characterized by its 
asymptotic
directions and by the sequence of its near-collisions with the nuclei
($\theta$--visits in the sense of Def.\ \ref{defi:visit}).
Conversely every admissible word will be shown to occur for some orbit.

Namely we set the angle parameter controlling the asymptotic directions
$\AD$ equal to
\beq
\vartheta\equiv \vartheta(E):=\min(C_7/\sqrt{E},\amin/2),
\Leq{theta:var}
with $C_7$ to be fixed in Thm.\ \ref{thm:sc}.

Similarly, the angle parameter controlling the near-collisions is fixed by 
\beq
\theta\equiv \theta(E):=\min(C_9/\sqrt{E},\amin/2)\qmbox{with} C_9:=4c_1. 
\Leq{theta:non:var}
In Lemma \ref{lem:lin} the parameter $c_1$ 
controlled what the regime of hard scattering. 
It will be fixed in Thm.\ \ref{thm:sc}, too.

For $x\in s_E$ let 
$w(x):=\emptyset$ if $\NCT_{\!\!\theta}(x)=\emptyset$. Otherwise
the orbit through $x$ enters the interaction zone, {\em i.e.} 
intersects the hypersurface 
$\pa\cD_E^-$ at a unique point $x'$, and we set 
\[w(x):= \uk(x')\in\cW,\]
where the trajectory $t\mapsto\Phi^t(x')$ $\theta$--visits 
the centres in succession $\uk$
(we start at $x'$ instead of $x$ since we want 
the $\theta$---visits to occur for positive times so
that $\uk\in\Adm_1^r$ for some $r$).
\begin{lemma} \label{lem:total:defl}
For $C_8>0$ large and $E>\Eth$
every orbit in $s_E(C_8/E)$ intersects the hypersurfaces 
$\pa\cD_E^\pm$ exactly once, so that the entrance, 
resp.\ exit times $T^\pm_E$ on $s_E(C_8/E)$ are uniquely defined by  
$\Phi(T^\pm_E(x),x)\in \pa\cD_E^\pm$.\\
The restrictions of $T^+_E$ and $T^-_E$ to $s_E(C_8/E)$ are smooth functions.
\end{lemma}
{\bf Proof.}
If a scattering orbit in $s_E$ does not meet the interaction zone, 
then by the virial estimate (\ref{qp}) 
there is a unique point $x\equiv \pq$ on this 
orbit with $\LA \p,\q\RA =0$.
By (\ref{V:small}) the speed is bounded below by $|\p|\geq \sqrt{E}$. Thus the
direction $\hat{p}=\p/|\p|$ differs from the asymptotic directions
$\hat{p}^\pm(x)$ (defined in (\ref{as:dir})) only by
\[\sphericalangle(\hat{p},\hat{p}^\pm(x)) = \cO(1/E), \]
using (\ref{O:pL}).
Choosing a large enough constant $C_8$, we see that
\beq
\sphericalangle (\hat{p} , \hat{p}^-(x)) +
  \sphericalangle (\hat{p} , \hat{p}^+(x)) < C_8/E
\qquad\mbox{if } E > \Eth. 
\Leq{total:defl:out}
Thus the orbits in $s_E(C_8/E)$ intersect $\pa\cD_E^\pm$.

For these scattering states we know from (\ref{qp}) that there is a unique
entrance resp.\ exit time.
As long as the orbits intersect the $C^\infty$-hypersurfaces $\pa\cD_E^\pm$
transversally, smoothness of $T^\pm_E$ follows from smoothness of
the flow $\Phi^t$. But this transversality can be enforced by further enlarging
$C_8$.
\hfill $\Box$
\begin{theorem} \label{thm:sc}
We assume that $V$ satisfies the decay estimates (\ref{smooth}). Then
for $C_7>0$ large in (\ref{theta:var}), and $E>\Eth$
the map  
\beqno
{\rm Diff}_E: s_E(\vartheta(E))&\ar& 
\AD(\vartheta(E))\times \cW\times \bR\\
x &\mapsto& (\hat{p}^-(x),\hat{p}^+(x),w(x),T^-_E(x))
\eeqno
is a diffeomorphism.
\end{theorem}
\begin{remarks}
{\bf 1)}
In particular the {\em orbits} in $s_E(\vartheta(E))$ are uniquely characterized by
their asymptotic directions and the sequence of near-collisions, and
any sequence is realized by an orbit.\\
{\bf 2)}
In this section and Sect.\ \ref{sect:diff:cross} we exclude from our 
consideration small cones of aperture $\cO(1/\sqrt{E})$ around the
axes through two nuclei.

This is indeed necessary, since there {\em non-universal phenomena occur}.
Also, the order $\cO(1/\sqrt{E})$ of aperture is optimal, as we show
now by example.

We consider the purely Coulombic
two-centre problem (see Appendix B) with $\s_1:=\bsm 1\\0\\0\esm$, 
$\s_2:=\bsm -1\\0\\0\esm$, $Z_1>0$ and $Z_2=-Z_1$.

We claim that for energy $E>\Eth$ there is no orbit 
colliding with $\s_1$ and having scattering angles $\theta(\infty)$
w.r.t.\ to the negative 1--axis smaller than
$\eh \sqrt{|Z_2|/E}$, whereas there are such orbits for 
$\theta(\infty) > 2\sqrt{|Z_2|/E}$.

Such an orbit would lie in a plane containing $\s_1$ and $\s_2$, say the 
$1-2$--plane. By reflection symmetry 
we need only consider the part of the orbit after collision, and we denote by 
$\tau$ the time of its last intersection with the 
$q_1\equiv0$--plane.
Then $V(\q(\tau))=0$ so that $p_1(\tau)=-\sqrt{2(E-p_2^2(\tau))}$. 
So at time $\tau$ the angle of the particle direction
$\p(\tau)$ with the negative $1$--axis equals
\[\theta(\tau)=\arcsin(|p_2(\tau)|/\sqrt{2E}).\]
As we are interested in small
scattering angles, we may assume $|p_1(\tau)|>|p_2(\tau)|$.
Moreover
\[q_2(\tau) = \frac{p_2(\tau)}{-p_1(\tau)} +\cO(1/E),\]
as follows from the estimates of Sect.\ \ref{sect:long:paths}.
Thus the angular momentum w.r.t.\  $\s_2$ equals
\[|\vec{L}_2|=|(\q(\tau)-\s_2)\times \p(\tau)|=|p_2(\tau)-q_2(\tau)p_1(\tau)|=
2|p_2(\tau)|+\cO(1/E).\]
As the scattering angle $\Delta\theta$ in a Coulombic potential of charge
$Z_2$ meets the relation
\[\frac{\sin\Delta\theta}{1+\cos\Delta\theta}= \frac{|Z_2|}{\sqrt{2E}|\vec{L}_2|},\]
we obtain for $1/\sqrt{E}\ll |p_2(\tau)|\ll \sqrt{E} $ the total scattering
angle
\[\theta(\infty)=\theta(\tau)+\Delta\theta\cdot(1+o(1))= 
(2E)^{-1/2}\l[|p_2| + \frac{|Z_2|}{|p_2|} \ri] \cdot(1+o(1)),\]
which is minimized by $|p_2|\approx \sqrt{|Z_2|}$, with value
\[\theta(\infty)=\sqrt{2|Z_2|/E}\cdot(1+o(1)).\]
\end{remarks}
\bigskip
{\bf Proof of Thm.\ \ref{thm:sc}.}\\
As remarked before Lemma \ref{lem:lin}, there are two regimes for 
Coulomb scattering, hard and soft scattering. 
A large choice of the constant $C_7$ in Thm.\ \ref{thm:sc} separates
the two regimes.\\ 
{\bf 1)}
In order to prove the theorem, we thus first have to
fix the constants $c_1,c_2$ which appear in Part 1) and 2) of Lemma \ref{lem:lin}, 
for $\delta=\eh$. The lemma describes the linearization of the Kepler 
Transformation, but its estimates are also valid for the motion in the 
potential $V$, using estimate (\ref{Pu:one}) and (\ref{Pu:two}).

We are interested in the effect of that linearized scattering transformation 
on cone fields of the form 
\beq
\cC(a):=\l\{\bsm \delta\vec{u} \\ \delta\vec{v}\esm\mid 
|\delta\vec{u}-\delta\vec{v}|\leq 
a |\delta\vec{u}+\delta\vec{v}|\ri\}\qquad(a>0).
\Leq{CF}
So we have $\cC(b)\subset\cC(a)$ for $b<a$
(note for comparison that in (\ref{cone}) we used the 
cone field with energy-dependent $a=C/E$).

\begin{itemize}
\item
$\vep>0$ is chosen small enough 
so that the (according to Lemma \ref{lem:C1} nearly free)
linearized flow within $\IZ(c_q)$ maps 
the cone $\cC(1+2\vep)$ into $\cC(2)$, and 
the cone $\cC(\eh)$ into $\cC(1-2\vep)$.

For such a small $\vep>0$ we then choose 
\item
$c_1>0$ large enough so that the linearized scattering transformation
maps the cone $\cC(2)$ into $\cC(\eh)$
\item
$c_2>0$ small enough so that the linearized scattering transformation
maps the cone $\cC(1+\vep)$ into $\cC(1+2\vep)$.
\end {itemize}
All these choices can indeed be made (apply the relevant
unperturbed matrices $\bsm R&R\\ R&R\esm$ and 
$\bsm \idty&0\\ s\idty&\idty\esm$ from (\ref{est:dvwp}) and (\ref{C1:near}),
resp.\  (\ref{dev:from:free})
on the vector $\bsm \delta\vec{u} \\ \delta\vec{v}\esm$). 

By fixing the constants $c_1$ and $c_2$,
we have defined soft and hard scattering.\\[2mm]
{\bf 2)} Now all $\theta$--visits in the sense of Def.\ \ref{defi:visit}
(meaning that 
the orbit of energy $E$ locally intersects the Poincar\'{e}
surface $\Hul(\theta(E))$),
are hard scattering events.

Namely, as $\theta(E)=4c_1/\sqrt{E}$ (see (\ref{theta:non:var})),
the scattering angle $\Delta\psi$ inside the ball $B_l(c_q)$ is
bounded below by
\[\Delta\psi\geq\theta(E)/4=c_1/\sqrt{E},\]
using (\ref{0:0}). It thus meets the criterion of Lemma \ref{lem:lin} for
hard scattering.\\[2mm]
{\bf 3)} 
We show now that, for large enough $C_7$ in (\ref{theta:var}) 
({\em i.e.} by excluding large cones of asymptotic directions
$\hat{p}^\pm(x)$ around the axes through two centres), 
${\rm Diff}_E$ is well-defined, that is, that
every scattering state $x\in s_E(\vartheta(E))$\ 
$\theta$--visits at least one nucleus, so that
$\NCT_{\!\!\theta}(x)\neq \emptyset$. 

In view of 
Lemma  \ref{lem:total:defl} the orbit through $x\in s_E(\vartheta(E))$ meets
the interaction zone. 
If the directions at the times $T^\pm_E(x)$ of entrance resp.\  exit
are denoted by $\hat{p}_i$ resp.\ $\hat{p}_o$, then
like in (\ref{total:defl:out}), we have the estimate 
\[\sphericalangle (\hat{p}_i , \hat{p}^-(x)) +
  \sphericalangle (\hat{p}_o , \hat{p}^+(x)) =\cO(1/E)\qquad
(E > \Eth). \]
Since $\sphericalangle (\hat{p}^-(x) , \hat{p}^+(x)) > \vartheta(E)$,
we conclude that for large enough $C_7$ 
the change of direction inside the interaction zone is at least
\beq
\sphericalangle (\hat{p}_i ,\hat{p}_o) > \eh\vartheta(E).
\Leq{g:dr:vi}
Lemma \ref{lem:C1} then tells us that the trajectory enters at least one 
ball of (energy-independent) radius $c_q$ around a singularity $\s_l$.

We count the number $N$ of such occurrences. If $N\geq3$, then
Prop.\ \ref{propo:drei:R} says that the orbit $\amin/2$--visits 
at least one nucleus. Comparing with the definition (\ref{theta:non:var})
of $\theta$, this is a $\theta$--visit.

But the same holds for $N\leq2$, since otherwise\\
$\bullet$
the changes of directions inside the $N$ $c_q$--balls are 
each smaller than $\theta$ (Lemma \ref{lem:defl}), and \\
$\bullet$ 
the $N+1$ components of the trajectory inside the interaction zone but 
outside the $c_q$-balls each contribute at most with $C/E$ to the change of 
direction (Lemma \ref{lem:C1}).\\
$\bullet$ 
Adding these contributions and using the definition (\ref{theta:var}) of 
$\theta$, we would get 
\[\sphericalangle (\hat{p}_i ,\hat{p}_o)<\frac{3C}{E}+\theta(E)= 
\frac{3C}{E}+\frac{C_9}{\sqrt{E}}
<\eh \frac{C_7}{\sqrt{E}}=\eh\vartheta(E)
\qquad\mbox{for } \Eth\mbox{ large},\]
in contradiction with (\ref{g:dr:vi}) if
\[C_7\geq4 C_9.\]
{\bf 4)}
We want to avoid intermediate scattering events with angles $\Delta\Psi$
meeting 
\beq
\frac{c_2}{\sqrt{E}}\leq \Delta\Psi\leq \frac{c_1}{\sqrt{E}}.
\Leq{intermediate}
This can be done, too, by choosing $C_7$ in (\ref{theta:var}) 
large enough.

From 3) we know that the orbit $\theta$--visits the centres in 
some succession $\uk=(k_1,\ldots,k_r)$. There can be at most {\em one} visit
of some ball $B_{k_0}(c_q)$ before and some ball $B_{k_{r+1}}(c_q)$ 
after this sequence of near-collisions
(Prop.\ \ref{propo:drei:R}).

Moreover, for large $C_7$, these visits, if they occur, are {\em soft} 
scattering events in the sense of Lemma \ref{lem:lin}. We need only prove 
this for $B_{k_0}(c_q)$, as the result for $B_{k_{r+1}}(c_q)$
then follows by time reversal.

So we consider the half-orbit $\Phi((-\infty,t_1],x)$ through $x$
with $t_1=\min \NCT_\theta(x)$, see Def.\ \ref{defi:visit}.
This half-orbit ends in the Poincar\'{e} surface
near $\s_{k_1}$, namely $\Phi(t_1,x)\in\cH_{k_1}(\theta)$.
More precisely, by def.\ (\ref{def:Hul}) the configuration
space distance is bounded by
\[|\q(t_1,x)-\s_{k_1}|\leq \frac{4\Zmax}{C_9\sqrt{E}}. \]
The half-trajectory $\q((-\infty,t_1],x)$ consists of three types of segments:
\begin{enumerate}
\item 
The segment outside the interaction zone $\IZ$
Here the total change of direction is of order $\cO(1/E)$
(Thm.\ \ref{thm:both:moeller}).
\item 
The (one or two) segments in $\IZ(c_q)$.
Here, too the total change of direction is of order $\cO(1/E)$
(Lemma \ref{lem:C1}).
\item 
The segment inside the ball $B_{k_0}(c_q)$, perhaps empty.
We know that the half-orbit has no hard collision here, {\em i.e.}
it misses $\cH_{k_0}(\theta)$ 
(otherwise $\uk$ would begin with $k_0$ instead of $k_1)$).
Thus the total change of direction inside $B_{k_0}(c_q)$ is bounded by
\[\Delta\psi<\theta=\frac{C_9}{\sqrt{E}}.\]
\end{enumerate}
Summing these contributions, for $\Eth$ large the total change of direction
on the half-orbit is bounded by
\beq
\sphericalangle(\hat{p}^-(x),\p(t))< \frac{2C_9}{\sqrt{E}}
\qquad(-\infty<t\leq\tau),
\Leq{wi1}
$\tau$ being the time where the half orbit enters
the ball $B_{k_1}(c_q)$.

On the other hand by our assumption $x\in s_E(\vartheta(E))$
\beq
\sphericalangle\l(\hat{p}^-(x),\hat{s}^{k_1,k_0}\ri)>\frac{C_7}{\sqrt{E}}.
\Leq{wi2}
$(\vec{v}^-,\vec{w}^-) = \l(\p(\tau,x)/\sqrt{2E},(\q(\tau,x)-\s_{k_1})/c_q\ri)$ 
are the coordinates (\ref{vw:coord})
of the end point $\Phi(\tau,x)$ of the half-orbit. 
By Lemma \ref{lem:large:dev} 
\[\sphericalangle(\vec{v}^-,-\vec{w}^-)<\frac{\pi\cdot\Eth}{\ev C_9\sqrt{E}},\]
so that with (\ref{wi1}) the relative position $\vec{w}^-$ meets
\[\sphericalangle(\vec{w}^-,-\hat{p}^-(x))\leq
\sphericalangle(\hat{p}^-(x),\vec{v}^-)+\sphericalangle(\vec{v}^-,-\vec{w}^-)
< \l.\l(2C_9+\frac{4\pi\cdot\Eth}{C_9}\ri)\ri/\sqrt{E}.\]
Together with (\ref{wi1}) and (\ref{wi2})
this means that for $C_7$ large the distance of $\q(t,x)$ from
the line $L$ through $\s_{k_0}$ and $\s_{k_1}$ is bounded by 
\[\dist(\q(t,x),L)\geq \eh\frac{C_7}{\sqrt{E}}\cdot |\q(t,x)-\s_{k_1}|
\qquad(-\infty<t\leq\tau).\]
But then 
\[\dist(\q(t,x),\s_{k_0})\geq \eh\frac{C_7}{\sqrt{E}}\cdot \dmin/2
\qquad(-\infty<t\leq\tau).\]
so that there is at most {\em soft} scattering inside $B_{k_0}(c_q)$,
if $C_7$ is large.
Thus there are no intermediate scattering events with angles 
meeting (\ref{intermediate}).
\\[2mm]
{\bf 5)}
The last remark implies that the symbolic sequence $x\mapsto w(x)$
is locally constant on $s_E(\vartheta(E))$.
So by Lemma \ref{lem:total:defl} and Thm.\ \ref{thm:smooth:moeller}
the map ${\rm Diff}_E$ is smooth.\\[2mm]
{\bf 6)}
Our next task is to show that ${\rm Diff}_E$ is onto. 
It suffices to find for given data 
\beq
(\htheta^-,\htheta^+,\uk)\in \AD(\vartheta)\times \cW
\Leq{data}
an orbit in $\SuE$ with these asymptotic directions and $\theta$--visits.

We proceed in a way similar to the construction of bounded orbits in 
Sections \ref{sect:ICF} and  \ref{sect:symbol}.
If the symbol sequence equals $\uk=(k_1,\ldots,k_r)$, then we
erect incoming resp.\ outgoing Poincar\'{e} surfaces $V(\htheta^-,k_1)$
resp.\ $W(k_r,\htheta^+)$ with 
\beqno
V(\htheta,l)&:=& 
\l\{x\equiv\pq\in \cD \l| |\q-\s_l|=c_q,
\frac{|\q-\q_i|}{c_q}< \frac{C_{10}}{\sqrt{H(x)}}\ri.\ri.,\\
& & \hspace*{4.5cm}\l.
\l|\frac{\p}{\sqrt{2H(x)}}\times\htheta\ri| < 
\frac{2C_{10}}{\sqrt{H(x)}}\ri\},
\eeqno
\beqno
W(l,\htheta)&:=&  
\l\{x\equiv\pq\in \cD \l|\ |\q-\s_l|=c_q,
\frac{|\q-\q_i|}{c_q}< \frac{C_{10}}{\sqrt{H(x)}}\ri.\ri.,\\
& & \hspace*{4.5cm}\l.
\l|\frac{\p}{\sqrt{2H(x)}}\times\htheta\ri| < 
\frac{2C_{10}}{\sqrt{H(x)}} \ri\},
\eeqno
$\htheta\in S^2$, 
near $\q_i:=\s_l-c_q\htheta$ resp.\ near $\q_o:=\s_l+c_q\htheta$
(compare this definition with the one of $\Poi^{k,l}_E$ in (\ref{def:PoiE})).

For all $\q^-$ in the configuration space projection of $V(\htheta^-,k_1)$
and all $\q^+$ in the configuration space projection of $W(k_r,\htheta^+)$
there exist $\p^\pm$ with $(\p^-,\q^-)\in V(\htheta^-,k_1)$,
$(\p^+,\q^+)\in W(k_r,\htheta^+)$ such that the orbit through $(\p^-,\q^-)$
$\theta$--visits the centres in succession $\uk$ and then meets $(\p^+,\q^+)$.
The proof of this assertion
uses the estimates of Sect.\ \ref{sect:ICF},
slightly modifying Prop.\ \ref{prop:bVb} 

We set $C_{10} := C_7/8$ and choose a large value of $C_7$.
As $\thmin(\htheta^\pm)> \vartheta(E)=C_7/\sqrt{E}$, 
by what we have proven in part 4), 
$V(\htheta^-,k_1)$ did not $\theta$--visit a nucleus in the past, and
the orbits through
$W(k_r,\htheta^+)$ will not $\theta$--visit a nucleus in the future.
Instead by Prop.\ \ref{propo:drei:R}, 
they leave the interaction zone in time $\mp\cO(\sqrt{E})$.

Using estimate (\ref{wi1} with $t:=\tau$, we see that for every 
$\q^-$ there is at least one $\p^-$ such that the orbit through
$(\p^-,\q^-)\in V(\htheta^-,k_1)$ has asyptotic direction
\[\lim_{t\ar-\infty} \hat{p}(t,(\p^-,\q^-))=\htheta^-.\]
By standard arguments this family of orbits contains one
$\theta$--visiting the centres in succession $\uk$ and then having limit
\[\lim_{t\ar\infty} \hat{p}(t,(\p^-,\q^-))=\htheta^+.\]
So ${\rm Diff}_E$ is onto.\\[2mm]
{\bf 7)}
Cone field estimates based on Part 1) now show that ${\rm Diff}_E$ is  
one to one and  smoothly invertible.    
Specifically we show that there exist invariant cone fields $\cC(a)$ 
(see (\ref{CF})) along the orbits constructed in Part 6).

The idea is simply that there is at least one hard scattering along the
scattering orbit, making the family of energy $E$ configuration space
trajectories with
initial asymptotic direction $\htheta^-$ divergent for large positive times. 

So for given $\htheta^-$ we consider the Lagrange manifold of points $x\in s_E$
having asymptotic direction $\hat{p}^-(x)=\htheta^-$. 
We claim that for $\Eth$ large the 
tangent space of this submanifold at $(\p^-,\q^-)\in V(\htheta^-,k_1)$
is contained in the (wide) cone $\cC(2)$. Namely 
the half-orbit ending in $(\p^-,\q^-)$ consist of at most four segments:
\begin{enumerate}
\item 
The segment outside the interaction zone $\IZ$, controlled by
Thm.\ \ref{thm:both:moeller}. For large $\Eth$ this is in a cone 
$\cC(1+\vep/2)$. 
\item 
The segment in $\IZ(c_q)$ before entering the ball $B_{k_0}(c_q)$, 
controlled by Lemma \ref{lem:C1}. Here the linearized Lagrange manifold
is contained in $\cC(1+\vep)$.
\item 
At most one segment in $B_{k_0}(c_q)$, controlled by the choice of 
the soft scattering constant $c_2$ in 
Part 1) of this proof. 
So after this soft collision, we are in $\cC(1+2\vep)$.
\item 
The segment in $\IZ(c_q)$ after leaving $B_{k_0}(c_q)$ and
before entering $B_{k_1}(c_q)$. 
By Part 1) here the linearized Lagrange manifold
is contained in $\cC(2)$.
\end{enumerate}
Our choice of the hard scattering constant $c_1$
implies that after leaving $B_{k_1}(c_q)$ we are in the narrow cone
$\cC(\eh)$.
At least the same is true after leaving $B_{k_r}(c_q)$, using the cone field 
estimates of Sect.\ \ref{sect:ICF}. 
Finally, we consider the positive half orbit starting at 
$(\p^+,\q^+)\in W(k_r,\htheta^+)$ and get a similar sequence 
\[\cC(\eh)\to \cC(1-2\vep)\to \cC(1-\vep)\to\cC(1-\eh\vep)\]
of cone fields, using Part 1) again.
For large energies this shows uniqueness of the orbit with data (\ref{data}).
\hfill $\Box$   
%
\Section{The Differential Cross Section} \label{sect:diff:cross}
%
The scattering transformation $S$ defined in (\ref{eq:S}) contains
complete information on the scattering process. As we have seen in
the previous sections, it exhibits many aspects of irregularity if
$n\geq 3$. Nevertheless, the scattering transformation is not
directly accessible in a (classical) scattering experiment.

Firstly, one typically cannot fix the initial angular momentum of the
test particle. 
Secondly (unlike in a quantum mechanical setting where interference effects
exist) it is hard to measure time delay. 

What {\em is} accessible is the differential cross section
$\DCSP$. Informally speaking, 
this is the (density of the) number of particles per second scattered in the 
final direction
$\htheta^+$, assuming a uniform flux of one particle 
per second and unit area of incoming particles of energy $E$ and
initial direction $\htheta^-$. 

One could expect to see some trace of irregularity in the 
differential cross section, and in fact for all systems considered up
to now numerical calculations of the cross section indicated the
existence of so-called rainbow singularities 
on a Cantor set of angles, see \cite{Ec,EJ,Ga,Sm,Te2}. 

For the simplest case of a Kepler potential we obtain the 
so-called Rutherford cross section, see (\ref{Ruth}) below. 
It is remarkable that the differential cross sections for the cases
$n=2$ and $n\geq 3$ turn out to be 
very similar to the Rutherford cross section 
(see (\ref{comp:Ru}) for the statement, and
Figure 12.1 of \cite{KK} for a numerical plot for $d=2$ dimensions). 
So the complicated
structure of the time delay and the scattering orbits
is not reflected in the cross section. 

The reason for that discrepancy is, roughly speaking, the following.
For $d=2$ the deflection functions $L^-\mapsto \htheta^+(E,\htheta^-,L^-)$
(depicted in \cite{KK}, Figs.\ 10.2-10.3) 
are strictly monotonic
w.r.t.\ the initial angular momentum $L^-$. It is clear from the definition
of $\DCS$ that extrema of the deflection function lead to
singularities in the differential cross section. Since there are no
extrema (except for the degenerate situation 
$L^{-}\ar\pm\infty$), we have a nonsingular $\DCS$ (except for the
forward direction).

Before stating our theorems, we shall recall a mathematically correct
definition of cross section.
In the physics literature the cross section is sometimes introduced
as a function, whereas it really is a measure. The difference is of some
importance
because in general that cross section measure is not absolutely
continuous w.r.t.\ Haar measure. In our context, we shall {\em
show} that under certain conditions the cross section measure {\em is}
absolutely continuous if one excludes the forward direction and the 
directions near $\hat{s}^{k,l}$, and that
the Radon-Nikodym derivative, {\em i.e.}\ the differential cross section,
is smooth.

We denote by 
\[\vec{a}:P_+\ar\bR^3\qmbox{,}\pq\mapsto \q-\LA \q,\hat{p}\RA\hat{p}\]
the {\em impact parameter}. For  
phase space points $x\in P_+$ projecting to a singularity $\s_l$
or with $x\equiv\pq=(\vec{0},\q)$ the impact parameter is defined by, say 
$\vec{a}(x):=\vec{0}$.

Now on the $\pm$--scattering states $s^\pm\subset\P_+$
the asymptotic impact parameters 
\[\vec{a}^\pm: s^\pm\ar\bR^3 \qmbox{,}
\vec{a}^\pm:=\lim_{t\ar\pm\infty} \vec{a}\circ\Pt\]
are well-defined continuous functions, since $\vec{a}\pq$ is continuous
outside the interaction zone (that is, for $|\q|\geq \Rvir(H\pq)$), and
there 
\beq
\vec{a}\pq =\frac{\p\times \vec{L}\pq}{2(H\pq-V\pq)}.
\Leq{eq:impact}
By Thm.\ \ref{thm:both:moeller} the r.h.s.\ of (\ref{eq:impact}) has a limit on 
scattering orbits.  
 
For energy $E>0$ we consider the maps 
\[A_E^\pm:s_E\ar T^*S^2\qmbox{,} x\mapsto \l(\vec{a}^\pm(x),\hat{p}^\pm(x)\ri)\]
from the set of $\pm$ scattering states of energy $E$ to their asymptotic data.
These can really be considered as points in the cotangent bundle 
\[N:=T^*S^2\]
of the two-sphere, since 
\begin{itemize}
\item
$|\hat{p}^\pm(x)|=1$,
\item
$\LA\vec{a}^\pm(x),\hat{p}^\pm(x)\RA=0$, and
\item
the bilinear map $\pq\mapsto \LA\p,\q\RA$ is the
natural pairing between vectors and co-vectors (no Riemannian metric involved).
\end{itemize}
The cotangent bundle $N$ carries the canonical symplectic two-form $\omega_N$ 
and the volume form
$\lambda_N:=\eh\omega_N\wedge\omega_N$.

The maps $A_E^\pm$ are constant on orbits, and for $y\in N$ the preimages
$(A_E^\pm)^{-1}(y)$ consist of at most one orbit.

\begin{proposition} \label{prop:MAA}
For energy $E>0$ the set $A_E^\pm(s_E)\subset N$ 
of asymptotic data of the scattering states is open, and its
complement 
$N\setminus A_E^\pm(s_E)$ is a compact set of $\lambda_N$-measure zero.

If $V$ satisfies the decay estimates (\ref{smooth}), then the map 
\beq
M_E:A_E^-(s_E)\ar A_E^+(s_E)\qmbox{,} 
A_E^-(x)\mapsto A_E^+(x)\qquad (x\in s_E)
\Leq{ME}
from initial to final data for energy $E$ is a smooth canonical transformation.
\end{proposition}
{\bf Proof.}
For energy $E>0$ and
radius $r>\Rvir(E)$ we consider the smooth Poincar\'{e} surfaces
\[U^\pm_{E,r}:=\{\pq\in \SuE\mid |\q|=r, \pm\LA\p,\q\RA>0\}.\]
These four-manifolds are transversal to the flow as 
$\{H,\q^{\,2}\}=\eh\{\p^{\,2},\q^{\,2}\}=-2\LA\p,\q\RA\neq0$.
According to Lemma 8.2 of McDuff and Salamon \cite{DS},
they are symplectic submanifolds of our phase space $P$, and
the Poincar\'{e} section map 
\beq
U^-_{E,r}\cap s_E\ar U^+_{E,r}\cap s_E
\Leq{UEr}
which send $x$ to the unique intersection point
of the orbit $\Phi(\bR,x)$ with $U^+_{E,r}$
is a symplectomorphism.

We use on $U^\pm_{E,r}$ the coordinates $(\vec{a},\hat{p})$, which
are maps 
\[\kappa^\pm_{E,r}:U^\pm_{E,r}\ar N.\] 
Indeed they are diffeomorphisms onto their common image 
\[N_r:=\{(\vec{a}',\hat{p})\in N\mid |\vec{a}'|<r\}.\]
Thus the map (\ref{UEr}) induces a diffeomorphism
\[M_{E,r}: N_r\cap\kappa^-_{E,r}(s_E) \ar N_r\cap\kappa^+_{E,r}(s_E)\]
which converges pointwisely to (\ref{ME}) as $r\ar\infty$, using 
Thm.\ \ref{thm:both:moeller}.
With assumption (\ref{smooth}), smoothness of $M_E$ is then implied by 
Thm.\ \ref{thm:smooth:moeller}.

The scattering states $s\subset P$ form an open subset (see Thm\ 2.3.3
of \cite{DG}), and thus $U^-_{E,r}\cap s_E$ is relatively open, too.
All orbits in $\SuE$ which do not meet the interaction zone $\IZ(E)$
are scattering. Thus $N\setminus A_E^\pm(s_E)$ is a compact set. 
Then using $A_E^\pm(s_E^\pm)=N$ and asymptotic completeness 
Corollary \ref{coro:complete}.\ref{asymptotically:complete},
we see that $\lambda_N(N\setminus A_E^\pm(s_E))=0$.
\hfill $\Box$\\[2mm]
For $\lambda_{S^2}$-almost all $\htheta^{-} \in S^{2}$ the map
\beq
\vet:T^*_{\htheta^{-}}S^2\ar S^{2},\qquad 
\vet(\vec{a}^{-}) := \hat{p}^+(E,\vec{a}^{-},\htheta^{-})
\Leq{def:vet}
is measurable. Here the restriction of the asymptotic direction
$\hat{p}^+: s^+\ar S^2$ (which is constant on orbits) to $s_E$
is considered as a map 
\[\hat{p}^+(E,\cdot,\cdot):A_E^-(s_E)\ar S^2,\] 
see Prop.\ \ref{prop:MAA} 
\begin{definition} \label{defi:diff:cross}
For $E>0$ and $\htheta^{-}\in S^2$
the {\bf cross section measure} $\set$ on $S^{2}$ is the image measure
\beq
\set := \vet\l( \lambda_{\htheta^{-}} \ri),
\Leq{csm}
$\lambda_{\htheta^{-}}$ being Lebesgue measure on
the cotangent plane $T^*_{\htheta^{-}}S^2$.

Assuming $\set$ on $S^{2}\setminus \{\htheta^{-}\}$ to be
absolutely continuous w.r.t.\ Haar measure $\lambda_{S^2}$,
the {\bf differential cross section} $\DCSP$ is the Radon-Niko\-dym 
derivative of $\set$.
\end{definition}
In Def.\ (\ref{csm}) we by using $\lambda_{\htheta^{-}}$ we
normalize the flux through unit area in configuration space to equal 
one.
\begin{example} 
By radial symmetry, for the simplest case (\ref{free:hamiltonian})
of scattering by a Kepler potential with $Z\equiv\Zi\neq 0$, the 
(Rutherford) differential cross section depends only on the angle
\[\Delta\theta := \sphericalangle(\htheta^-,\htheta^+)\]
between the initial and final direction.
Using formula (\ref{ecce}) for the eccentricity 
$e=+1/\sin(\eh\Delta\theta)$
of the Kepler hyperbola, we see that the modulus $a$ of the impact parameters
$\vec{a}^\pm$ equals
\[a=\frac{|Z|}{{2E}}\cot(\eh\Delta\theta),\]
so that
\[\l|\frac{da}{d\Delta\theta}\ri| = 
\frac{|Z|}{4E\sin^2(\eh\Delta\theta)}.\]
We may assume that $\htheta^-=(0,0,1)$ so that the 3-component of
$\vec{a}^-=\p^-\times\vec{L}^-$ vanishes and
$\lambda_{\htheta^{-}}$
corresponds to integration with the two-form $dL_1^-\wedge dL_2^-$.
Introducing polar coordinates $(L,\vv^-)$ in the plane 
$T^*_{\htheta^{-}}S^2$, and expressing the volume element at $\htheta^+$
on $S^2$ in the form 
\[\sin(\Delta\theta)\cdot d\vv^-\wedge d\Delta\theta,\]
we obtain the familiar expression
\beq
\l(\DCSP\ri)_{\rm Ru} =\l|\frac{L}{2E}
\frac{dL}{d\Delta\theta}\ri|=\l(\frac{Z}{4E\sin^2(\eh\Delta\theta)}\ri)^2
\Leq{Ruth}
for the {\em Rutherford cross section}. 
Note that it depends only on the modulus of the charge $Z$.
\end{example}
For $(\htheta^-,\htheta^+)\in\thmin^{-1}(\theta)$, $\uk\in\cW$ let
\[\vec{a}_{\uk}^{-}\l(E,\htheta^{-},\htheta^{+}\ri):=
\lim_{t\ar-\infty}\vec{a}({\rm Diff}_E^{-1}(t,\htheta^-,\htheta^+,\uk)),\]
with ${\rm Diff}_E$ defined in Thm.\ \ref{thm:sc}.
\begin{theorem} \label{theorem:classify:scattering} 
Let $V$ be a Coulombic potential satisfying the decay estimates (\ref{smooth}).

Then for energy $E>\Eth$ and $\vartheta=\min(c E^{-\delta},\amin)$ with
$0\leq\delta\leq \eh$,
on $\thmin^{-1}(\vartheta)$ the differential cross section is smooth,
of the form 
\beq
\DCSP = 
\sum_{\uk \in \cW}
\l|\det\l( \frac{d\vet}{d\vec{a}} 
(\vec{a}_{\uk}^{-}(E,\htheta^{-},\htheta^{+}))\ri) \ri|^{-1},
\Leq{eq:dcs:sum}
and differs from the Rutherford cross section (\ref{Ruth})
for charge $Z :=\sqrt{{\textstyle\sum_{l=1}^n Z_l^2}}$ only by 
\beq
\DCSP =\l(\DCSP\ri)_{\rm Ru} \cdot(1+\cO(E^{2\delta-1}))
\Leq{comp:Ru}
uniformly in $(\htheta^-,\htheta^+)\in\thmin^{-1}(\vartheta)$.
\end{theorem}
\begin{remarks}
{\bf 1)} 
In particular the relative difference w.r.t.\ Rutherford cross section 
is only of order $\cO(1/E)$ if one excludes cones
of an energy-independent aperture $\vartheta$.\\
{\bf 2)}
In 2-dim.\ and purely Coulombic potentials the differential cross section
is even smooth 
(up to the forward direction $\htheta^+=\htheta^-$) for all positive
energies \cite{KK}.
As shown in \cite{Kn3}, this smoothness is rather exceptional. 

Also here for $n>1$ centres we cannot have smoothness of the differential 
cross section (for {\em no} $V$ and {\em no} $E>\Eth$),
if we add to $\AD(\vartheta)$ 
the neighbourhood of {\em any} direction $\hat{s}^{i,k}$.
This follows from the observation that (contrary to the 2D attracting case)
in 3 dimensions any hard collision with a nucleus
changes the degree. By a limit argument there must then exist points
where the degree is zero. At these points the differential 
cross section diverges. 
\end{remarks}
{\bf Proof.}
Thm.\ \ref{thm:sc} says that the orbits with data $(E,\htheta^{-},\htheta^{+})$
are enumerated by $\cW$. So if the r.h.s.\ of
(\ref{eq:dcs:sum}) converges, then by the Transformation Theorem for
Lebesgue measure (\ref{eq:dcs:sum}) follows from the definition of 
the cross section measure in (\ref{csm}).

Estimates (\ref{DO:p}) and (\ref{DO:L}) of Thm.\ \ref{thm:smooth:moeller}
imply that for
$|\q_0| \geq \Rvir(E)$, $\pm\LA\q_0,\p_0\RA \geq 0$ and
multi-indices $\gamma:=(\alpha,\beta)\in\bN_0^3\times\bN_0^3$
\[\pa^\gamma_{x_0}(\hat{p}^{\pm}(x_0)-\hat{p}_0)=
\cO\l(E^{-1-\eh|\alpha|}\ri)\  \]
and  
\[\pa^\gamma_{x_0}(\vec{a}^{\pm}(x_0)-\vec{a}(x_0))=
\cO\l(E^{-1-\eh |\alpha| })\ri),\]
the last estimate being obtained with the help of (\ref{eq:impact}).

Similar statements are true for orbit segments in $\IZ(c_q)$. 

So up
to an error of order $\cO(1/E)$, all variations of the asymptotic data come 
from the single scattering processes within the balls of radius $c_q$
around the singularities. 

We switch to $(\y, \z)$-coordinates.

There are two types of such contributions:
\begin{enumerate}
\item
The ones coming from the hard collisions 
($\theta$--visits in succession $\uk$).
These lead to factors of the form  (\ref{Tx}) in the product formula 
for the linearization:
\beqn
T_x\Po = f(k_{-1},k_0,k_1)E\cdot\bem{cc}
\idty&\idty\\\idty&\idty\eem+\cO(E^0)
\label{expansion}
\eeqn
with
\[f(k_{-1},k_0,k_1):=
\frac{2d^{k_{-1},k_0}\cos^2(\eh\alpha(k_{-1},k_0,k_1))}{-Z_{k_0}}.\]
As $\Delta\psi>c E^{-\delta}$, the relative error in this 
estimate is of the order $\cO(E^{2\delta-1})$.
\item
Visits of $c_q$-balls around some singularities, which are not 
hard collisions ($\theta$--visits).
By what we have shown, there can be at most two such events,
one before the $\uk$--visits, one after.

So these visits meet the hypothesis (\ref{large:dev}) of 
Lemma \ref{lem:large:dev}
\[\l|\frac{\vec{v}^-}{|\vec{v}^-|} + \w^-\ri|\geq c E^{-\delta}),\] 
{\em i.e.} with angle $\vartheta=\cO(E^{-\delta'})$, where 
$\delta':=1-\delta\in[\eh,1]$.

So formula (\ref{dev:from:free}) of Lemma \ref{lem:lin} says that the 
relative deviation 
of the linearized flow from free motion during these soft collisions 
is of order $\cO(E^{1-2\delta'})=\cO(E^{2\delta-1})$:
\[\delta\vec{v}^+ = \delta\vec{v}^- + \hspace*{34mm}
\cO(c^2E^{1-2\delta'})\cdot (|\delta\vec{v}^-|+|\delta\w^-|),\]
\[\delta\vec{w}^+ = 2u \delta\vec{v}^- + \delta\vec{w}^- + 2\vec{v}^- du +
\cO(c^2E^{1-2\delta'})\cdot (|\delta\vec{v}^-|+|\delta\w^-|).\]
\end{enumerate}
Now we see that the formal sum (\ref{eq:dcs:sum}) converges:
\begin{itemize}
\item
There are exactly $n\cdot (n-1)^{l-1}$ words $\uk\in\cW$ of length $l$.
\item
For word length $l+1$ of $\uk\in\cW$ the term 
\[ \l|\det\l( \frac{d\vet}{d\vec{a}} 
(\vec{a}_{\uk}^{-}(E,\htheta^{-},\htheta^{+}))\ri) \ri|\]
in (\ref{eq:dcs:sum})
is only of relative order $\cO(E^{-(d-1)})=\cO(1/E^2)$, compared to the term of 
the word shortened by one letter, as there is one extra factor $E$
coming from (\ref{expansion}). Here $d$ denotes the dimension, so $d=3$.
\end{itemize}

So if $E>\Eth$ and the threshold $\Eth$ is suitably chosen, this decay 
outweighs the exponential proliferation of words with given word length.

The comparison in (\ref{comp:Ru})
with Rutherford cross section (\ref{Ruth})
for squared charge $Z^2 ={\textstyle\sum_{l=1}^n Z_l^2}$ follows by
adding the contributions in (\ref{eq:dcs:sum}) of word length one.
In $d=3$ dimensions the leading errors of order $\cO(E^{2\delta-1})$ 
come from  1) and 2) above, whereas neglecting the contributions of the 
longer words is only of order $\cO(E^{-2})$.
\hfill $\Box$\\[2mm]
In $d$ dimensions the Rutherford cross section equals
\beq
\l(\DCSP\ri)_{\rm Ru} =\l(\frac{|Z|}{4E\sin^2(\eh\Delta\theta)}\ri)^{d-1}.
\Leq{Rud}
In 12.4 of \cite{KK} we remarked that for $d=2$ dimensions, $Z_l>0$
and $\theta^+\neq\theta^-$ 
the differential cross section of the $n$--centre problem converges to the $d=2$
Rutherford cross section as $E\ar\infty$. The charge $Z$ in (\ref{Rud}), 
however, must be chosen
as $Z=\pm\sum_{i=1}^n Z_i$ (and not $\Zi$, as wrongly stated in \cite{KK}). 
This result can be sharpened:
\begin{corollary}
The analog of formula (\ref{comp:Ru}) holds true in $d=2$ dimensions
if one sets $Z:=\sum_{i=1}^n |Z_i|$.
\end{corollary}
{\bf Proof.}
Up to error terms the formulae (\ref{Tx}) and (\ref{dev:from:free}) 
for the linearization of the
flow which we used to derive the result (\ref{comp:Ru}) in $d=3$ are invariant 
under rotations. Thus (\ref{comp:Ru}) is true in $d=2$ dimensions, too.
\hfill$\Box$
%
\Section{The Collinear Case} \label{sect:collinear}
%
In this section we show by counterexample that the non-collinearity conditions
cannot be dropped in Thm.\ \ref{thm:homeo} and Thm.\ \ref{thm:sc}.\\[2mm]
{\bf 1)} We first consider the set $\buE$ of bounded orbits of energy $E$
for {\em attracting} Coulombic potentials $V$ which are rotationally 
symmetric w.r.t.\ some axis $A\subset
\Mu$. Thus in particular the nuclei are situated on that axis:
\beq
\s_1,\ldots,\s_n\in A.
\Leq{on:axis}
Conversely, that condition ensures that $V$ is rotationally symmetric around
$A$, if it is a purely  Coulombic potential.

By symmetry trajectories with initial conditions tangential to
a two-plane $F\subset\Mu$ containing $A$ stay in $F$, and we may thus  
consider the restricted two-dimensional motion on $F$.

The axis $A$ is divided by $\s_1,\ldots,\s_n$ into $n+1$ closed
intervals meeting only in their endpoints $\s_l$. These intervals 
correspond for energies $E>\Vmax$ 
to trajectories which are reflected by the nuclei. Thus two of these
trajectories are unbounded and $n-1$ bounded.

These special trajectories are of course invariant under rotations
around $A$. However, there cannot be any further energy $E$ trajectories in $F$
having this property.

On the other hand it has been shown in \cite{KK} that for $n\geq 3$
nuclei above some energy
$\Eth$ there is a Cantor set of bounded trajectories in $F$.
So in particular there is an uncountable number of trajectories in $F$
which are not moving tangential to $A$ and thus give rise to {\em
one-parameter}
families of trajectories for the full motion in $\Muh$.\\[2mm]
{\bf 2)} 
In general there are bounded orbits which do not lie in any plane $F$
containing the axis $A$. We observe that by rotational symmetry the component
$\LA\vec{L},\hat{s}\RA$ of angular momentum in the direction $\hat{s}$
of $A$ is preserved. We now indicate that 
for certain collinear configurations there exist {\em two-parameter} families of
bounded orbits of a given energy, parametrized by that angular momentum 
component and its conjugate angle.

For the sake of simplicity we consider a 3-centre potential 
\[V(\q):= -\sum_{l=1}^3 \frac{Z_l}{|\q-\s_l|}\]
with $s_1:=\vec{0}$, $\s_2:=d\cdot\hat{s}$, $\s_3:=-d\cdot\hat{s}$, 
$\hat{s}:=(0,0,1)$ and $Z_2=Z_3$ which, in
addition of being axially symmetric w.r.t.\ the 3-axis $A=\bR\cdot\hat{s}$,
is mirror-symmetric w.r.t.\ reflection by the plane 
$F_{12}:=\{ \q\in\Mu \mid q_3=0\}$. 
We first consider periodic 
trajectories with angular momentum component $L_3=0$ in, 
say, the plane $F\cong F_{13} :=\{ \q\in\Mu \mid q_2=0\}$
which are invariant under reflection by the plane $F_{12}$:
\[q_1(-t)=q_1(t)\qmbox{,}q_2(t)=0\qmbox{,}q_3(-t)=-q_3(t)
\qquad(t\in\bR).\]
By symbolic dynamics arguments combined with Thm.\ 6.11 and Remark 11.2.2
of \cite{KK} for all energies $E>0$ there exists a countable
infinity of these orbits, all being hyperbolic and having index 0 (as orbits
in the two-plane $F_{13}$. 
Thus by invoking an implicit function
argument one may show the existence of smooth family of energy-$E$
bounded orbits starting on the 1-axis near $\q(0)$ and parametrized by
$L_3$. Rotating these orbits around the axis $A$ then yields a two-parameter
family of bounded orbits.\\[2mm]
{\bf 3)} 
For {\em repelling} axially symmetric potentials the situation is completely
different. W.l.o.g.\ we again consider potentials $V$ which are invariant
w.r.t\ rotations around the axis $A=\bR\cdot\hat{s}$, with $\hat{s}=(0,0,1)$.
But now we assume that
\[\LA\nabla V(\q),\q-\langle\q,\hat{s}\rangle \hat{s}\RA < 0 \qquad
(\q\in\Mu\setminus A).\] 
This condition is met, e.g., 
by repelling ($Z_l<0$) purely Coulombic potentials meeting (\ref{on:axis}).

Now consider a trajectory starting at $(\p(0),\q(0))$ with 
\beq\q(0)-\langle\q(0),\hat{s}\rangle \hat{s}\neq \vec{0}\qmbox{and} 
\LA\nabla \p(0),\q(0)-\langle\q(0),\hat{s}\rangle \hat{s}\RA\geq 0.
\Leq{not:on:axis}
Then 
\[\LA\p(t),\q(t)-\langle\q(t),\hat{s}\rangle \hat{s}\RA > 0 \qquad (t>0)\]
and is monotonically inreasing in $t$, since 
\[\frac{d}{dt}\LA\p,\q-\langle\q,\hat{s}\rangle \hat{s}\RA 
= -\LA\nabla V(\q),\q-\langle\q,\hat{s}\rangle \hat{s}\RA+
\LA\p,\p-\langle\p,\hat{s}\rangle \hat{s}\RA>0.\]
But this means that the orbit leaves the interaction zone in finite time and
thus is not bounded.
The second of the conditions in (\ref{not:on:axis}) is not restrictive, since
otherwise we may consider negative times. 

We conclude that the only bounded orbits lie on the axis. Thus for $E>\Vmax$
there are exactly $n-1$ bounded orbits, compared to the
uncountable infinity of bounded orbits 
in the case of NC configurations and $n\geq3$.\\[2mm]
{\bf 4)} 
Finally we consider {\em scattering} orbits. Already for $n=2$ nuclei and
attracting
purely Coulombic potentials we have {\em one-parameter} families of 
orbits of a given energy $E$ scattering from a direction
parallel to the axis $A$ through the positions $\s_1,\s_2$
and to the backward direction. These are obtained
by rotating a given solution around $A$.
There are infinitely many such families, as can be seen from the explicit 
Jacobi solution of the two-centre problem, described in Appendix B, or from
\cite{KK}, Thm.\ 12.1.
\appendix
%
\Section{Aspects of Geometry and Global Analysis}
%
This article on the 3-dimensional $n$--centre problem
is based upon analytical perturbation estimates.
To the contrary the 2-dim.\ $n$--centre problem 
(and similarly, the 2-dim.\ periodic potentials of \cite{Kn1,Kn2})
was treated in \cite{KK}
using techniques from Riemannian geometry and global analysis.

In this appendix both approaches are compared.

It is known that  
the trajectories of energy $E>\sup_\q V(\q)$
generated by a Hamiltonian function $H:T^*M\ar\bR$ of the (local) form 
\[H(\p,\q) = \eh \sum_{k,l=1}^d g^{k,l}(\q)p_k p_l +V(\q)\]
on a $d$-dimensional Riemannian manifold $(M,g)$ coincide
(up to a time re\-para\-metrization) with the geodesics in the so-called 
{\em Jacobi metric} $g_E$ on $M$, conformally equivalent to $g$:
\beq
g_E(\q):= (1-V(\q)/E) g(\q)
\Leq{def:jacobi}
(here we assume for simplicity $E> 0$).

In the simplest case covered by the paper $V(\q)=-Z/|\q|$ with $Z>0$, 
{\em i.e.} the attracting
Coulomb potential. There, using a formula from Spivak (\cite{Sp}, p.\ 337), 
we obtain the expression
\beq
K_{1,2}(\q)= \frac{Z}{2E}\frac{-1+3\l(1+\frac{Z}{2E|q|}\ri)\frac{q_3^2}
{|\q|^2}}{\l(|\q|+Z/E \ri)^3}
\Leq{sect:curv} 
for the sectional curvature of the Jacobi metric in the 1-2-tangent plane at
$\q$.
We can learn several things from that formula:
\begin{itemize}
\item
Setting $q_3=0$ and thus considering planar motion, 
\[K_{1,2}(\q)=-\frac{Z}{2E(|\q|+Z/E)^3}<0\]
for positive $E$, and this expression is bounded below, by
\beq
K_{1,2}(\q)\geq -\frac{E^2}{2Z^2}.
\Leq{b:bel}
This fact was used extensively in \cite{Kn1}, \cite{Kn2} and \cite{KK} 
to analyze planar motion by going to the {\em smooth}
branched covering surface 
\[\Mo := \l\{ (q,Q) \in \bC\times\bC 
                        \l| Q^{2} = \prod_{l=1}^{n}(q-s_{l})\ri. \ri\} \]
of the configuration plane, equipped with the lifted
Jacobi metric (the branched covering being given by projection to the first
factor $q$).
\item
For $d=3$ dimensions the sectional curvature (\ref{sect:curv}) is neither
uniformly bounded in the $\q$ variable nor definite.
In fact, for $\q=(0,0,q_3)$
\beq K_{1,2}(\q)\sim \frac{3E}{4Z|\q|}\ar +\infty\qquad (q\ar 0).
\Leq{three:d}
We thus consider here geodesic motion in mixed sectional curvature.
Although the $E$-dependence of (\ref{b:bel}) is quadratic, whereas
(\ref{three:d}) is only linear in $E$, negative curvature 
does never dominate positive curvature in our estimates. Namely we
have seen in Lemma \ref{lem:defl} that
the minimal distances of the bounded orbits from the nuclei are of the 
order $1/E$, so that effectively (\ref{three:d}), too, goes like $E^2$.

Mixed curvature dynamics is rather intractable in general.
However, in the case considered here the motion near the singularities 
can be treated as a perturbation of Keplerian motion, and this allows us to
control the motion in the high energy limit.
\end{itemize}
Whereas the estimates of this papers are somewhat optimal in the high
energy limit $E>\Eth$, nothing much could be said about 
the energy region $0<E\leq\Eth$ (the exception being Sect.\ 
\ref{sect:moeller}).

To the contrary, for two dimensions 
many results were shown for all positive energies, using 
the negativity of Gaussian curvature and the topology of the branched
covering surface $\Mo$ (whose fundamental group is non-abelian for $n\geq3$).

The branched covering $\Mo\ar\bC$ globalizes the so-called Levi-Civita
transform $Q\mapsto Q^2$ of celestial mechanics.
So it is natural to pose the question whether there exists a 
similar globalization of the Hopf map
\[\bC^2\ar\bR^3,\qquad z\mapsto 
\LA z,\vec{\sigma}z\RA\] 
used in the Kustaanheimo-Stiefel regularization of a 3-dim.\
Coulomb singularity. As already mentioned, this was done in \cite{HS} 
for $n=2$ centres. However, a generalization to arbitrary $n$
seems to be unknown.

We expect that the corresponding manifolds, {\em i.e.} four-dimensional
analogs of Riemann surfaces, should have interesting topological 
properties.

One last aspect concerns structural stability. Both in $d=2$
and three dimensions the compact set $\buE\subset\SuE$ of bounded orbits
is hyperbolic and thus structurally stable.

Thus if we continuously move the singularities by 
suitable maps $S^1\ar \bR^d$, $u\mapsto \s_l(u)$, we obtain a family of 
Coulombic potentials $V_u$ and correspondingly a one-parameter
family of bound states $b_{E,u}$. As the parameter varies cyclically, we 
obtain a permutation of the bounded orbits in $\buE\equiv b_{E,1}$.

For $d=2$ this action of the {\em braid group on $n$ strands} of $\bR^2$ by 
permutations is nontrivial in general (see Remark 6.12 of \cite{KK}).

Although in $d=3$ dimensions the manifold 
$\NC$ of non-collinear configurations is not simply connected for $n\geq 3$,
the action of the fundamental group $\pi_1(\NC)$ on $\buE$
is trivial.

\Section{The Two-Centre Problem}
%
Here we shortly discuss the purely Coulombic two-centre problem,
{\em i.e.}
\[V(\q)=\frac{-Z_1}{|\q-\s_1|}+\frac{-Z_2}{|\q-\s_2|}.\]
W.l.o.g.\ we assume that the two centres are at $\s_1:=\bsm 1\\0\\0\esm$
and $\s_2:=\bsm -1\\0\\0\esm$. 
As is well-known (see e.g.\ Thirring \cite{Th}), 
the problem is analytically integrable,
see \cite{GKM} for an application to satellite
motion and \cite{SR} for an application to semiclassics of the hydrogen 
molecule.

The motion is integrated using the prolate ellipsoidal coordinates 
$(\xi,\eta,\vv)\in\bR^+\times [0,\pi)\times[0,2\pi)$
with
\[\q\equiv \bsm q_1\\ q_2\\ q_3 \esm=
\bsm \cosh(\xi)\cos(\eta)\\  \sinh(\xi)\sin(\eta)\cos(\vv)\\
                             \sinh(\xi)\sin(\eta)\sin(\vv) \esm.\]
As in these coordinates $H$ is independent of $\vv$, the conjugate momentum
\[p_\vv= q_2p_3-q_3p_2\] 
is a constant of the motion, equal to the first component of angular momentum.
We set $l_1:=p_\vv(x_0)$ for initial conditions $x_0$.

Then by going to extended phase space and using a new time parameter
$s$ defined by 
\[\frac{dt}{ds} = 2\l(\cosh^2(\xi)-\cos^2(\eta)\ri),\]
the new Hamiltonian function separates:
\[\cH:=\frac{dt}{ds}(H-E) = H_1+H_2\] 
with 
\[H_1(p_\xi,\xi):= p_\xi^2+V_1(\xi)\qmbox{with}
V_1(\xi):= \frac{l_1^2}{\sinh^2(\xi)}-2Z_+\cosh(\xi)-2E\cosh^2(\xi)\]
\[H_2(p_\eta,\eta):= p_\eta^2+V_2(\eta)\qmbox{with}
V_2(\eta):=\frac{l_1^2}{\sin^2(\eta)}+2Z_-\cos(\eta)+2E\cos^2(\eta), \]
where $Z_\pm:=Z_2\pm Z_1$.

The motion on $\cH^{-1}(0)$ coincides --- up to time parameterization ---
with the motion on $H^{-1}(E)$. 
Setting 
\[K:=H_1(x_0)=-H_2(x_0),\] 
we have three generally independent 
constants of the motion $H,H_1$ and $l_1$, whose values are denoted by 
$E,K$ and $l_1$, respectively. 

The {\em bifurcation set} is then given by the set of values for which the 
mapping from phase space to the constants of the motion is not 
{\em locally trivial} (see \cite{AM}, Sect.\ 4.5). The most interesting
subset is the one for $l_1=0$, {\em i.e.} two-dimensional motion.

\begin{figure}[ht]
\centerline{
\epsfig{file=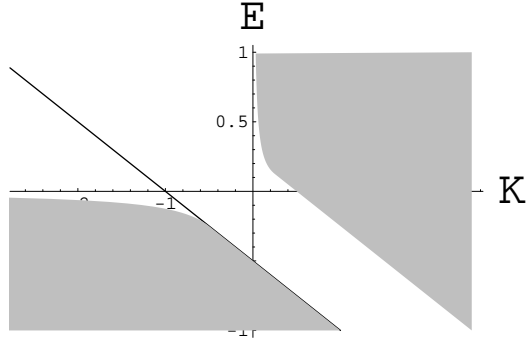,height=5cm,
bbllx=90,bblly=510,bburx=380,bbury=710,clip=}%
}
\caption{Bifurcation diagramme}\label{fig2}
\end{figure}

By inspection of the extrema of the $V_i$ one sees that
for $l_1=0$ the image of $(H, H_1)$ is the region in $\bR^2$ bounded by the 
curves
\[K_+(E) := \l\{\begin{array}{ll}
\frac{Z_+^2}{2E} & ,0>E>-Z_+/2\\
-2(Z_+ + E)      & ,E\leq -Z_+/2\end{array}\ri.\]
and
\[K_-(E) := \l\{\begin{array}{ll}
\frac{Z_-^2}{2E} & ,E>|Z_-|/2\\
2(|Z_-|-E)       & ,E\leq |Z_-|/2\end{array}
\ri. .\] 
The bifurcation diagramme (see Fig.\ \ref{fig2}) is the union of $K_-$, $K_+$, and the lines
\[E=0\qmbox{,}K=0\qmbox{and} K_0(E):= -2(Z_+ + E)\] 
inside the image of $(H, H_1)$. 

The line $K_0$ corresponds to the $(H,H_1)$--values 
of the closed orbit wandering
between the centres $\s_1$ and $\s_2$, {\em i.e.} having coordinate $\xi=0$.

The following relations are useful for the scattering problem.
In the (1-2)--plane, {\em i.e.} for $\vv=0$, we have the polar coordinates
$q_1=r\cos\phi$, $q_2=r\sin\phi$. Then $r^2=\cosh^2(\xi)-\sin^2(\eta)$ 
and $\tan(\phi)=\tanh(\xi)\tan(\eta)$ so that in the $r\ar\infty$  limit
$\eta$ coincides with $\phi$.

In the same limit $p_\eta=q_1p_2-q_2p_1-e^{-\xi}(p_1\sin\eta+p_2\cos(\eta))$ 
coincides with the angular momentum $q_1p_2-q_2p_1$ of the (1-2)--plane.

This suffices to relate the asymptotic data 
$(p_\eta^\pm,\eta^\pm):=\lim_{s\ar\s^\pm}(p_\eta(s),\eta(s))$
with the ones in the original system (the times $s^+>s^-$ being defined by
$\lim_{s\ar\s^\pm}\xi(s)=\infty$). The constant $K$ is then given by
$(p_\eta^\pm)^2+V_2(\eta^\pm)$.
Equalling the elliptic integrals 
\[s^+-s^-= \l|\int_{\eta^-}^{\eta^+}\frac{d\eta}{\sqrt{-K-V_2(\eta)}}\ri|\]
respectively
\[s^+-s^-= 2\int_{\xi_{\rm min}}^{\infty}\frac{d\xi}{\sqrt{K-V_1(\xi)}}
\qmbox{with}\cosh(\xi_{\rm min}) = 
\frac{-Z_+}{2E}+\sqrt{\l(\frac{Z_+}{2E}\ri)^2-\frac{K}{2E}}\]
then suffice in principle to calculate analytically the 
scattering transformation, but the expressions become rather lengthy.
\addcontentsline{toc}{section}{References}
{\bf Note added in proof.}

The anonymous referee informed me about the interesting related article:\\ 
S.V. Bolotin, P. Negrini: 
Regularization and topological  entropy for the spatial $n$-center problem,
which meanwhile appeared in 
{\em Ergodic Theory and Dynamical Systems} {\bf 21}, 383--399 (2001).\\
Concerning the context of the $n$-centre-problem on $\bR^3$, the
authors succeed to construct a global regularization of $n$ attracting
singularities, based on the local KS transform. Furthermore, they prove
that for $n\geq3$ the topological entropy is strictly positive for
{\em all} energies $E\geq0$, whereas the present paper is only dealing with all
energies above a positive threshold energy. 
Additionally, the authors prove several results
for configuration manifolds different from $\bR^3$.

I do not believe, however, that these topological methods could be used
to substantially simplify the proofs of the analytic results 
given in the present paper.
\end{document}